\documentclass[a4paper, 10pt]{article}
\usepackage{cancel, ulem}

    \makeatletter
      
\usepackage{fancyhdr}
\usepackage{bm}
\usepackage{cancel}
\usepackage{latexsym}
\usepackage{tikz}
\usepackage{color}
    
\usepackage{amssymb,amsmath,amsfonts,bbm,pifont,upgreek,bbold
}  
\usepackage[colorlinks=true]{hyperref}
\hypersetup{urlcolor=blue, citecolor=red}

\usepackage{graphicx}
\definecolor{Salmon}{cmyk}{0,0.53,0.38,0}
\definecolor{MidnightBlue}{cmyk}{0.98,0.13,0,0.43}
\definecolor{BurntOrange}{cmyk}{0,0.51,1,0}
\definecolor{OliveGreen}{cmyk}{0.64,0,0.95,0.40}
\definecolor{Emerald}{cmyk}{1,0,0.51,0}
\definecolor{airforceblue}{rgb}{0.36, 0.54, 0.66}
\definecolor{bleudefrance}{rgb}{0.19, 0.55, 0.91}
\definecolor{celestialblue}{rgb}{0.29, 0.59, 0.82}
\definecolor{cerulean}{rgb}{0.0, 0.48, 0.65}
\definecolor{cornflowerblue}{rgb}{0.39, 0.58, 0.93}
\definecolor{cyan(process)}{rgb}{0.0, 0.72, 0.92}
\definecolor{ao}{rgb}{0.0, 0.0, 1.0}
	\definecolor{red}{rgb}{1.0, 0.0, 0.0}
	\definecolor{amethyst}{rgb}{0.6, 0.4, 0.8}
\newcommand\tred[1]{\textcolor{red}{#1}}

\newcommand\tblue[1]{\textcolor{ao}{#1}}

\newcommand\beq[1]{\begin{equation}\label{#1} }
\newcommand{\eeq}{\end{equation} }
\newcommand{\beqno}{\[ }
\newcommand{\eeqno}{\] }
\newcommand\beqa[1]{\begin{eqnarray} \label{#1}}
\newcommand{\eeqa}{\end{eqnarray} }
\newcommand{\beqano}{\begin{eqnarray*} }
\newcommand{\eeqano}{\end{eqnarray*} }
\newcommand\arr[1]{\left\{\begin{array}{l}#1\end{array}\right.}

\renewcommand{\theequation}{\arabic{section}.\arabic{equation}}

\newtheorem{theorem}{Theorem}[section]
\newtheorem{definition}{Definition}[section]
\newtheorem{proposition}{Proposition}[section]
\newtheorem{lemma}{Lemma}[section]
\newtheorem{conjecture}{Conjecture}[section]

\newtheorem{sublemma}{Sublemma}[section]
\newtheorem{remark}{Remark}[section]
\newtheorem{notationalremark}{Notational Remark}[section]
\newtheorem{corollary}{Corollary}[section]
\newtheorem{assumption}{Assumption}[section]
\newtheorem{claim}{Claim}[section]

\newtheorem{tools}{$\negsp\negsp$}[subsection]

\newcommand\thm[1]{\begin{theorem}\label{#1}}
\newcommand\thmtwo[2]{\begin{theorem}[#1]\label{#2}}
\newcommand\ethm{\end{theorem} }
\newcommand\dfn[1]{\begin{definition}\label{#1} \rm}
\newcommand\dfntwo[2]{\begin{definition}[#1]\label{#2} \rm}
\newcommand\edfn{\end{definition} }
\newcommand\pro[1]{\begin{proposition}\label{#1}}
\newcommand\protwo[2]{\begin{proposition}[#1]\label{#2}}
\newcommand\epro{\end{proposition} }
\newcommand\lem[1]{\begin{lemma}\label{#1}}
\newcommand\lemtwo[2]{\begin{lemma}[#1]\label{#2}}
\newcommand\elem{\end{lemma} }
\newcommand\sublem[1]{\begin{sublemma}\label{#1}}
\newcommand\sublemtwo[2]{\begin{sublemma}[#1]\label{#2}}
\newcommand\esublem{\end{sublemma} }
\newcommand\rem[1]{\begin{remark}\label{#1} \rm}
\newcommand\erem{\end{remark} }
\newcommand\notrem[1]{\begin{notationalremark}\label{#1} \rm}
\newcommand\enotrem{\end{notationalremark} }
\newcommand\cor[1]{\begin{corollary}\label{#1}}
\newcommand\cortwo[2]{\begin{corollary}[#1]\label{#2}}
\newcommand\ecor{\end{corollary} }
\newcommand\asmp[1]{\begin{assumption}\label{#1}}
\newcommand\asmptwo[2]{\begin{assumption}[#1]\label{#2}}
\newcommand\easmp{\end{assumption} }
\newcommand\clm[1]{\begin{claim}\label{#1}}
\newcommand\eclm{\end{claim} }
\newcommand{\proof}{\par\medskip\noindent{\bf Proof\ }}
\newcommand\equ[1]{{\rm (\ref{#1})}}

\expandafter\chardef\csname pre amssym.def
at\endcsname=\the\catcode`\@
\catcode`\@=11
\def\undefine#1{\let#1\undefined}
\def\newsymbol#1#2#3#4#5{\let\next@\relax
 \ifnum#2=\@ne\let\next@\msafam@\else
 \ifnum#2=\tw@\let\next@\msbfam@\fi\fi
 \mathchardef#1="#3\next@#4#5}
\def\mathhexbox@#1#2#3{\relax
 \ifmmode\mathpalette{}{\m@th\mathchar"#1#2#3}%
 \else\leavevmode\hbox{$\m@th\mathchar"#1#2#3$}\fi}
\def\hexnumber@#1{\ifcase#1 0\or 1\or 2\or 3\or 4\or 5\or 6\or 7\or
8\or
 9\or A\or B\or C\or D\or E\or F\fi}
\ifcase\@ptsize
 \font\tenmsb=msbm10
 \font\sevenmsb=msbm7
 \font\fivemsb=msbm5
\or
 \font\tenmsb=msbm10 scaled \magstephalf
 \font\sevenmsb=msbm7 scaled \magstephalf
 \font\fivemsb=msbm5  scaled \magstephalf
\or
 \font\tenmsb=msbm10 scaled \magstep1
 \font\sevenmsb=msbm7 scaled \magstep1
 \font\fivemsb=msbm5 scaled \magstep1
\fi
\newfam\msbfam
\textfont\msbfam=\tenmsb
\scriptfont\msbfam=\sevenmsb
\scriptscriptfont\msbfam=\fivemsb
\edef\msbfam@{\hexnumber@\msbfam}
\def\Bbb#1{\fam\msbfam\relax#1}
\def\widehat#1{\setboxz@h{$\m@th#1$}%
 \ifdim\wdz@>\tw@ em\mathaccent"0\msbfam@5B{#1}%
 \else\mathaccent"0362{#1}\fi}
\def\widetilde#1{\setboxz@h{$\m@th#1$}%
 \ifdim\wdz@>\tw@ em\mathaccent"0\msbfam@5D{#1}%
 \else\mathaccent"0365{#1}\fi}

\def\RIfM@{\relax\ifmmode}
\def\nonmatherr@#1{\errmessage{\string#1\space allowed only in math mode}}
\def\Bbb{\RIfM@\expandafter\Bbb@\else
 \expandafter\nonmatherr@\expandafter\Bbb\fi}
\def\Bbb@#1{{\Bbb@@{#1}}}
\def\Bbb@@#1{\fam\msbfam\relax#1}
\def\setboxz@h{\setbox\z@\hbox}
\def\wdz@{\wd\z@}
\catcode`\@=\csname pre amssym.def at\endcsname
\newcommand{\ii}{{\rm i}  }
\newcommand{\ie}{{\rm i.e.\,}}
\newcommand{\eg}{{\rm e.g.\,}}

\newcommand{\Giu}{{\bigskip\noindent}}
\newcommand{\nl}{{\smallskip\noindent}}
\newcommand{\noi}{{\noindent}}

\newcommand{\qed}{\hskip.5truecm
\vrule width 1.7truemm height 3.5truemm depth 0.truemm
\par\Giu}

\newcommand{\negsp}{\hspace{-.09truecm}}  

\newcommand{\dst}{\displaystyle}
\newcommand\ovl[1]{\overline {#1} }

\newcommand\su[1]{\frac{1}{{#1}} }

\newcommand{\tr}{{\, {\rm tr}\, }}

\newcommand{\torus}{{\Bbb T}   }
\renewcommand{\natural}{{\Bbb N}   }
\newcommand{\real}{{\Bbb R}   }
\newcommand{\integer}{{\Bbb Z}   }
\newcommand{\complex}{{\Bbb C}   }

\renewcommand{\a }{{\alpha}   }
\renewcommand{\b}{{\beta}   }
\newcommand{\g}{{\gamma}   }
\newcommand{\G}{{\Gamma}   }

\newcommand{\D}{{\Delta}   }
\newcommand{\e }{{\varepsilon}   }

\renewcommand{\k}{{\kappa}   }
\renewcommand{\l}{{\lambda}   }
\renewcommand{\L}{{\Lambda}   }
\newcommand{\m}{{\mu}   }
\newcommand{\n}{{\nu}   }

\newcommand{\p}{{\pi}   }
\renewcommand{\P}{{\Pi}   }
\renewcommand{\r}{{\rho}   }
\newcommand{\s}{{\sigma}   }

\renewcommand{\t}{{\tau}   }
\newcommand{\f}{{\varphi}   }

\renewcommand{\o}{{\omega}   }
\renewcommand{\O}{{\Omega}   }

\renewcommand{\Im}{{\, \rm Im\, }}
\renewcommand{\Re}{{\, \rm Re\, }}

\newcommand{\cA}{{\cal A} }
\newcommand{\cB}{{\cal B} }
\newcommand{{\cE}}{{\cal  E} }
\newcommand{\cT}{{\cal T} }
\newcommand{\cR}{{\cal R} }
\newcommand{{\cH}}{{\cal H} }
\newcommand{{\cK}}{{\cal K} }
\newcommand{\cC}{{\cal C} }
\newcommand{\cD}{{\cal D} }

\newcommand{\cG}{{\cal G} }
\newcommand{{\cJ}}{{\cal J}}
\newcommand{\cL}{{\cal L} }
\newcommand{\cM}{{\cal M} }
\newcommand{\cP}{{\cal P} }
\newcommand{\cI}{{\cal I} }

\newcommand{\cS}{{\cal S} }
\newcommand{\cU}{{\cal U} }

\newcommand\bx{{\mathbf x}}

\newcommand\by{{\mathbf y}}

\newcommand\ppu{{(1) }}
\newcommand\ppd{{(2) }}

\newcommand\pph{{(h) }}

\newcommand\ppi{{(i) }}

\newcommand\ppo{{(0) }}

\newcommand\ul{{\uplambda}}
\newcommand\ux{{\upxi}}
\newcommand\uh{{\upeta}}
\newcommand\up{{\rm p}}
\newcommand\uq{{\rm q}}
\newcommand\uz{{\rm z}}

\newcommand\DD{{\rm D}}
\newcommand{\HH}{\cH}

\newcommand\ZZ{{\rm Z}}

\usepackage{color}
\definecolor{yellow}{rgb}{0.99, 0.93, 0.0}

\newcommand\meas{{\, \rm meas\,}}

\newcommand\id{{\, \rm id \,}}

\newcommand\td{{n}}
\newcommand\tk{{k }}
\newcommand\pertnorm{{E}}
\newcommand\KAM{{\hat E}}

\begin{document}
\title{
\bf Perturbation theory  and canonical coordinates\\ in celestial mechanics\thanks{
Notes of two courses given, respectively, at the 18th School of Interaction Between Dynamical Systems and Partial Differential Equations  (Barcelona, June 27--July 1, 2022) and at the XLVII Summer School on Mathematical Physics (Ravello, 12--24 September, 2022).
I warmly thank the Centre de Recerca Matematica of Bellaterra (Barcelona) and Istituto Nazionale di Alta Matematica and the Gruppo Nazionale per la Fisica Matematica for
their  kind hospitality and especially A. Delshams, M. Guardia, T. Ruggeri, G. Saccomandi and T. M--Seara for their interest.

Sections \ref{sec: K map}, \ref{The reduction of perihelia}, \ref{P-map vs rotations and reflections}, \ref{Global Kolmogorov tori} and \ref {Coexistence of stable and whiskered tori} are based on work done while the author was funded by the ERC grant 677793 StableChaoticPlanetM (2016--2022). 

{\bf MSC2000 numbers:}
primary:
34C20, 70F10,  37J10, 37J15, 37J40;
secondary: 
34D10,  70F07, 70F15, 37J25, 37J35.
}
}
\author{Gabriella Pinzari\footnote{Department of Mathematics, University of Padua, e--mail address: {\tt pinzari@math.unipd.it}}}
\date{September 16, 2022}
\maketitle

\begin{abstract}\footnotesize{KAM theory owes most of its success to its initial motivation: the application to problems of celestial mechanics. The masterly application was offered by V.I.Arnold in the 60s who worked out a theorem, that he named the “Fundamental Theorem” (FT), especially designed for the planetary problem.  However, FT could be really used at that purpose only when, about 50 years later, a set of coordinates constructively taking the invariance by rotation and close--to--integrability into account was used. Since then, some progress has been done in the symplectic assessment of the problem, and here we review such results. 
}
\end{abstract}

\maketitle

\tableofcontents


\renewcommand{\theequation}{\arabic{equation}}
\setcounter{equation}{0}

\section{Some sets of canonical coordinates  for many--body problems}
\subsection{$(1+n)$--body problem, Delaunay--Poincar\'e coordinates and Arnold's theorem}\label{sec: AT intro}
In the masterpiece \cite{arnold63}, a young a brilliant mathematician, named Vladimir Igorevich Arnold, stated, and partly proved, the following result.
\begin{theorem}{\bf ``Theorem of stability of planetary motions'', \cite[Chapter III, p. 125]{arnold63}}\label{Arnold Theorem}
For the majority of initial conditions under which the instantaneous orbits of the planets are close to circles lying in a single plane, perturbation of the planets on one another produces, in the course of an infinite interval of time, little change on these orbits provided the masses of the planets are sufficiently small. {\rm[...]}
In particular {\rm[...]} in the n-body problem there exists a set of initial conditions having a positive Lebesgue measure and such that, if the initial positions and velocities of the bodies belong to this set, the distances of the bodies from each other will remain
perpetually bounded. 
\end{theorem}

\vskip.1in
\noindent
Let us summarize the main ideas behind the statement above. \\
After the symplectic reduction of the linear momentum, the $(1+n)$--body problem with masses $m_0$, $m_1$, $\ldots$, $m_n$ is governed by the $3n$--degrees of freedom Hamiltonian (see Appendix \ref{appendix})
\beqa{Helio}{\cal H}&=\,&\sum_{1\leq i\leq n}\left(\frac{|\by_i|^2}{2\mu_i}-\frac{\mu_i M_i}{|\bx_i|}\right)+\sum_{1\leq i<j\leq n}\left(\frac{\by_i\cdot \by_j}{ m_0}-\frac{m_i m_j}{|\bx_i-\bx_j|}\right)
\eeqa 
where $\bx_i$ represent the difference between the position of the $i^{\rm th}$ planet and the mass $m_0$, $\by_i$ are the associated symplectic momenta,  $\bx\cdot \by=\sum_{1\le i\le 3}x_i y_i$  and $|\bx|:=(\bx\cdot \bx)^{1/2}$ denote, respectively,  the standard inner product in $\real^3$ and the Euclidean norm;
\beq{masses}
\mu_i:=\frac{m_0 m_i}{m_0+ m_i}\,,\qquad M_i:=m_0+ m_i
\eeq
The phase space is the ``collisionless'' domain of $ \real^{3n}\times\real^{3n}$ 
\beqa{P6n}
\Big\{(\by,\bx)=\big((\by_1,\dots,\by_n), (\bx_1,\dots,\bx_n)\big) \ {\rm s.t.} \ \ \ 0\ne \bx_i\ne \bx_j\ ,\ \forall \ i\neq j\Big\}\ ,
\eeqa
endowed with the standard symplectic form 
$$\o=\sum_{i=1}^n d\by_i\wedge d\bx_i=\sum_{i=1}^n \sum_{j=1}^3 d\by_{ij}\wedge d\bx_{ij}$$
where $\by_{ij}$,  $\bx_{ij}$ denote the $j^ {\rm th}$ component of $\by_i$,  $\bx_i$.\\
The {\it planetary case} is when $m_1$, $\ldots$, $m_n$ are of the same order, and much smaller that $m_0$. In such a case, letting $m_i\to \mu m_i$, $\by_i\to \mu \by_i$, with $0<\mu\ll 1$, one obtains
\beqa{HelioNEW}{\cal H}&=\,&\sum_{1\leq i\leq n}\left(\frac{|\by_i|^2}{2\mu_i}-\frac{\mu_i M_i}{|\bx_i|}\right)+\mu\sum_{1\leq i<j\leq n}\left(\frac{\by_i\cdot \by_j}{ m_0}-\frac{m_i m_j}{|\bx_i-\bx_j|}\right)
\eeqa 
with \beq{massesNEW}
\mu_i:=\frac{m_0 m_i}{m_0+ \mu m_i}\,,\qquad M_i:=m_0+ \mu m_i
\eeq
\vskip.1in
\noindent
Consider the {\it  two--body} Hamiltonians \beqa{KeplerHam}h_i(\by_i, \bx_i):= \frac{|\by_i|^2}{2\m_i}-\frac{\m_i M_i}{|\bx_i|}\,.\eeqa
Assume that $h_i(\by_i, \bx_i)<0$ so that the Hamiltonian flow $\phi^t_{h_i}$ evolves on a Keplerian ellipse ${\cal E}_i$ and assume that the eccentricity  $e_i\in (0,1)$. Let $a_i$,  ${\mathbf P}_i$  denote, respectively,  the {\it  semimajor axis}  and the   {\it  perihelion} of ${\cal E}_i$.
Let ${\mathbf C}_i$ denote the $i^{\rm th}$ angular momentum
\beqa{Ci}\mathbf{C}_i(\by_j, \bx_j):=\bx_i\times \by_i\,.\eeqa Define the {\it Delaunay nodes}
\beqa{barni}\bar{\bm n}_i:=\mathbf k\times \mathbf C_i\eeqa
and,
for $u,v\in\real^3$ lying in the plane orthogonal to a vector $w$,  let $\a_w(u,v)$ denote the positively oriented angle (mod $2\p$) between $u$ and $v$  (orientation follows  the ``right hand rule'').
\vskip.1in
\noindent
The {\it  Delaunay action--angle variables}   
\beqa{Delaa}{\cal D}_{e\ell, aa}:= ({\mathbf Z}, {\mathbf G}, {\bm\Lambda}, \bm\zeta, {\mathbf g}, \bm\ell)\eeqa
with
\beqano 
\begin{array}{lll}
\dst{\mathbf Z}=(Z_1,\ldots,Z_{n}),\quad &\bm\zeta=(\zeta_1,\ldots,\zeta_{n})\\\\
\dst {\mathbf G}=(G_1,\ldots, G_{n}),& {\mathbf g}=(g_1,\ldots,g_{n})\\\\   
 \bm\L=(\L_1,\ldots,\L_n), &\bm\ell=(\ell_1,\ldots,\ell_{n})
\end{array}
 \eeqano
 are defined as
\beqa{Delaunay variables}
\left\{\begin{array}{l} \L_i:=\mu_i\sqrt{M_i a_i}\\
\ell_i:= {\rm mean\ anomaly\ of}\ \bx_i 
\ {\rm on}\ {\cal  E}_i
\end{array}\right. &&\left\{\begin{array}{l}
G_i:=|\mathbf{C}_i|=\L_i\sqrt{1-e_i^2}\\ 
g_i:=\a_{\mathbf{C}_i}(\bar{\bm n}_i, \mathbf P_i)
\end{array}
\right.\nonumber\\ \ \nonumber \\
&&\left\{\begin{array}{l}
Z_i:={\mathbf C}_i\cdot \mathbf{k}\\ 
\zeta_i:=\a_{\mathbf{k}}(\mathbf{i},\bar{\bm n}_i)
\end{array}
\right.\eeqa

\vskip.1in
\noindent
The {\it Poincar\'e variables} 
\beqano{\cal P}_{oinc}:= \big((\bm\uh, {\bm\up}, \bm\L), (\bm\ux,  {\bm\uq}, \bm\l)\big)\eeqano
with
\beqano 
\begin{array}{lll}
\dst{\bm\uh}=(\uh_1,\ldots,\uh_{n}),\quad &\bm\ux=(\ux_1,\ldots,\ux_{n})\\\\
\dst {\bm\up}=(\up_1,\ldots, \up_{n}),& {\bm\uq}=(\uq_1,\ldots,\uq_{n})\\\\   
 \bm\L=(\L_1,\ldots,\L_n), &\bm\l=(\l_1,\ldots,\l_{n})
\end{array}
 \eeqano
with the $\L_i$'s as in \equ{Delaunay variables} and
 \beqa{Poinc reg}
 \ul_i=\ell_i+ {g_i}+\theta_i\qquad&&\arr{\uh_i=\sqrt{2(\L_i-G_i)}\ \cos{(\zeta_i+g_i)}\\ \ux_i=-\sqrt{2(\L_i-G_i)}\ \sin{(\zeta_i+g_i})}
\nonumber\\
\\
&&
 \arr{\up_i=\sqrt{2(G_i-Z_i)}\ \cos{\zeta_i}\\ \uq_i=-\sqrt{2(G_i-Z_i)}\ \sin{\zeta_i}} 
 \nonumber
 \eeqa

\begin{figure}
\begin{tikzpicture}
\draw [ultra thick, ->] (0,0,0) -- (4,0,0);
\node (j) at (4.2,0,0) {$\mathbf j$};
\draw [ultra thick,->] (0,0,0) -- (0,4,0);
\node (k) at (0,4.2,0) {$\mathbf k$};
\draw [ultra thick,->] (0,0,0) -- (0,0,5);
\node (i) at (0,0,5.2) {$\mathbf i$};
\node (inew) at (0,0,1.8) {};
\draw [thick, ->] (0,0,0) -- (3.5,3.5,3.5);
\node (C) at (3.8,3.8,3.8) {$\mathbf{ C}_i$};
\node (Cnew) at (1.9,2.2,1.9) {$\tred{G_i}$};
\draw [dashed, -] (3.5,3.5,3.5) -- (3.5,0,3.5);
\draw [dashed, -] (0,0,0) -- (3.5,0,3.5);
\draw [dashed, -] (0,3.5,0) -- (3.5,3.5,3.5);
\draw [->] (-3.0,0,3.0) -- (3.0,0,-3.0);
\node (gamma) at (3.2,0,-3.2) {$\bar{\bm n}_i
$};
\node (gammanew) at (0.8,0,-0.8) {} ;
\draw [-latex, bend right] (inew) to (gammanew);
\node (Zeta) at (-0.2,2,0.2) {$\tred{Z_i}$} ;
\node (zeta) at (1.2,0,1.2) {$\tred{\zeta_i}$} ;
 \end{tikzpicture}
 \caption{Delaunay coordinates $Z_i$, $\zeta_i$, $G_i$.}
 \begin{tikzpicture}
\draw [ultra thick, ->] (0,0,0) -- (4,0,0);
\node (j) at (4.7,0,0) {$\mathbf{C}_i\times \bar{\bm n}_i$};
\draw [ultra thick,->] (0,0,0) -- (0,4,0);
\node (k) at (0,4.2,0) {$\mathbf{ C}_i$};
\draw [ultra thick,->] (0,0,0) -- (0,0,5);
\node (i) at (0,0,5.5) {$\bar{\bm n}_i$};
\node (inew) at (0,0,1.8) {};
\draw [dashed, -] (0,0,0) -- (5.6,0,7.3) 
;
\draw [thick, ->] (0,0,0) -- (3.57,0,2.75) 
;
\node (x2) at (4.0,0,3.0) 
{$\mathbf{x}_i$};
\node (ell2) at (2.0,0,1.5) 
{};
\node (C) at (6.1,0,8.0) 
{$\mathbf{P}_i$};
\node (g) at (1.0,0,2.5) {$\tred{g_i}$};
\node (g1) at (2,0,2.7) {};
\node (ell3) at (2.5,0,2.7) {$\tred{\ell_i}$};
\node (gammanew) at (0.8,0,-0.8) {} ;
\node (Zeta) at (-0.2,2,0.2) {$\tred{G_i}$} ;
\node (zeta) at (1.2,0,1.2) {} ;
\draw [-latex, bend right] (inew) to (zeta);
\draw [-latex, bend right] (g1) to (ell2);
\draw plot [variable=\t, domain=-78:63,
samples=50]
({3.5*(cos(\t)-0.2)+3.5*0.3*sin(\t)}, {-3.5*(cos(\t)-0.2)+3.5*0.3*sin(\t)});
 \end{tikzpicture}
 \caption{Delaunay coordinates $G_i$, $g_i$, $\ell_i$.}
 \end{figure}

\nl
In Poincar\'e coordinates  the Hamiltonian \equ{HelioNEW} takes the
form
\beq{prop deg}
\cH_{\textrm{\scshape p}}(\L,\ul,\uz)=h_{\textrm{\scshape k}} (\L)+ \mu f_{\textrm{\scshape p}}(\L,\ul, \uz)\ ,\ \ \uz:=(\uh,\up,\ux,\uq)\in\real^{4n}
\eeq
where  $(\L,\l)\in\real^n\times\torus^n$;
the ``Kepler'' unperturbed term $h_{\textrm{\scshape k}}$, coming from $h_{\rm plt} $ in \equ{Helio}, becomes
\beq{Kep}
h_{\textrm{\scshape k}}:=\sum_{i=1}^nh^\ppi_{\textrm{\scshape k}}(\L)=-\sum_{i=1}^n\frac{ \m_i^3 M_i^2}{2\L_i^2}\ .
\eeq

\nl
Because of rotation (with respect the ${\mathbf k}$--axis)  and reflection (with respect to the coordinate planes) invariance   of the Hamiltonian  \equ{Helio}, the perturbation $f_{\textrm{\scshape p}}$ in \equ{prop deg}
satisfies  well known symmetry relations called {\it d'Alembert rules}, see  \cite{chierchiaPi11c}. 
 By such symmetries, in particular, the 
averaged  perturbation
\beq{fav pl}f^{\rm av}_{\textrm{\scshape p}}(\L, \uz):=\su{(2\p)^n}\int_{\torus^n}f_{\textrm{\scshape p}}(\L,\ul,\uz) d\l\eeq
  is even around the origin $\uz=0$ and its expansion in powers of $\uz$  has the form\footnote{${\cal Q}\cdot u^2$ denotes the 2--indices contraction $\sum_{i,j}{\cal Q}_{ij}u_i u_j$ (${\cal Q}_{ij}$, $u_i$ denoting the entries of ${\cal Q}$, $u$). }
\beq{f'aaav}
f^{\rm av}_{\textrm{\scshape p}}=C_0(\L)+{\cal Q}_h(\L)\cdot\frac{{\uh}^2+{\ux}^2}{2}+{\cal Q}_v(\L)\cdot\frac{{\up}^2+{\uq}^2}{2}+{\rm O}(|\uz|^4)\ ,
\eeq 
where ${\cal Q}_h$, ${\cal Q}_v$ are suitable quadratic forms. The explicit expression of  such qua\-dra\-tic forms can be found, \eg, in   \cite[(36), (37)]{fejoz04}.

\nl
By such expansion, the (secular) origin $\uz=0$  is an {\it elliptic equilibrium} for $f^{\rm av}_{\textrm{\scshape p}}$ and 
corresponds to co--planar and co--circular motions.
It is therefore natural to put  \equ{f'aaav} into Birkhoff Normal Form (BNF, from now on) in a small neighborhood of the secular origin; see, \eg, \cite{hoferZ94} for general information on BNFs for Birkhoff theory for rotational invariant Hamiltonian systems.  

\nl
As a preliminary step, one can diagonalize \equ{f'aaav}, \ie, find 
a symplectic transformation defined by $\L\to\L$ and
 \beq{diag Poinc1}\ul=\tilde\ul+\varphi(\L,\tilde\uz),\ \uh=\r_h(\L)\tilde\uh,\ \ux=\r_h(\L)\tilde\ux,\ \up=\r_v(\L)\tilde\up,\ \uq=\r_v(\L)\tilde\uq\ ,
\eeq
with  $\r_h$, $\r_v\in {\rm SO}(n)$ diagonalizing ${\cal Q}_h$, ${\cal Q}_v$. In this way,  \equ{prop deg} takes the form 
\beq{planetary diag}
\tilde\cH_{\textrm{\scshape p}}(\L,\tilde\ul,\tilde\uz)=h_{\textrm{\scshape k}} (\L)+\m \tilde f(\L,\tilde\ul, \tilde\uz)\ ,
\eeq
with the average over $\tilde\ul$ of $\tilde f^{\rm av}$ given by
\beqa{f'aaavdiag}\tilde f^{\rm av}(\L,\tilde \uz)=C_0(\L)+\sum_ {i=1}^ {m}\Omega_i(\L)\frac{{\tilde u}_i^2+{\tilde v}_i^2}{2}+{\rm O}(|\tilde\uz|^4), \quad \tilde\uz=(\tilde u, \tilde v)=\big((\tilde\uh,\tilde\up)\,,\ (\tilde\ux,\tilde\uq)\big).\eeqa


\nl
with $m=2n$, and the   vector $\O(\L):=
(\s_1(\L),\ldots,\s_n(\L),\varsigma_1(\L), \ldots, \varsigma_n(\L))$ 
being formed by the eigenvalues of the matrices ${\cal Q}_h$ and ${\cal Q}_v$.

 \begin{theorem}[Birkhoff]\label{BNF} Let ${\cal H}$ be a  Hamiltonian having the form in \equ{planetary diag}--\equ{f'aaavdiag}. 
 Assume that there exists $\tilde\varepsilon>0$ $\cA\subset {\mathbb R}^n$ and $s\in \natural$ such that ${\cal H}$ is smooth on an open set $\tilde\cM^ {2m+2n}_{\varepsilon}=\cA\times\torus^n\times B^{2m}_{\tilde\varepsilon}$ and that 
\beqa{nonres}\sum_{i=0}^m\Omega_i(\Lambda)k_i\ne 0\quad\forall\ k=(k_1\,,\ldots\,,\ k_m)\in \integer^m:\ 0<|k|_1\le 2s\,,\ \forall\ \L\in \cA\,. \eeqa
 Then there exists $0<\varepsilon\le\tilde\varepsilon$ and a symplectic map (``Birkhoff transformation'')
 \beq{birkhoff transf}\Phi_{\textrm{\scshape b}}:\quad (\bm\L,{\mathbf l},\bar{\mathbf{w}})\in\cM^ {2m+2n}_{\varepsilon}\to (\L,\tilde\ul,\tilde\uz)\in \Phi_{\textrm{\scshape b}}(\cM^ {2m+2n}_{\varepsilon})\subseteq\cM^ {2m+2n}_{\tilde\varepsilon}
 \eeq
 which puts the  Hamiltonian \equ{planetary diag} into the form
 \beq{birkhoff planetary}
\cH_{\textrm{\scshape b}}(\bm\L,{\mathbf l},\bar{\mathbf{w}}):=\tilde\cH_{\textrm{\scshape p}}\circ\Phi_{\textrm{\scshape b}}=h_{\textrm{\scshape k}} (\L)+\m f_{\textrm{\scshape b}}(\L,l, w)\eeq
where the average $f_{\textrm{\scshape b}}^{\rm av}(\L,w):=\int_{\torus^n}f_{\textrm{\scshape b}}d l$ is in BNF of order $s$:
\beq{fb}
f_{\textrm{\scshape b}}^{\rm av}(\L,w)=C_0+\O\cdot r+{\rm P}_s(r)+{\rm O}(|w|^{2s+1})\quad w:=(u,v)\quad r_i:=\frac{u_i^2+v_i^2}{2}\ ,
\eeq
${\rm P}_s$ being homogeneous polynomial in $r$ of order $s$, with coefficients depending on $\L$.\\
In particular, if \equ{nonres} holds with $s=4$, 
\beq{fbNEW}
f_{\textrm{\scshape b}}^{\rm av}(\L,w)=C_0(\L)+\O(\L)\cdot r+r\cdot \tau(\L) r+{\rm O}(|w|^{5})\quad w:=(u,v)\quad r_i:=\frac{u_i^2+v_i^2}{2}\ ,
\eeq
with some square matrix $\tau(\L)$ of order $m$ (``torsion'', or ``second-order Birkhoff invariants'').
%
\end{theorem}

\begin{theorem}\label{FT}
{\bf (``The Fundamental Theorem'', V. I. Arnold, \cite{arnold63})}  If the Hessian matrix of ${\rm h}$ and the matrix $\t(\L)$  do not vanish identically, and if $\m$ is suitably small with respect to $\varepsilon$, the system affords a positive measure set $\cK_{\m, \varepsilon}$ of quasi--periodic motions in phase space such that its density goes to one as $\varepsilon\to 0$.
\end{theorem}

\begin{remark}[Arnold, Herman]\rm
It turns out that such invariants satisfy identically the following two  {\it secular resonances} 
\beq{Herman resonance}
\varsigma_ {n}(\Lambda)\equiv0\ ,\qquad\qquad\sum_ {i=1}^ {n}(\s_i(\Lambda)+\varsigma_i(\Lambda))\equiv 0
\eeq
Such resonances strongly violate the assumption \equ{nonres} of Theorem \ref{BNF}. \end{remark}
We remark that the former equality in \equ{Herman resonance} is  mentioned in \cite{arnold63}, while the latter been pointed out by M. Herman in the 1990s. 
Note that \equ{Herman resonance} do not appear in the planar problem, because the matrix ${\cal Q}_v$, hence the $\varsigma_i$'s, do not exist in that case. Being aware of such difficulty, Arnold completely proved Theorem \ref{Arnold Theorem} via Theorem \ref{FT} in the case of the planar three--body problem, checking explicitly the non vanishing of the $2\times 2$ torsion matrix for that case. However, in the case of the spatial problem, the question remained open until 2004, when M. Herman and J. F\'ejoz \cite{fejoz04}  proved Theorem \ref{Arnold Theorem} via a completely different strategy, which does need Birkhoff normal form. We refer to \cite{chierchiaPi14} for more details.

\subsection{The rotational degeneracy}
In \cite{arnold63},  Arnold wrote -- without giving the details --  that the former resonance in \equ{Herman resonance} was to be ascribed to  the conservation of the total angular momentum of the system:
\beqa{C}{\mathbf C}=\sum_{j=1}^n{\mathbf C}_j\,,\qquad {\mathbf C}_j=\bx_j\times \by_j\,.\eeqa 
An argument which clearly shows this goes as follows. Using Poincar\'e coordinates, the planets' angular momenta have the expressions
\beqano \mathbf C_j
&=&\left(
\begin{array}{ccc}
-\uq_j\sqrt{\Lambda_j-\frac{\uh^2_j+\ux^2_j}{2}-\frac{\up^2_j+\uq^2_j}{4}}\\
-\up_j\sqrt{\Lambda_j-\frac{\uh^2_j+\ux^2_j}{2}-\frac{\up^2_j+\uq^2_j}{4}}\\
\Lambda_j-\frac{\uh^2_j+\ux^2_j}{2}-\frac{\up^2_j+\uq^2_j}{2}
\end{array}
\right)\nonumber\\
&=&\left(
\begin{array}{ccc}
-\sqrt{\Lambda_j}\uq_j+{\rm O}(|\uz|^3)\\
-\sqrt{\Lambda_j}\up_j+{\rm O}(|\uz|^3)\\
\Lambda_j+{\rm O}(|\uz|^2)
\end{array}
\right)
\eeqano
In particular, the two former components of the total angular momentum \equ{C} are given by
\beqa{C1C2***}C_1=-\sum_{j=1}^n\sqrt{\Lambda_j}\uq_j+{\rm O}(|\uz|^3)\,,\qquad C_2=-\sum_{j=1}^n\sqrt{\Lambda_j}\up_j+{\rm O}(|\uz|^3)\eeqa
On the other hand, it is possible to find a canonical transformation 
\beqa{checktransf}(\Lambda, \check\ul, \check\uh, \check\up, \check\ux, \check\uq)\to (\Lambda, \ul, \uh, \up, \ux, \uq)\eeqa
having the form \equ{diag Poinc1} with $\rho_h=\id$ and $\rho_v\in SO(n)$ chosen such in a way that the last raw of $\rho^{-1}_v$ is
\beqa{vector}N(\L)\big(\sqrt{\Lambda_1}\,,\ldots\,,\sqrt{\Lambda_n}\big)
\eeqa
where $N(\L)=\frac{1}{\sqrt{\sum_{i=1}^n \Lambda_i}}$ fixes the Euclidean norm of \equ{vector} to $1$. With such choice, we have
$$\check\up_n=\rho^{-1}_v\left(
\begin{array}{cc}
\up_1\\
\vdots\\
\up_n
\end{array}
\right)_n=N(\L)\sum_{j=1}^n\sqrt{\Lambda_j}\up_j$$
and, similarly,
$$\check\uq_n=N(\L)\sum_{j=1}^n\sqrt{\Lambda_j}\uq_j$$
Therefore,  \equ{C1C2***} become
\beqa{C1C2NEW} C_1=-N(\L)^{-1}\check q_n+{\rm O}(|\check\uz|^3)\,,\qquad C_2=-N(\L)^{-1}\check p_n+{\rm O}(|\check\uz|^3)\eeqa
Now,  as the projection of the transformation \equ{checktransf} on $\check\ul$'s is a $\check\ul$--independent translation, the averaged perturbing function using the new coordinates can be obtained applying such transformation  to the function in \equ{f'aaav}. We denote it as
\beqano
\check f^{\rm av}=C_0(\L)+\check{\cal Q}_h(\L)\cdot\frac{\check{\uh}^2+\check{\ux}^2}{2}+\check{\cal Q}_v(\L)\cdot\frac{\check{\up}^2+\check{\uq}^2}{2}+{\rm O}(|\check\uz|^4)\ ,
\eeqano
with $\check{\cal Q}_h(\L)={\cal Q}_h(\L)$ and $\check{\cal Q}_v(\L)=\rho_v(\L)^{-1}{\cal Q}_v(\L)\rho_v(\L)$. 
Note that $\check{\cal Q}_v(\L)$ has the same eigenvalues as ${\cal Q}_v(\L)$, as $\rho_v\in SO(n)$.
Let us now use 
\beqa{commutation}\{\check f^{\rm av}, C_1\}=0=\{\check f^{\rm av}, C_2\}\eeqa
which hold because they are true for $f$, and $\mathbf C$ is $\check\ul$--independent.
Using \equ{C1C2NEW}, it is immediate to see that \equ{commutation} imply that
the quadratic form
$$\check{\cal Q}_v(\L)\cdot\frac{\check{\up}^2+\check{\uq}^2}{2}$$
is independent of $\check\up_n$, $\check\uq_n$. Hence, the $n^{\rm th}$ raw
and column of $\check{\cal Q}_v(\L)$ vanish identically. This implies that $\check{\cal Q}_v(\L)$, hence ${\cal Q}_v(\L)$,  has an identically vanishing eigenvalue, which is $\varsigma_n(\L)$ in \equ{Herman resonance}.
\subsection{Jacobi reduction of the nodes}
In the case $n=2$, Arnold in \cite{arnold63} suggested to get rid of the rotation invariance (described in the previous section) by means of the classical so--called {\it Jacobi reduction of the nodes}. This is a classical procedure with a remarkable geometric meaning, which goes as follows. 
 Let us consider a reference frame $({\mathbf i}, {\mathbf j}, {\mathbf k})$ whose third axis ${\mathbf k}$ is along the direction of the total angular momentum ${\mathbf C}={\mathbf C}_1+{\mathbf C}_2$, while ${\mathbf i}$ coincides with the intersection of the planes orthogonal to  ${\mathbf C}_1$, ${\mathbf C}_2$. Such intersection is well defined provided that ${\mathbf C}_1\not\parallel{\mathbf C}_2$, namely, when the problem is not planar. With such a choice of the reference frame, one cannot fix Delaunay coordinates completely freely. Indeed, by the choice of $\mathbf i$, we have that
 the $\zeta_j$ satisfy 
\beqa{Jacobi1}\zeta_2-\zeta_1=\p\,.\eeqa
 Moreover,
 a geometrical analysis of the triangle formed by ${\mathbf C}_1$, ${\mathbf C}_2$ and ${\mathbf C}$
 shows that the coordinates $Z_j$ satisfy
 \begin{figure}[htp]
		\center{\includegraphics[width=0.5\linewidth]{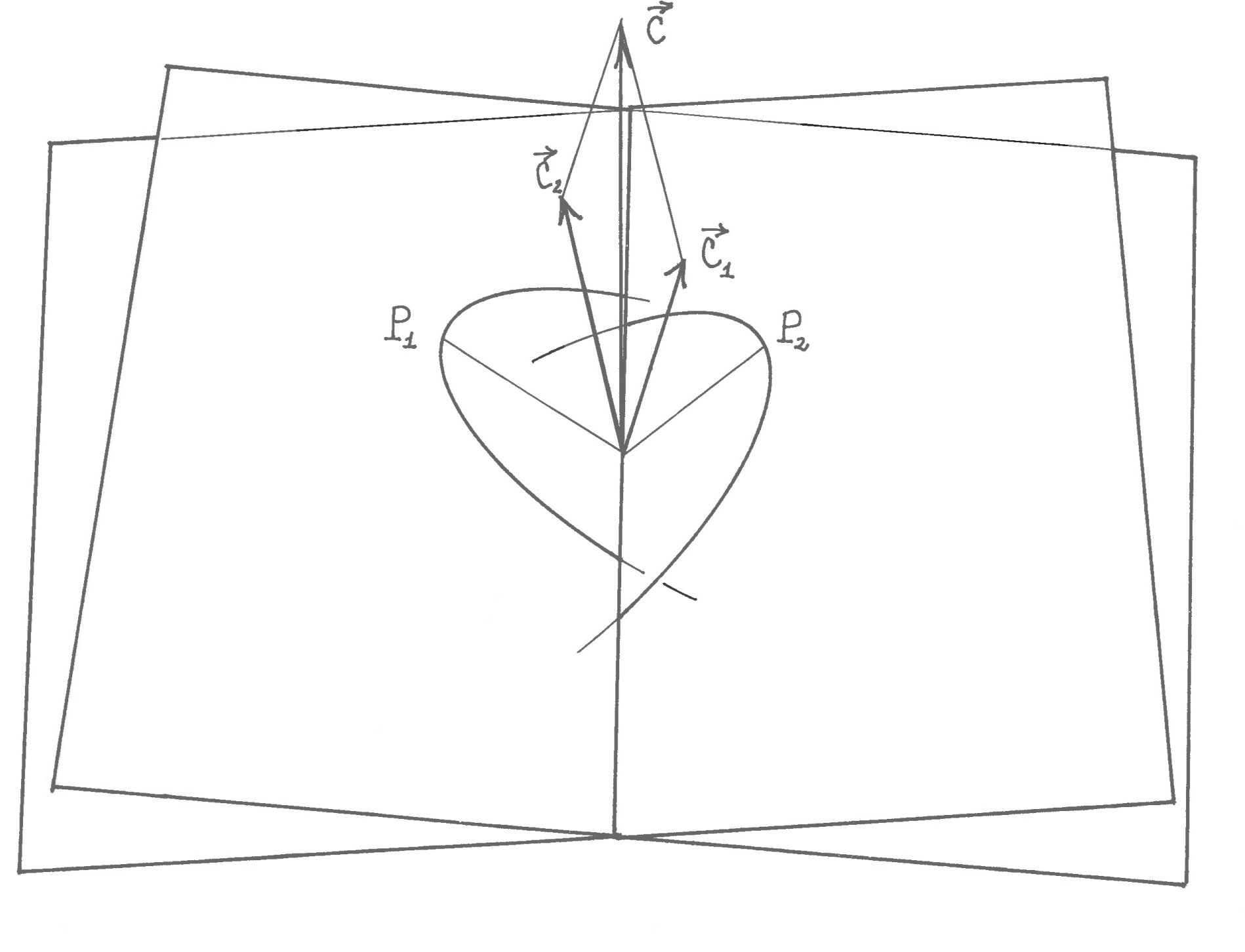}}
		\caption{The construction underlying Jacobi reduction of the nodes.}
	\end{figure}
\beq{Jacobi2}Z_1=\frac{G }{2}+\frac{G_1^2-G_2^2}{2G }\ ,\quad Z_2=\frac{G }{2}-\frac{G_1^2-G_2^2}{2G }\eeq
where
${G}:=|{\mathbf C}|=\sqrt{{ C}_1^2+{ C}_2^2+{ C}_3^2}$ is the Euclidean norm of ${\mathbf C}$. As $\mathbf i$ moves, the following fact is not obvious at all -- in fact proved by R. Radau.
\begin{theorem}[R. Radau, 1868, \cite{radau1868}]
Replacing relations \equ{Jacobi1}--\equ{Jacobi2} inside the Hamiltonian \equ{Helio} with $n=2$ written in Delaunay coordinates, one obtains a function, depending on  $(\Lambda_j, \ell_j, G_j, g_j)$ ($j=1$, $2$) and $G$, whose Hamilton equations relatively to $(\Lambda_j, \ell_j, G_j, g_j)$
generate the motions of the coordinates $(\Lambda_j, \ell_j, G_j, g_j)$ referred to the  rotating frame under the action of the Hamiltonian \equ{Helio} with $n=2$. The motion of $Z_j$ and $\zeta_j$ can be recovered via  \equ{Jacobi1}--\equ{Jacobi2}.
\end{theorem}

\subsection{Deprit coordinates}
Arnold commented on the general problem of rotational degeneracy  as follows:\vskip.1in \noindent
\cite[\bf Chap.III, \S 5, n. 5]{arnold63} {\it In the case of more than three bodies {\rm[$n> 2$]} there is no such {\rm[analogue to Jacobi reduction of the nodes]} elegant method of reducing the number of degrees of freedom} [...].

\vskip.1in \noindent
However, exactly 20 years later, in 1983, A. Deprit \cite{deprit83} discovered a set of canonical coordinates which, after a simple transformation, do the desired job and reduce to Jacobi's when $n=2$. Let us describe them.
\\
Consider the ``partial angular momenta''
 \beq{partial sums}
{\mathbf S}_j(\by, \bx):=\sum_{i=1}^j {\mathbf C}_j\ ;
\eeq
with $\mathbf C_i$ as in \equ{Ci}.
Notice that ${\mathbf S}_n={\mathbf C}$ is the total angular momentum of the system. Define the ``Deprit nodes'' 
\beqa{nodes}
\left\{\begin{array}{l}
\bm\n_{i+1}:=
{\mathbf S}_{i+1}\times {\mathbf C}_{i+1}\ ,\qquad\ \ 1\le i\le n-1\\ 
\bm\n_1:={\mathbf S}_{2}\times {\mathbf C}_{1}=-\bm\n_2\\
\bm\n_{n+1}:=
{\mathbf k}\times {\mathbf C}=:\bar{\bm\n}\ .
\end{array}
\right.
\eeqa

\noindent
If $n\ge 2$,
Deprit's coordinates \beqa{Deprit coordinates}{\cal D}_{ep}=(\mathbf R, \mathbf G,\bm\Psi, \mathbf r,\bm\varphi,\bm\psi)\eeqa
with
\beqa{Deparg}
&&\mathbf R=(R_1, \ldots, R_n)\,,\ \bm\Psi=(\Psi_1, \ldots, \Psi_n)\,,\ \mathbf G=(G_1, \ldots, G_n)\,,\nonumber\\
&& \mathbf r=(r_1, \ldots, r_n)\,,\ \ \ \  \bm \psi=(\psi_1, \ldots, \psi_n)\,,\ \ \bm\varphi=(\varphi_1, \ldots, \varphi_n)\,.
\eeqa
are defined as follows (compare also Figures \ref{Deprit1}, \ref{Deprit2} and \ref{Deprit3}):

\beqa{Deprit variables}
&&\left\{\begin{array}{l}\dst R_i:=\by_i\cdot\frac{\bx_i}{|\bx_i|}\\\\
\dst r_i:= |\bx_i| 
\end{array}\right.\qquad\left\{\begin{array}{l}
\dst G_i:=|{\mathbf C}_i|\\ \\
\dst \varphi_i:=\a_{{\mathbf C}_i}(\bm\n_i,\bx_i)
\end{array}
\right.\nonumber 
\\  \\
&&\Psi_i:=\left\{\begin{array}{l}
\dst |\mathbf S_{i+1}|\phantom{AAAAAAAA}\ 1\le i\le n-2\ (n\ge 3)
\\ \\
\dst C:=|\mathbf C|\phantom{AAAAAA}\ \ i=n-1\\\\
\dst Z:={\mathbf C}\cdot \mathbf k\phantom{AAAA.} \ \ \ i=n
\end{array}\right. \nonumber\\ \nonumber\\
&&\psi_i:=\left\{\begin{array}{l}\dst \a_{{\mathbf S}_{i+1}}(\bm\n_{i+2},\bm\n_{i+1})\phantom{AAAa} 1\le i\le n-2\ (n\ge 3)\\\\
\gamma:= \a_{\mathbf C}(\bar{\bm\n},\bm\n_{n})\phantom{AAAAA}  i= n-1\\\\
\dst\zeta:=\a_{\mathbf k}(\mathbf i, \bar{\bm\n}) \phantom{AAAAAA}\ i=n
\end{array}\right.\nonumber
\eeqa

\begin{figure}
\begin{tikzpicture}
\draw [ultra thick, ->] (0,0,0) -- (4,0,0);
\node (j) at (4.2,0,0) {$\mathbf j$};
\draw [ultra thick,->] (0,0,0) -- (0,4,0);
\node (k) at (0,4.2,0) {$\mathbf k$};
\draw [ultra thick,->] (0,0,0) -- (0,0,5);
\node (i) at (0,0,5.2) {$\mathbf i$};
\node (inew) at (0,0,1.8) {};
\draw [thick, ->] (0,0,0) -- (3.5,3.5,3.5);
\node (C) at (3.8,3.8,3.8) {$\mathbf C$};
\node (Cnew) at (1.9,2.2,1.9) {$\tred C$};
\draw [dashed, -] (3.5,3.5,3.5) -- (3.5,0,3.5);
\draw [dashed, -] (0,0,0) -- (3.5,0,3.5);
\draw [dashed, -] (0,3.5,0) -- (3.5,3.5,3.5);
\draw [->] (-3.0,0,3.0) -- (3.0,0,-3.0);
\node (gamma) at (3.2,0,-3.2) {$\bar{\bm\nu}$};
\node (gammanew) at (0.8,0,-0.8) {} ;
\draw [-latex, bend right] (inew) to (gammanew);
\node (Zeta) at (-0.2,2,0.2) {$\tred Z$} ;
\node (zeta) at (1.2,0,1.2) {$\tred\zeta$} ;
 \end{tikzpicture}
 \caption{Deprit coordinates $Z$, $C$ and $\zeta$ fix the angular momentum in the initial reference frame $({\mathbf i}, {\mathbf j}, {\mathbf k})$.}\label{Deprit1}
 \begin{tikzpicture}
\draw [ultra thick, ->] (0,0,0) -- (4,0,0);
\node (j) at (5.0,0,0) {$\mathbf S_{i+1}\times \bm\nu_{i+2}$};
\draw [ultra thick,->] (0,0,0) -- (0,4.2,0);
\node (k) at (-0.4,2.2,0) {$\mathbf S_{i+1}$};
\node (k) at (0.4,2.2,0) {$\tred{\Psi_i}$};
\draw [ultra thick,->] (0,0,0) -- (0,0,5);
\node (i) at (0,0,5.5) {$\bm\nu_{i+2}$};
\node (inew) at (0,0,1.8) {};
\draw [thick, ->] (0,0,0) -- (3.5,3.5,3.5);
\node (C) at (3.0,2.5,3.8) {$\mathbf C_{i+1}$};
\draw [thick, ->] (3.5,3.5,3.5) -- (0,4.2,0);
\node (Gi) at (2.3,4.5,3.7) {$\tred{\Psi_{i-1}}$};
\node (Ginew) at (3.0,4.8,3.7) {$\mathbf S_i$};
\node (Cnew) at (1.7,2.3,1.9) {$\tred{G_{i+1}}$};
\draw [dashed, -] (3.5,3.5,3.5) -- (3.5,0,3.5);
\draw [dashed, -] (0,0,0) -- (3.5,0,3.5);
\draw [dashed, -] (0,3.5,0) -- (3.5,3.5,3.5);
\draw [->] (-3.0,0,3.0) -- (3.0,0,-3.0);
\node (gamma) at (3.2,0,-2.0) {$\bm\nu_{i+1}$};
\node (gammanew) at (0.8,0,-0.8) {} ;
\draw [-latex, bend right] (inew) to (gammanew);
\node (Zeta) at (-0.2,2,0.2) {} ;
\node (zeta) at (1.2,0,1.2) {$\tred{\psi_i}$} ;
 \end{tikzpicture}
 \caption{{The frames $\DD_{i+1}$ and the} coordinates $\Psi_i$, $\Psi_{i-1}$, $G_{i+1}$ and $\psi_{i}$.}\label{Deprit2}
 \begin{tikzpicture}
\draw [ultra thick, ->] (0,0,0) -- (4,0,0);
\node (j) at (4.7,0,0) {$\mathbf C_i\times { \bm\nu_i}$};
\draw [ultra thick,->] (0,0,0) -- (0,4,0);
\node (k) at (0,4.2,0) {$\mathbf C_i$};
\draw [ultra thick,->] (0,0,0) -- (0,0,5);
\node (i) at (0,0,5.5) {${\bm\nu_i}$};
\node (inew) at (0,0,1.8) {};
\draw [dashed, -] (0,0,0) -- (5.6,0,7.3) 
;
\draw [thick, ->] (0,0,0) -- (3.57,0,2.75) 
;
\node (x2) at (4.0,0,3.0) 
{$\mathbf x_i$};
\node (x2new) at (3.0,0,2.25) 
{};
\node (ell2) at (2.0,0,1.5) 
{};
\node (C) at (6.1,0,8.0) 
{$\mathbf P_i$};
\node (g) at (1.0,0,2.5) {$\tred{g_i}$};
\node (g1) at (2,0,2.7) {};
\node (ell3) at (2.5,0,2.7) {$\tred{\ell_i}$};
\node (varphi) at (1.9,0,3.7) {$\tred{\varphi_i}$};
\node (gammanew) at (0.8,0,-0.8) {} ;
\node (Zeta) at (-0.2,2,0.2) {$\tred{G_i}$} ;
\node (zeta) at (1.2,0,1.2) {} ;
\draw [-latex, bend right] (inew) to (zeta);
\draw [-latex, bend right] (g1) to (ell2);
\draw [-latex, bend right] (inew) to (x2new);
\draw plot [variable=\t, domain=-78:63,
samples=50]
({3.5*(cos(\t)-0.2)+3.5*0.3*sin(\t)}, {-3.5*(cos(\t)-0.2)+3.5*0.3*sin(\t)});
 \end{tikzpicture}
  \caption{{The frames ${\rm H}_{i}$ and the}  coordinates $g_i$, $G_{i}$, $\ell_i$.}\label{Deprit3}
\end{figure}

\noindent
We have
\begin{theorem}[A. Deprit, 1983, \cite{deprit83}]\label{Deprit theorem}
$\sum_{i=1}^n\by_i\cdot d\bx_i={\mathbf R}\cdot d\mathbf r+\mathbf \Psi\cdot d\bm\psi+\mathbf G\cdot d\bm\varphi$ for all $n\in \natural$.
\end{theorem}
For later need, we formulate an equivalen statement of Theorem \ref{Deprit theorem}. We consider the coordinates
\beqa{Del}{\cal D}_{e\ell}:= (\mathbf Z, \mathbf G, \mathbf R, \bm\zeta, \bm\phi, \mathbf r)\eeqa
with
\beqa{Delarg}
&&\mathbf Z=(Z_1, \ldots, Z_n)\,,\ \mathbf G=(G_1, \ldots, G_n)\,,\ \mathbf R=(R_1, \ldots, R_n)\nonumber\\
&&\bm \zeta=(\zeta_1, \ldots, \zeta_n)\,,\ \ \  \bm\phi=(\phi_1, \ldots, \phi_n)\,,\ \mathbf r=(r_1, \ldots, r_n)
\eeqa
where  $Z_i$, $G_i$,  $\zeta_i$,  are as in
\equ{Delaunay variables},  $R_i$, $r_i$ are as in \equ{Deprit variables}, and, finally,
\beqano \phi_i:=\a_{\mathbf{C}_i}({\mathbf n}_i, \mathbf x_i)\,.\eeqano
Let
\beq{elementary rotations}\cR_{1}(i)=\left(
\begin{array}{cccc} 
1&0&0\\
0&\cos i&-\sin i\\
0&\sin i&\cos i
\end{array}
 \right)\ ,\qquad \cR_{3}(\theta)=\left(
\begin{array}{ccc} 
\cos\theta&-\sin\theta&0\\
\sin\theta&\cos\theta&0\\
0&0&1
\end{array}
 \right)\eeq
 and
 $$\bx=\cR_3(\theta)\cR_1(i)\bar{\bx}\ ,\quad \by=\cR_3(\theta)\cR_1(i)\bar {\by}\ ,\quad {\mathbf C}:=\bx\times \by\ ,\quad \bar{\mathbf C}:=\bar{\bx}\times \bar{\by}\ ,\quad \mathbf i=\left(
\begin{array}{ccc}
1\\
0\\
0
\end{array}
\right)\ ,\quad \mathbf k=\left(
\begin{array}{ccc}
0\\
0\\
1
\end{array}
\right)$$
with $\bx,\bar{\bx}, \by, \bar{\by}\in\real^3$. The proof of the following fact is left to the reader.
\begin{lemma}\label{LemmaD}
$\by\cdot d\bx={\mathbf C}\cdot \mathbf kd\theta+\bar{\mathbf  C}\cdot \mathbf i di+\bar{\by} \cdot d\bar{\bx}$.

\end{lemma}

\noindent
Lemma \ref{LemmaD} immediately implies
\begin{lemma}\label{LemmaDel}
$\mathbf y_j\cdot d\mathbf x_j=Z_jd\zeta_j+G_jd\phi_j+R_jdr_j\quad \forall\ j=1\,,\ldots\,, n\,,\ \forall n\in \natural$.
\end{lemma}
Indeed, we have
$$\left\{\begin{array}{lll}\mathbf x_j={\cal R}_3(\zeta_j){\cal R}_1(i^*_j){\mathbf x}_j^*\\
\mathbf y_j={\cal R}_3(\zeta_j){\cal R}_1(i^*_j){\mathbf y}_j^*
\end{array}
\right.
\quad j=1\,,\ldots\,,\ n$$
where 
$i^*_j$ is the convex angle formed by $\mathbf k$ and $\mathbf C_j$ and, finally, 
\beqa{xjyjpl}{\mathbf x}_j^*=\left(
\begin{array}{ccc}
r_j\cos\phi_j\\
r_j\sin \phi_j\\
0
\end{array}
\right)\,,\quad {\mathbf y}_j^*=\left(
\begin{array}{ccc}
R_j\cos\varphi_j-\frac{G_j}{r_j}\sin\phi_j\\
R_j\sin \varphi_j+\frac{G_j}{r_j}\cos\phi_j\\
0
\end{array}
\right) \eeqa
verify, as well known,
\beqa{Kepler symplectic}{\mathbf y}_j^*\cdot d{\mathbf x}_j^*=R_j dr_j+G_j d\phi_j\,.\eeqa
 Then, by Lemma \ref{LemmaD}, \equ{Kepler symplectic} and as ${\mathbf C_j}^*\cdot \mathbf i=0$, we have
\beqano
\mathbf y_j\cdot d\mathbf x_j&=&{\mathbf C}_j\cdot \mathbf kd\zeta_j+{\mathbf C_j}^*\cdot \mathbf i di_j+{\mathbf y}_j^* \cdot d{\mathbf x}_j^*\nonumber\\
&=&Z_jd\zeta_j+G_jd\phi_j+R_jdr_j\,.\qquad \square
\eeqano
We denote as 
$$\phi_{\cD_{e\ell}}^{{\cal D}_{ep}}:\quad {\cal D}_{e\ell}=(\mathbf Z, \mathbf G, \mathbf R, \bm\zeta, \bm\phi, \mathbf r)\to {\cal D}_{ep}=(\bm\Psi, \mathbf G, \mathbf R, \bm\psi, \bm\varphi, \mathbf  r)$$
the map which relates $\cD_{e\ell}$ and $\cD_{ep}$ and as
$$\widehat\phi_{\cD_{e\ell}}^{{\cal D}_{ep}}:\quad \widehat{\cal D}_{e\ell}=(\mathbf Z, \mathbf G, \bm\zeta, \bm\phi)\to \widehat{\cal D}_{ep}=(\bm\Psi, \mathbf G,  \bm\psi, \bm\varphi)$$
is the natural projections on the coordinates above. It is easy to check that $\widehat\phi_{\cD_{e\ell}}^{{\cal D}_{ep}}$ is independent of $\mathbf R$ and $\mathbf r$.
Indeed, $\widehat\phi_{\cD_{e\ell}}^{{\cal D}_{ep}}$ has the expression
\beqa{reduced Deprit} 
G_j&=&G_j\ ,\nonumber\\ \nonumber \\
\varphi_j&=&\phi_j+\a_{{\mathbf C}_i}(\bm{\n}_j,\bar{\bm \n}_j)\ {\rm with}\ \bar{\bm \n}_j={\mathbf k}\times {\mathbf C}_j,\nonumber\\ \nonumber \\
 \Psi_j&=&\left\{\begin{array}{lll}|\mathbf S_{j+1}|\quad &j\ne n\\
Z_1+\ldots+Z_n\ &j=n\end{array}\right.\quad \nonumber\\ 
\psi_j&=&\left\{\begin {array}{lll}\a_{\mathbf S_{j+1}}(\bm{\n}_{j+2},\bm{\n}_{j+1})\quad &j\ne n\\
\a_{\mathbf  k}(\mathbf i,\bar{\bm\nu})&j=n
\end{array}
\right.
\eeqa
where {${\mathbf S_{j+1}}$, $\bm{\n}_{j}$, $\ovl{\bm{\n}}_{j}$ at} the right hand sides are to be written as functions of ${\cal D}_{e\ell}$ ({see \equ{nodes} and \equ{Delaunay variables}}):
$$\left\{\begin{array}{lll}\dst{{\mathbf S_{j+1}}=\sum_{i=1}^{j+1}G_i\cR_3(\zeta_i)\cR_1(i^{\cD el}_i){\mathbf k}}\\\\
\dst{\bm\n_{j+1}=
\left(\sum_{i=1}^{j+1}G_i\cR_3(\zeta_i)\cR_1(i^{\cD el}_i){\mathbf k}\right)\times G_{j+1}\cR_3(\zeta_{j+1})\cR_1(i^{\cD el}_{j+1}){\mathbf k}\ ,\qquad\ \ 1\le j\le n-1}\\ \\
\dst{\bm\n_1=-\bm\n_2=\left(\sum_{i=1}^{j+1}G_i\cR_3(\zeta_i)\cR_1(i^{\cD el}_i){\mathbf k}\right)\times G_1\cR_3(\zeta_1)\cR_1(i^{\cD el}_1){\mathbf k}}\\\\
\dst {\bm\n_{n+1}=\bar{\bm\n}
{\mathbf k}\times \left(\sum_{i=1}^{n}G_i\cR_3(\zeta_i)\cR_1(i^{\cD el}_i){\mathbf k}\right)}\ .
\end{array}
\right.$$
{with $ i_i^{\cD el}=\cos^{-1}\frac{G_i}{Z_i}$.}
 As the right hand sides are defined only in terms of $\mathbf C_j$, so they are functions of $\mathbf Z$, $\bm\zeta$ and $\mathbf G$, while are independent of $\mathbf R$ and $\mathbf r$.

\begin{theorem}\label{Deprit theoremNEW}
Theorem \ref{Deprit theorem} is equivalent to stress that 
\beqa{New statement}\ \widehat\phi_{\cD_{e\ell}}^{{\cal D}_{ep}}\ verifies:\quad \mathbf Z\cdot d\bm\zeta+\mathbf G\cdot d\bm\phi=\mathbf \Psi\cdot d\bm\psi+\mathbf G\cdot d\bm\varphi\ \  for\ all\  n\in \natural\,.\eeqa\end{theorem}
\proof Use Lemma \ref{LemmaDel} and that the coordinates $(\mathbf R, \mathbf r)$ are shared by
${\cal D}_{ep}$ and ${\cal D}_{e\ell}$. $\quad \square$

\vskip.1in
\noindent
We prove Theorem \ref{Deprit theorem} ($\Longleftrightarrow$ \equ{New statement})  by induction on $n$, with $n\ge 2$, as in \cite{pinzari-th09}.

\vskip.1in
\noindent
{\bf Base step} We prove the statement \ref{Deprit theorem} with $n=2$. We first observe that, in such case, $(\by_j, \bx_j)$ are expressed, through $(\mathbf R, \bm \Psi, \mathbf G, \mathbf r, \bm\psi, \bm\varphi)$ via the formulae
$$\left\{
\begin{array}{lll}\mathbf x_j={\cal R}_3(\zeta){\cal R}_1(i)\cR_3(\gamma)\cR_1(i_j){\mathbf x_j}_{\rm pl}\\
\mathbf y_j={\cal R}_3(\zeta){\cal R}_1(i)\cR_3(\gamma)\cR_1(i_j){\mathbf y_j}_{\rm pl}
\end{array}
\right.\quad j=1\,,\ 2$$
where 
$i$ is the convex\footnote{The expressions of $i_1$, $i_2$ and $i$ -- not needed here -- can easily be deduced by the analysis of the triangle formed by $\mathbf C_1$, $\mathbf C_2$ and $\mathbf C$: see Figure \ref{Deprit2}} angle formed by $\mathbf k$ and $\mathbf C$; $i_j$ is the convex angle formed by $\mathbf C$ and $\mathbf C_j$
and, finally, ${\mathbf x_j}_{\rm pl}$, $ {\mathbf y_j}_{\rm pl}$ are as in \equ{xjyjpl}, with $\phi_j$ replaced by $\varphi_j$.\\
Using  Lemma \ref{LemmaD} twice, one easily finds
\beqa{C1}
\by_j\cdot d{\bx}_j=&&{{\mathbf C}_j}\cdot {\mathbf k}\,d\zeta+{\bar{\mathbf C}_j}\cdot {\mathbf i}\,di +{\bar{\mathbf C}_j}\cdot {\mathbf k}\,d\gamma+\ {{\mathbf C}_j}_{\rm pl}\cdot {\mathbf i}\,d(i_j)\nonumber\\
&&+{\by_j}_{\rm pl}\cdot d{{\bx}_j}_{\rm pl}\nonumber\\
=&&{{\mathbf C}_j}\cdot {\mathbf k}\,d\zeta+{{\mathbf C}_j}\cdot {\mathbf e}_1\,di +{{\mathbf C}_j}\cdot {\mathbf e}_3\,d\gamma+{\by_j}_{\rm pl}\cdot d{{\bx}_j}_{\rm pl}
\eeqa
 We have used ${{\mathbf C}_j}_{\rm pl}\cdot {\mathbf i}=0$, ${{\mathbf C}_j}=\cR_3(\zeta)\cR_1(i ){\bar{\mathbf C}_j}$ and we have let \beq{e1e3}\mathbf e_1:=\cR_3(\zeta)\cR_1(i )\mathbf i\,,\qquad \mathbf e_3:=\cR_3(\zeta)\cR_1(i )\mathbf k\ .\eeq
Taking the sum of \equ{C1} with $j=1$, $2$ and using \equ{Kepler symplectic} and  recognizing that
$$\arr{({{\mathbf C}_1}+{\mathbf C}_2)\cdot {\mathbf k}={\mathbf C}\cdot {\mathbf k}=Z\\
 ({{\mathbf C}_1}+{\mathbf C}_2)\cdot {\mathbf e}_1= {\mathbf C}\cdot {\mathbf e}_1=0\\
 ({{\mathbf C}_1}+{\mathbf C}_2)\cdot {\mathbf e}_3={\mathbf C}\cdot {\mathbf e}_3=C}$$
we have the proof. \qed
\vskip.1in
\noindent
{\bf Induction} The inductive step is made on the statement  \equ{New statement}.
The map $\widehat\phi_{\cD_{e\ell}}^{{\cal D}_{ep}}$ in \equ{New statement}  
will be named $\widehat\phi_n$. We assume that \equ{New statement}
holds for a given $n\ge 2$ and prove it for $n+1$.  
Consider the map
$$\phi^*_{n+1}:\quad \widehat{\cal D}_{e\ell, n+1}=(\mathbf Z, \mathbf G,  \bm\zeta, \bm\phi)\to \widetilde{\cal D}_{ep, n+1}=(\bm\Psi^*, \mathbf G^*,  \bm\psi^*, \bm\varphi^*)$$
defined as follows. If
 $$\mathbf Z=\big(\widetilde{\mathbf Z}, Z_{n+1}\big)\,,\ \mathbf G=\big(\widetilde{\mathbf G}, G_{n+1}\big)\,,\ \bm \zeta=\big(\widetilde{\bm \zeta}, \zeta_{n+1}\big)\,,\ \bm \phi=\big(\widetilde{\bm \phi}, \phi_{n+1}\big)$$ where the tilded arguments have dimension  {$n$}, we let
$$(\widetilde{\bm \Psi}, \widetilde{\mathbf G},\widetilde{\bm \psi}, \widetilde{\bm \varphi})=\phi_{n}(\widetilde{\mathbf Z}, \widetilde{\mathbf G},\widetilde{\bm \zeta}, \widetilde{\bm \phi})$$
and then
$$\phi^*_{n+1}(\mathbf Z, \mathbf G,  \bm\zeta, \bm\phi):=\big((\widetilde{\bm\Psi}, Z_{n+1}), (\widetilde{\mathbf G}, G_{n+1}),  (\widetilde{\bm\psi}, \zeta_{n+1}), (\widetilde{\bm\varphi}, \phi_{n+1})\big)=:(\bm\Psi^*, \mathbf G^*,  \bm\psi^*, \bm\varphi^*)$$
By the inductive assumption, $\phi_n$ verifies
$$\widetilde{\mathbf Z}\cdot d\widetilde{\bm\zeta}+\widetilde{\mathbf G}\cdot d\widetilde{\bm\phi}=\widetilde{\mathbf \Psi}\cdot d\widetilde{\bm\psi}+\widetilde{\mathbf G}\cdot d\widetilde{\bm\varphi}$$
and hence $\phi^*_{n+1}$ verifies
\beqa{canonical1}
{\mathbf Z}\cdot d{\bm\zeta}+{\mathbf G}\cdot d{\bm\phi}&=&{\mathbf \Psi}^*\cdot d{\bm\psi}^*+{\mathbf G}^*\cdot d{\bm\varphi}^*\nonumber\\&=&\widetilde{\mathbf \Psi}\cdot d\widetilde{\bm\psi}+\widetilde{\mathbf G}\cdot d\widetilde{\bm\varphi}+Z_{n+1}d\zeta_{n+1}+G_{n+1}d\phi_{n+1}\nonumber\\
&=&\left(\sum_{j=1}^{n-2}\widetilde{ \Psi}_j\cdot d\widetilde{\psi}_j+\widetilde{\mathbf G}\cdot d\widetilde{\bm\varphi}\right)
+\widetilde{ \Psi}_{n-1}\cdot d\widetilde{\psi}_{n-1}+\widetilde{ \Psi}_{n}\cdot d\widetilde{\psi}_{n}\nonumber\\
&+&Z_{n+1}d\zeta_{n+1}+G_{n+1}d\phi_{n+1}\quad {\rm split\ the\ rhs}
\eeqa
having split
\beqano
\widetilde{\bm\Psi}=(\widetilde\Psi_1, \ldots, \widetilde\Psi_n)\,,\ \widetilde{\bm \psi}=(\widetilde\psi_1, \ldots, \widetilde\psi_n)\,.\eeqano
We moreover define a map $\phi_{*, n+1}$ on $(\bm\Psi^*, \mathbf G^*,  \bm\psi^*, \bm\varphi^*)$ acting as 
$$(\bm\Psi_*, \mathbf G_*,  \bm\psi_*, \bm\varphi_*)=\phi_{2}\big((\widetilde{\Psi}_{n}, Z_{n+1}), (\widetilde{ \Psi}_{n-1}, G_{n+1}),  (\widetilde{\psi}_n, \zeta_{n+1}), (\widetilde{\psi}_{n-1}, \phi_{n+1})\big)$$
on the designed variables, and as the identity on the remaining ones. Note that the arguments at left hand side have dimension $2$,  that
$\mathbf G_*=(\widetilde{ \Psi}_{n-1}, G_{n+1})$, and put $\bm\varphi_*=(\varphi_{*, 1}, \varphi_{*, 2})$.
Again by the inductive assumption, we have
\beqa{canonical2}\widetilde{ \Psi}_{n-1}\cdot d\widetilde{\psi}_{n-1}+\widetilde{ \Psi}_{n}\cdot d\widetilde{\psi}_{n}
+Z_{n+1}d\zeta_{n+1}+G_{n+1}d\phi_{n+1}={\mathbf \Psi}_*\cdot d{\bm\psi}_*+{\mathbf G}_*\cdot d{\bm\varphi}_*\eeqa
Let us now look at the composition \beqa{composition}\phi_{*, n+1}\circ\phi^*_{n+1}\eeqa It acts as
\beqano
 (\mathbf Z, \mathbf G,  \bm\zeta, \bm\phi)&\to& \big(
(\widetilde\Psi_1, \ldots, \widetilde\Psi_{n-2},\widetilde\Psi_{n-1}, \bm\Psi_*), (\widetilde{\mathbf G}, G_{n+1}), (\widetilde\psi_1, \ldots, \widetilde\psi_{n-2},{\varphi}_{*1},  \bm\psi_*), (\widetilde{\bm\varphi}, {\varphi}_{*2})\nonumber\\
&&=:(\bm\Psi, \mathbf G,  \bm\psi, \bm\varphi)
\big)
\eeqano
and, by \equ{canonical1} and \equ{canonical2}, verifies
\beqano
{\mathbf Z}\cdot d{\bm\zeta}+{\mathbf G}\cdot d{\bm\phi}&=&\left(\sum_{j=1}^{n-2}\widetilde{ \Psi}_j\cdot d\widetilde{\psi}_j+\widetilde{\mathbf G}\cdot d\widetilde{\bm\varphi}\right)
+{\mathbf \Psi}_*\cdot d{\bm\psi}_*+{\mathbf G}_*\cdot d{\bm\varphi}_*
\nonumber\\
&=&{\mathbf \Psi}\cdot d{\bm\psi}+{\mathbf G}\cdot d{\bm\varphi}\,.
\eeqano
It is not difficult to recognize -- using \equ{reduced Deprit}  -- that the map \equ{composition} coincides with $\phi_{n+1}$. For the details, we refer to \cite{pinzari-th09, chierchiaPi11a}. $\quad \square$

\paragraph{{The Deprit map}}
{In this section, we provide the explicit expression of the map 
\beqa{cartesianToDeprit}\phi_{\cC}^{{\cal D}_{ep}}:\quad {\cal C}=(\mathbf y_1, \ldots, \mathbf y_n, \mathbf x_1, \ldots, \mathbf x_n)\to {\cal D}_{ep}=(\bm\Psi, \mathbf G, \mathbf R, \bm\psi, \bm\varphi, \mathbf  r)\,.\eeqa
The discussion in the previous section shows that each orbital frame $\HH_i$, $i=1$, $\ldots$, $n$, can be reached via  a sequence of transformations which overlap the $\DD_{n+1}:=({\mathbf i}, {\mathbf j}, {\mathbf k})$ to $\HH_i$ through the following diagram (named {\it tree} by Deprit):
\beqano
\begin{array}{cccccccccccccccccc}
\dst \DD_{n+1}&\to&\DD_n&\to& \DD_{n-1}&\to&\cdots&\to& \DD_2&\to{\rm H}_1\\ \\
\dst&&\downarrow&&\downarrow&&\vdots&&\downarrow&&\\ \\
\dst&&{\rm H}_{n}&&{\rm H}_{n-1}&&\vdots&&{\rm H}_2&&&\\ \\
\end{array}
\eeqano
In turn, 
\begin{itemize}
\item[--] the transition $\DD_{n+1}\to \rm D_{n}$ is described by the sequence of rotations $\cR_3(\psi_{n})\cR_1(i_{n})$, with $\cos i_n=\frac{Z}{\rm G}=\frac{\Psi_{n}}{\Psi_{n-1}}$ (see figure \ref{Deprit1});
\item[--] the transitions $\DD_{i+1}\to \rm H_{i+1}$, $i+1=n-1$, $\ldots$, $2$,  are described by the sequence of rotations $\cR_3(\psi_i)\cR_1(i_i)$, with $\cos i_i=\frac{\Psi_{i}^2+G_i^2-\Psi_{i-1}^2}{2\Psi_{i} G_{i+1}}$ (see figure \ref{Deprit2});
\item[--] the transitions $\DD_{i+1}\to \DD_{i}$, $i+1=n-1$, $\ldots$, $1$,  are related by the sequence of rotations $\cR_3(\psi_i+\pi)\cR_1(i^*_i):=\cR_3(\psi_i^*)\cR_1(i^*_i)$, with $\cos i^*_i=\frac{\Psi_{i}^2-G_{i+1}^2+\Psi_{i-1}^2}{2\Psi_{i-1} \Psi_{i}}$  (see figure \ref{Deprit2}, noticing that $\mathbf S_{i+1}\times \mathbf C_{i+1}=-\mathbf S_{i+1}\times \mathbf S_{i}$).
\end{itemize}
Then we find that \equ{cartesianToDeprit} has the expression
$$\left\{
\begin{array}{lll}
\dst\by_i=\cR_i^n \by^*_i\\\\
\dst\bx_i=\cR_i^n \bx^*_i
\end{array}
\right.$$
with
$$\cR_i^n:=\cR_3(\psi_{n})\cR_1(i_{n})\cR_3(\psi^*_{n-1})\cR_1(i^*_{n-1})\cdots \cR_3(\psi^*_{i})\cR_1(i^*_{i})\cR_3(\psi_{i-1})\cR_1(i_{i-1})$$
and $\by^*_i$, $\bx^*_i$ as in \equ{xjyjpl}.
}
\subsection{The map ${\cal K}$}\label{sec: K map}

\begin{figure}[htp]
\begin{tikzpicture}
\draw [ultra thick, ->] (0,0,0) -- (4,0,0);
\node (j) at (4.7,0,0) {${\mathbf S}_{j}\times \hat{\bm\nu}_j$};
\draw [ultra thick,->] (0,0,0) -- (0,4,0);
\node (k) at (0,4.2,0) {${\mathbf S}_{j}$};
\node (knew) at (0,2.2,0) {};
\node (xnew) at (2.2, 2.2, 2.2) {};
\node (incli) at (1.5, 3.2, 2.2) {$\tred{{\rm i}_{j}}$};
\draw [-latex, bend right] (xnew) to (knew);
\draw [ultra thick,->] (0,0,0) -- (0,0,5);
\node (i) at (0,0,5.5) {$\hat{\bm\nu}_j$};
\node (inew) at (0,0,1.8) {};
\draw [blue, ultra thick, ->] (0,0,0) -- (3.5,3.5,3.5);
\draw [dashed, -] (3.5,3.5,3.5) -- (4.3,4.3,4.3);
\draw [dashed, -] (4.3,4.3,4.3) -- (5.5,2.8,3.5);
\node (C) at (3.0,3.5,3.8) {$\tblue{\mathbf x_j}$};
\draw [thick, ->] (3.5,3.5,3.5) -- (5.5,2.8,3.5);
\node (y) at (4.4,3.0,3.8) {$\mathbf y_j$};
\node (Cnew) at (1.9,2.2,1.9) {$\tred{r_j}$};
\node (R) at (3.7,4.2,4) {$\tred{R_j}$};
\draw [dashed, -] (3.5,3.5,3.5) -- (3.5,0,3.5);
\draw [dashed, -] (0,0,0) -- (3.5,0,3.5);
\draw [dashed, -] (0,3.5,0) -- (3.5,3.5,3.5);
\draw [-] (-3.0,0,3.0) -- (0,0,0);
\draw [blue, ultra thick,->] (0,0,0) -- (3.0,0,-3.0);
\node (gamma) at (3.2,0,-2.0) {$\tblue{\hat{\mathbf{ n}}_j}
$};
\node (gammanew) at (0.8,0,-0.8) {} ;
\draw [-latex, bend right] (inew) to (gammanew);
\node (Zeta) at (-0.2,2,0.2) {} ;
\node (zeta) at (1.4,0,1.2) {$\tred{\hat\k_{j-1}}$} ;
\node (gammanewnew) at (1.6,0,-1.6) {} ;
\node (inewnew) at (0,0,2.4) {};
\node (varphi) at (2.2,0,4.2) {$\tred{\hat\k_{j-1}-\frac{\p}{2}}$};
\node (gammanewnewnew) at (2.2,-1.2,-1.6) {$\tred{\frac{\pi}{2}}$} ;
\node (x2new) at (3.0,0,3.0) 
{};
\draw [-latex, bend right] (inewnew) to (x2new);
\draw [-latex, bend right] (x2new) to (gammanewnew);

 \end{tikzpicture}
 \caption{The reference frames $\hat{\rm F}_j$ and the $\cK$--coordinates $\hat\k_{j-1}$, $r_j$, $R_j$  $j=2$, $\ldots$, $n$.}\label{K1}
 \begin{tikzpicture}
\draw [ultra thick, ->] (0,0,0) -- (4,0,0);
\node (j) at (4.7,0,0) {$\mathbf C_1\times \hat{\bm\nu}_1$};
\draw [ultra thick,->] (0,0,0) -- (0,4,0);
\node (k) at (0,4.2,0) {$\mathbf C_1$};
\draw [ultra thick,->] (0,0,0) -- (0,0,5);
\node (i) at (0,0,5.5) {$\hat{\bm\nu}_1$};
\node (inew) at (0,0,1.8) {};
\draw [blue, ultra thick, ->] (0,0,0) -- (3.57,0,2.75) 
;
\node (x2) at (4.0,0,3.0) 
{$\tblue{\mathbf x_1}$};
\node (x2newnew) at (2.0,0,1.5) 
{};
\node (ell2) at (2.0,0,1.5) 
{};
\node (g1) at (2,0,2.7) {};
\node (ell3) at (2.5,0,2.7) {
};

\draw [-] (-3.0,0,3.0) -- (0,0,0);
\draw [blue, ultra thick,->] (0,0,0) -- (3.0,0,-3.0);
\node (gamma) at (1.0,0,-1.0) {};
\node (gammaNEW) at (3.2,0,-3.2) {$\tblue{\hat{\mathbf n}_1}$};
\node (gammanew) at (0.8,0,-0.8) {} ;
\node (gammanewnew) at (1.6,0,-1.6) {} ;
\node (inewnew) at (0,0,2.4) {};
\node (varphi) at (2.2,0,4.2) {$\tred{\hat\vartheta_1-\frac{\p}{2}}$};
\node (gammanewnewnew) at (2.5,0,-0.8) {$\tred{\frac{\pi}{2}}$} ;
\node (x2new) at (3.0,0,2.25) 
{};
\draw [-latex, bend right] (inewnew) to (x2new);
\draw [-latex, bend right] (x2new) to (gammanewnew);
\node (Zeta) at (-0.2,2,0.2) {$\tred{\hat\Theta_1}$} ;
\node (zeta) at (1.2,0,1.2) {} ;
 \end{tikzpicture}
 \caption{The reference $\hat{\rm F}_1$ and the $\cK$--coordinates $\hat\Theta_1$, $\hat\vartheta_1$.}\label{K2}
 \begin{tikzpicture}
\draw [ultra thick, ->] (0,0,0) -- (4,0,0);
\node (j) at (4.7,0,0) {$\mathbf x_j\times \hat{\mathbf n}_j$};
\draw [ultra thick,->] (0,0,0) -- (0,4,0);
\node (k) at (0,4.2,0) {$\mathbf x_j$};
\draw [ultra thick,->] (0,0,0) -- (0,0,5);
\node (i) at (0,0,5.5) {$\hat{\mathbf n}_j$};
\node (inew) at (0,0,1.8) {};
\draw [blue, ultra thick, ->] (0,0,0) -- (3.5,3.5,3.5);
\node (C) at (3.8,3.8,3.8) {$\tblue{\mathbf S_{j-1}}$};
\node (Cnew) at (1.7,2.2,1.9) {$\tred {\hat\chi_{j-2}}$};
\draw [dashed, -] (3.5,3.5,3.5) -- (3.5,0,3.5);
\draw [dashed, -] (0,0,0) -- (3.5,0,3.5);
\draw [dashed, -] (0,3.5,0) -- (3.5,3.5,3.5);
\draw [-] (-3.0,0,3.0) -- (0,0,0);
\draw [blue, ultra thick,->] (0,0,0) -- (3.0,0,-3.0);
\node (gamma) at (3.4,0,-3.2) {$\tblue{\hat{\bm\nu}_{j-1}}$};
\node (gammanew) at (0.8,0,-0.8) {} ;
\draw [-latex, bend right] (inew) to (gammanew);
\node (Zeta) at (-0.4,2,0.2) {$\tred{\hat\Theta_j}$} ;
\node (zeta) at (1.2,0,1.2) {$\tred{\hat\vartheta_j}$} ;
%
%
\node (knew) at (0,2.2,0) {};
\node (xnew) at (2.2, 2.2, 2.2) {};
\node (incli) at (1.5, 3.2, 2.2) {$\tred{\iota_{j-1}}$};
\draw [-latex, bend right] (xnew) to (knew);
%
%
\node (gammanewnew) at (1.6,0,-1.6) {} ;
\node (inewnew) at (0,0,2.4) {};
\node (varphi) at (2.2,0,4.2) {$\tred{\hat\vartheta_{j}-\frac{\p}{2}}$};
\node (gammanewnewnew) at (2.2,-1.2,-1.6) {$\tred{\frac{\pi}{2}}$} ;
\node (x2new) at (3.0,0,3.0) 
{};
\draw [-latex, bend right] (inewnew) to (x2new);
\draw [-latex, bend right] (x2new) to (gammanewnew);

 \end{tikzpicture}
 \caption{The reference frames $\hat{\rm G}_{j}$ and the $\cK$--coordinates $\hat\vartheta_j$, $\hat\Theta_j$, $\hat\chi_{j-1}$, $j=1$, $\ldots$, $n$. When $j=2$, take $\hat\chi_0:=\hat\Theta_1$; when $j=1$, disregard ${\mathbf S}_0$, $\hat{\bm\nu}_0$ and $\hat\chi_{-1}$.}\label{K3}
 \end{figure}

 \noindent
 The $\cK$--coordinates have been described in \cite{pinzari13} for $n=2$ and generalized to any $n\in \natural$, $n\ge 2$ in \cite{pinzari18}. Here, for sake of uniformity with the coordinates $\cD_{ep}$, we change\footnote{ The main changes regard the coordinates that in \cite{pinzari18} are called $\tilde\Theta_0$, $\tilde\chi_{n-1}$, which here al called $\hat\chi_n$, $\hat\Theta_1$. The other coordinates just underwent a different numbering:  $(\tilde\Theta_j)_{1\le j\le n-1}$, $\tilde\chi_0$, $(\tilde\chi_j)_{1\le j\le n-2}$,  $\Lambda_j$ here are denoted, respectively, as $(\hat\Theta_{n-j+1})$, $\hat\chi_{n-1}$, $(\hat\chi_{n-j-1})$, $\hat\L_{n-j+1}$. An analogue change of notations holds of course for the conjugated coordinates.    } notations a little bit compared to \cite{pinzari18}. We let
$$\cK=(\hat{\bm\Theta}, \hat{\bm\chi}, {\mathbf R}, \hat{\bm\vartheta}, \hat{\bm\k}, {\mathbf r})$$
where $\mathbf R$, $\mathbf r$ are as in \equ{Deprit variables}, while
\beqano 
\begin{array}{lll}
\dst\hat{\bm\Theta}=(\hat\Theta_1,\ldots,\hat\Theta_{n}),\quad &\hat{\bm\vartheta}=(\hat\vartheta_1,\ldots,\hat\vartheta_{n})\\\\
\dst\hat{\bm\chi}=(\hat\chi_1,\ldots,\hat\chi_{n}),&\hat{\bm\k}=(\hat\k_1,\ldots,\hat\k_{n})
\end{array}
 \eeqano
are defined as follows. Let $\mathbf S_j$ be as in \equ{partial sums}.
Define the {\it$\cK$-nodes} \beq{good nodesNEW} \hat{\bm\n}_j:=\left\{
    \begin{array}
        {llll} \dst {\mathbf k}\times {\mathbf C}\quad &j=n\\
        \\
        \dst {\mathbf x}_{j+1}\times {\mathbf S}_{j} &j=1,\ldots, n-1
    \end{array}
    \right.\qquad \hat{\mathbf n}_j:=\dst{\mathbf S}_{j}\times {\mathbf x}_{j} \qquad j=1,\ldots, n. \eeq    
  and then the $\cK$coordinates  as follows.  
    
    \beqa{belle*NEW}
    \begin{array}
        {llllrrr} \dst \hat\Theta_{j}:=\left\{
        \begin{array}
            {lrrr} \dst {\mathbf S}_{j}\cdot \frac{{\mathbf x}_{j}}{|{\mathbf x}_{j}|} \\
            \\
            \dst |{\mathbf C}_{1}|
        \end{array}
        \right.& \hat\vartheta_{j}:=\left\{
        \begin{array}
            {lrrr} \dst\a_{{\mathbf x}_{j}}(\hat{\mathbf n}_{j}, \hat{\bm\nu}_{j-1})\qquad& 2\le j\le n\\
            \\
            \dst \a_{{\mathbf C}_{1}}(\hat{\bm \nu}_{1}, \hat{\mathbf n}_1)&j=1
        \end{array}
        \right.\\
        \\
        \dst\hat\chi_{j}:=\left\{
        \begin{array}
            {lrrr}Z:={\mathbf C}\cdot {\mathbf k}\\
            \\
           C:=|{\mathbf C}
            |\\
            \\
            |{\mathbf S}_{j+1}|
        \end{array}
        \right. & \hat{\k}_{j}:=\left\{
        \begin{array}
            {lrrr}
            \zeta:=\a_{{\mathbf k}}({\mathbf i}, \hat{\bm\n}_n)
       \qquad& j=n\\
            \\
           \gamma:= \a_{{\mathbf S}_{n}}(\hat{\bm\n}_{n}, \hat{\mathbf n}_{n})&j=n-1\\
            \\
            \a_{{\mathbf S}_{j+1}}(\hat{\bm\n}_{j+1}, \hat{\mathbf n}_{j+1})&1\le j\le n-2\ (n\ge 3)

        \end{array}
        \right. 
    \end{array}
    \eeqa
    \begin{remark}\rm
Note that the node $\hat{\bm\n}_n$ coincides with $\ovl{\bm\n}={\bm\n}_{n+1}$ in \equ{nodes};  the coordinates $Z$ and $\zeta$ are the same as in \equ{Deprit variables} and, finally, the coordinates $\bm\chi$ coincide with the coordinates $\hat{\bm\Psi}$ in \equ{Deprit variables}. In particular, $\cD_{ep}$ and $\cK$ share the construction in Figure \ref{Deprit1}. The geometrical meaning of the other $\cK$--coordinates is pointed out in the next section.
\end{remark}

\paragraph{\bf A chain of reference frames}
We consider the following chain of vectors {\small\beqa{chainNEW}
\begin{array}
    {cccccccccccccccccc} \dst {\mathbf k}&\to&{\mathbf S}_n={\mathbf C}&\to& {\mathbf x}_n&\to&\cdots&\to&{\mathbf S}_j&\to&{\mathbf x}_j&\to& {\mathbf S}_{j-1}&\to&\cdots&\to& {\mathbf S}_1={\mathbf C}_1\\
    \\
    \dst&&\Downarrow&&\Downarrow&&\vdots&&\Downarrow&&\Downarrow&&\Downarrow&&\vdots&&\Downarrow\\
    \\
    \dst&&\hat{\bm\n}_n&&\hat{\mathbf n}_1&&\vdots&&\hat{\bm\n}_j&&\hat{\mathbf n}_j&&\hat{\bm\n}_{j-1}&&\vdots&&\hat{\bm\n}_1\\
\end{array}
\eeqa}%
where $\hat{\bm\n}_j$, $\hat{\mathbf n}_j$ are the {\it $\cK$-nodes} in \equ{good nodesNEW}, given by the skew-product of the two consecutive vectors in the chain.

\nl We associate to this chain of vectors the following chain of frames
 \beqa{P chain}
\begin{array}
    {cccccccccccccccccc} \dst \hat{\rm G}_{n+1}&\to&\hat{\rm F}_n&\to& \hat{\rm G}_{n}&\to&\cdots&\to& \hat{\rm F}_j&\to& \hat{\rm G}_j&\to& \hat{\rm F}_{j-1}&\to&\cdots&\to &\hat{\rm G}_1
\end{array}
\eeqa where $\hat{\rm G}_{n+1}=({\mathbf i}, {\mathbf j}, {\mathbf k})$ is the initial prefixed frame and the frames, while $\hat{\rm F}_j$, $\hat{\rm G}_j$ are frames defined via \beqa{FGNEW}&&\hat{\rm F}_j=(\hat{\bm\n}_j, \ \cdot , {\mathbf S}_j)\quad\hat{\rm G}_j=(\hat{\mathbf n}_j,\ \cdot, {\mathbf x}_j)\qquad j=1,\cdots, n. \eeqa By construction, each frame in the chain has its first axis coinciding with the intersection of  horizontal plane with the horizontal plane of the previous frame (hence, in particular, $\hat{\bm\n}_j\perp {\mathbf S}_j$ and $ \hat{\mathbf n}_j\perp {\mathbf x}_j$). 
\paragraph{\bf Explicit expression of the $\cK$--map} We now derive the explicit formulae of the map which relates the coordinates \equ{belle*NEW} to the coordinates $(\by_1, \ldots, \by_n, \bx_1, \ldots, \bx_n)$. We shall prove that such map has the expression
\beqa{Kmap}\left\{
\begin{array}{ll}\dst \bx_j=\bx^n_j:={\cal R}_j^n\tilde \bx_j\\\\
\dst \by_j=\by^n_j:={\cal R}_j^n\tilde \by_j
\end{array}
\right.\eeqa where \beqa{Rr}
\left\{\begin{array}{lll}
\dst{\cal R}^n_j:=\hat{\cal T}_n\hat{\cal S}_n\cdots \hat{\cal T}_{j+1}\hat{\cal S}_{j+1}\hat{\cal T}_{j}\hat{\cal S}_{j}\\\\
 \dst\tilde \bx_j:=r_j {\mathbf k}\\\\
\dst\tilde \by_j:=R_j{\mathbf k}+\frac{1}{r_j}\tilde {\mathbf C}_j\times {\mathbf k}\\\\
\dst \tilde {\mathbf C}_j:=\left\{\begin{array}{ll}\dst\hat{\cal S}^{-1}_j \Big(\hat\chi_{j-1}{\mathbf k}-\hat\chi_{j-2}\hat{\cal S}_{j}\hat{\cal T}_{j-1}{\mathbf k}\Big)=\tilde \bx_j\times \tilde \by_j\quad &i=2, \ldots, n\\\\
\dst \hat\Theta_1\hat{\cal S}^{-1}_1{\mathbf k} &j=1
\end{array}
\right.
\end{array}
\right. \eeqa 
where
$\hat\cT_j$, $\hat\cS_j$ have  the expressions 
\beqa{CC} && \hat{\cal T}_{j}:=\left\{\begin{array}{ll}\cR_3(\zeta)\cR_1(\iota_n)\quad &j=n\\\\
\cR_3(\hat\vartheta_{j+1})\cR_1(\iota_{j})&1\le j\le n-1
	\end{array}
\right.
\qquad \hat{\cal S}_{j}:=\left\{
\begin{array}{lll}\dst\cR_3(\hat\k_{j-1})\cR_1({\rm i}_j),\quad &2\le j\le n\\\\
\dst\cR_3(\hat\vartheta_1)\cR_1(\frac{\p}{2}),\quad &j=1
\end{array}
\right. \eeqa
with
\beqa{good incli*}\dst&&\left\{\begin{array}{lll}\dst\cos\iota_{n}=\frac{Z}{\hat\chi_{n-1}}\quad &\\\\
\dst\cos\iota_{j}=\frac{\hat\Theta_{j+1}}{\hat\chi_{j-1}}& 2\le j\le n-1\ (n\ge 3)\\\\
\dst\cos\iota_{1}=\frac{\hat\Theta_{2}}{\hat\Theta_{1}}& 
\end{array}
\right.\nonumber\\
&&\left\{\begin{array}{lll}\dst \cos{\rm i}_j:=\frac{\hat\Theta_{j}}{\hat\chi_{j-1}},\quad2\le j\le n\\\\
\dst {\rm i}_1=\frac{\p}{2}\end{array}
\right.\eeqa

\vskip.1in \noindent
Indeed, $\hat{\cal T}_j$ is the rotation matrix which describes the change of coordinates from $\hat{\rm G}_{j+1}$ to $\hat{\rm F}_j$, while  $\hat{\cal S}_j$ describes the change of coordinates from $\hat{\rm F}_j$ to $\hat{\rm G}_{j}$, as it follows from the definitions of $(\hat\Theta, \hat\chi,\hat\vartheta, \hat\k)$ in \equ{belle*NEW} (see also Figures \ref{K1}, \ref{K2} and \ref{K3}). 
The formulae \equ{Kmap}--\equ{good incli*} are obtained 
considering the following sequence of transformations \beqano
\begin{array}
    {cccccccccccccccccc} &\hat{\cal T}_n&&\hat{\cal S}_n& &&\cdots&&&\hat{\cal S}_j&&\hat{\cal T}_{j-1}&& &\cdots&\hat{\cal S}_1& \\
    \\
    \dst \hat{\rm G}_{n+1}&\to&\hat{\rm F}_n&\to& \hat{\rm G}_{n}&\to&\cdots&\to& \hat{\rm F}_j&\to& \hat{\rm G}_{j}&\to& \hat{\rm F}_{j-1}&\to&\cdots&\to &\hat{\rm G}_1
\end{array}
\eeqano connecting $\hat{\rm G}_n$ to any other frame in the chain. From this, and the definitions of the frames \equ{FGNEW}, one finds $${\mathbf S}_j=\left\{\begin{array}{lll}
\dst\chi_{j-1}\hat{\cal T}_n\hat{\cal S}_n\cdots \hat{\cal T}_{j+1}\hat{\cal S}_{j+1}\hat{\cal T}_j {\mathbf k}\quad &j=2\,\ldots\,n\\\\
\hat\Theta_1\hat{\cal T}_n\hat{\cal S}_n\cdots \hat{\cal T}_{2}\hat{\cal S}_{2}{\cal T}_1{\mathbf k}&j=1
\end{array}
\right.\qquad {\mathbf x}_j=r_j\hat{\cal T}_n\hat{\cal S}_n\cdots \hat{\cal T}_{j+1}\hat{\cal S}_{j+1}\hat{\cal T}_{j}\hat{\cal S}_{j}{\mathbf k}$$ whence 
$${\mathbf C}_j=\left\{\begin{array}{lll}\dst{\mathbf S}_j-{\mathbf S}_{j-1}=\hat{\cal T}_n\hat{\cal S}_n\cdots \hat{\cal T}_{j+1}\hat{\cal S}_{j+1}\hat{\cal T}_j \Big(\hat\chi_{j-1}{\mathbf k}-\hat\chi_{j-2}\hat{\cal S}_{j}\hat{\cal T}_{j-1}{\mathbf k}\Big)\quad &j=2\,,\ \ldots\,, n\\\\
\dst{\mathbf S}_1=\hat\Theta_1\hat{\cal T}_n\hat{\cal S}_n\cdots \hat{\cal T}_{2}\hat{\cal S}_{2}\hat{\cal T}_1{\mathbf k}&j=1
\end{array}
\right.$$
and finally
\beqano
\by_j=\frac{{R}_j}{{r}_j}\bx_j+\frac{1}{r_j^2}{\mathbf C}_j\times \bx_j
\eeqano
 Collecting such formulae, one finds \equ{Kmap}--\equ{good incli*}.

\paragraph{\bf Canonical character of $\cK$}
\begin{lemma}\label{Kcanonical}
    \label{lem: good change} $\cK$ preserves the standard Liouville 1-form: \beq{1form}\sum_{j=1}^n\by_j\cdot d\bx_j=\hat{\bm\Theta}\cdot d\hat{\bm\vartheta}+\hat{\bm\chi}\cdot d\hat{\bm\k}+{\mathbf R}\cdot d{\mathbf r}.\eeq
\end{lemma}
\nl The proof of Lemma \ref{lem: good change} again relies in Lemma \ref{LemmaD}.

\proof We use the expression in \equ{Kmap}. We  also define
$${\mathbf C}_j^n:={\cal R}_j^n\tilde {\mathbf C}_j\,,\qquad \ovl{\mathbf C}^n_j:=\ovl{\cal R}_j^n\tilde {\mathbf C}_j\,,\quad \ovl{\cal R}_j^n:=\hat{\cal T}_n^{-1}{\cal R}_j^n$$
Applying Lemma \ref{LemmaD} twice, we get
$$\by^n_j\cdot d\bx^n_j={\mathbf C}^n_j\cdot {\mathbf k}\,d\zeta+\ovl{\mathbf C}^n_j\cdot {\mathbf i}\,d\iota_n+\ovl{\mathbf C}^n_j\cdot {\mathbf k}\, d\hat\k_{n-1}+{\mathbf C}^{n-1}_j\cdot {\mathbf i}\, d {\rm i}_n+\by^{n-1}_j\cdot d\bx^{n-1}_j\,.$$
Continuing in this way, after  $n-j+1$ iterates we arrive at
\beqa{yjdxj}\by_j\cdot d\bx_j&=&{\mathbf C}^n_j\cdot {\mathbf k}\,d\zeta+\ovl{\mathbf C}^n_j\cdot {\mathbf i}\,d\iota_n+\ovl{\mathbf C}^n_j\cdot {\mathbf k}\, d\hat\k_{n-1}+{\mathbf C}^{n-1}_j\cdot {\mathbf i}\, d {\rm i}_n\nonumber\\
&+&\sum_{k=j}^{n-1}\Big({\mathbf C}^k_j\cdot {\mathbf k}\,d\hat\vartheta_{k+1}+\ovl{\mathbf C}^k_j\cdot {\mathbf i}\,d\iota_k+\ovl{\mathbf C}^k_j\cdot {\mathbf k}\, d\hat\k_{k-1}+{\mathbf C}^{k-1}_j\cdot {\mathbf i}\, d {\rm i}_k\Big)\nonumber\\
&+&\widetilde\by_j\cdot d\widetilde\bx_j\eeqa
with
$${\rm i}_1:=\frac{\p}{2}\,,\ \k_{0}:=\hat\vartheta_1\,,\quad {\mathbf C}^{j-1}_j:=\tilde{\mathbf C}_j=\widetilde\bx_j\times\widetilde\by_j\,.$$
We take the sum of \equ{yjdxj} with $j=1$, $\ldots$, $n$.
 Exchanging the sums
 $$\sum_{j=1}^n\sum_{k=j}^{n-1}=\sum_{k=1}^{n-1}\sum_{j=1}^k$$
 and recognizing that
 $$\left\{
 \begin{array}{lll}
 \dst\sum_{j=1}^k{\mathbf C}^k_j=\left\{
 \begin{array}{lll}\hat{\cal S}_{k+1}^{-1}\hat{\cal T}_{k+1}^{-1}\cdots\hat{\cal S}_n^{-1}\hat{\cal T}_n^{-1}{\mathbf S}_k=\chi_{k-1}\hat{\cal T}_k{\mathbf k} \quad &1\le k\le n-1\\\\
 \dst {\mathbf S}_n=\chi_{n-1}\hat{\cal T}_n{\mathbf k} &k=n
  \end{array}
 \right.\\\\
 \dst\sum_{j=1}^k{\mathbf C}^{k-1}_j=\left\{
 \begin{array}{lll}\hat{\cal S}_{k}^{-1}\hat{\cal T}_{k}^{-1}\cdots\hat{\cal S}_n^{-1}\hat{\cal T}_n^{-1}{\mathbf S}_k=\chi_{k-1}\hat{\cal S}_{k}^{-1}{\mathbf k} \quad &1\le k\le n-1\\\\
 \dst  \hat{\cal S}_n^{-1} \hat{\cal T}_n^{-1}{\mathbf S}_n=\chi_{n-1}\hat{\cal S}_{n}^{-1}{\mathbf k}  &k=n
  \end{array}
 \right.\\\\
 \dst
 \sum_{j=1}^k\ovl{\mathbf C}^k_j=\left\{
 \begin{array}{lll}\hat{\cal T}_k^{-1}\hat{\cal S}_{k+1}^{-1}\hat{\cal T}_{k+1}^{-1}\cdots\hat{\cal S}_n^{-1}\hat{\cal T}_n^{-1}{\mathbf S}_k=\chi_{k-1}{\mathbf k} \quad &1\le k\le n-1\\\\
 \dst \hat{\cal T}_n^{-1} {\mathbf S}_n=\chi_{n-1}{\mathbf k} &k=n
  \end{array}
 \right.
 %
 %
 %
 %
 \end{array}
 \right.$$
 with $\hat\chi_0:=\hat\Theta_1$ and that, by
 \equ{Rr}, the last term in \equ{yjdxj} is  $$\tilde \by_j\cdot d\tilde \bx_j=R_j dr_j$$ 
 we get
\beqano
\sum_{j=1}^n\by_j\cdot d\bx_j&=&\sum_{j=1}^{n}\Big({\mathbf C}^n_j\cdot {\mathbf k}\,d\zeta+\ovl{\mathbf C}^n_j\cdot {\mathbf i}\,d\iota_n+\ovl{\mathbf C}^n_j\cdot {\mathbf k}\, d\hat\k_{n-1}+{\mathbf C}^{n-1}_j\cdot {\mathbf i}\, d {\rm i}_n\Big)\nonumber\\
&+&\sum_{k=1}^{n-1}\sum_{j=1}^k\Big({\mathbf C}^k_j\cdot {\mathbf k}\,d\hat\vartheta_{k+1}+\ovl{\mathbf C}^k_j\cdot {\mathbf i}\,d\iota_k+\ovl{\mathbf C}^k_j\cdot {\mathbf k}\, d\hat\k_{k-1}+{\mathbf C}^{k-1}_j\cdot {\mathbf i}\, d {\rm i}_k\Big)\nonumber\\
&+&\sum_{j=1}^{n}R_j dr_j\nonumber\\\nonumber\\
&=&\hat\chi_{n-1}\hat{\cal T}_{n}{\mathbf k}\cdot {\mathbf k}\,d\zeta+\hat\chi_{n-1}{\mathbf k}\cdot {\mathbf i}\,d\iota_n+\hat\chi_{n-1}{\mathbf k}\cdot {\mathbf k}\, d\hat\k_{n-1}
+\hat\chi_{n-1}{\mathbf k}
\cdot \hat{\cal S}_n{\mathbf i}\, 
d {\rm i}_n\nonumber\\
&+&\sum_{k=1}^{n-1}\Big(\hat\chi_{k-1}\hat{\cal T}_{k}{\mathbf k}\cdot {\mathbf k}\,d\hat\vartheta_{k+1}+\hat\chi_{k-1}{\mathbf k}\cdot {\mathbf i}\,d\iota_k+\hat\chi_{k-1}{\mathbf k}\cdot {\mathbf k}\, d\hat\k_{k-1}
+\hat\chi_{k-1}{\mathbf k}
\cdot \hat{\cal S}_k{\mathbf i}\, 
d {\rm i}_k\Big)\nonumber\\
&+&\sum_{j=1}^{n}R_j dr_j\nonumber\\
&=&\sum_{k=1}^{n}\hat\Theta_{k}d\hat\vartheta_k+\sum_{k=1}^{n}\hat\chi_{k}d\hat\k_{k}+\sum_{j=1}^{n}R_j dr_j
\eeqano
having used
 $$\hat{\cal T}_{k}{\mathbf k}\cdot {\mathbf k}=\cos\iota_k{=\frac{\hat\Theta_{k+1}}{\hat\chi_{k-1}}}\,\quad \hat{\cal S}_{k}{\mathbf i}\cdot {\mathbf k}=0\,,\quad {\mathbf k}\cdot {\mathbf k}=1\,,\quad {\mathbf i}\cdot {\mathbf k}=0\,.$$
In the following section, we shall use the following byproduct of Lemma \ref{Kcanonical}. Recall the coordinates $\cD_{e\ell}$ in \equ{Del}
and denote
$$\phi_{\cD_{e\ell}}^\cK:\quad \cD_{e\ell}=(\mathbf Z, \mathbf G, \mathbf R, \bm\zeta, \bm\phi, \mathbf r)\to{\cal K}=(\hat{\bm\Theta}, \hat{\bm\chi}, {\mathbf R}, \hat{\bm\vartheta}, \hat{\bm\k}, {\mathbf r})$$
Consider the family of projections
\beqa{projection1}\hat\phi_{\cD_{e\ell}}^\cK:\quad \cD_{e\ell}=(\mathbf Z, \mathbf G,  \bm\zeta, \bm\phi)\to{\cal K}=(\hat{\bm\Theta}, \hat{\bm\chi},  \hat{\bm\vartheta}, \hat{\bm\k})\eeqa
which, as it is immediate to see, is independent of $\mathbf r$ and $\mathbf R$.
\begin{lemma}\label{projection2}
The projections \equ{projection1}
verify
$${\mathbf Z}\cdot d\bm\zeta+{\mathbf G}\cdot d\bm\phi=\hat{\bm\Theta}\cdot d\hat{\bm\vartheta}+\hat{\bm\chi}\cdot d\hat{\bm\k}\quad \forall\ \mathbf r$$
\end{lemma}

\subsection{The reduction of perihelia $\cP$}\label{The reduction of perihelia}

\begin{figure}[htp]
\begin{tikzpicture}
\draw [ultra thick, ->] (0,0,0) -- (4,0,0);
\node (j) at (4.7,0,0) {${\mathbf S}_{j}\times {\bm\nu}_j$};
\draw [ultra thick,->] (0,0,0) -- (0,4,0);
\node (k) at (0,4.2,0) {${\mathbf S}_{j}$};
\node (knew) at (0,2.2,0) {};
\node (xnew) at (2.2, 2.2, 2.2) {};
\node (incli) at (1.5, 3.2, 2.2) {$\tred{{\rm i}_{j}}$};
\draw [-latex, bend right] (xnew) to (knew);
\draw [ultra thick,->] (0,0,0) -- (0,0,5);
\node (i) at (0,0,5.5) {${\bm\nu}_j$};
\node (inew) at (0,0,1.8) {};
\draw [blue, ultra thick, ->] (0,0,0) -- (3.5,3.5,3.5);
\node (C) at (3.0,3.5,3.8) {$\tblue{\mathbf P_j}$};
\node (Cnew) at (1.9,2.2,1.9) {
};
\draw [dashed, -] (3.5,3.5,3.5) -- (3.5,0,3.5);
\draw [dashed, -] (0,0,0) -- (3.5,0,3.5);
\draw [dashed, -] (0,3.5,0) -- (3.5,3.5,3.5);
\draw [-] (-3.0,0,3.0) -- (0,0,0);
\draw [blue, ultra thick,->] (0,0,0) -- (3.0,0,-3.0);
\node (gamma) at (3.2,0,-2.0) {$\tblue{{\mathbf{ n}}_j}
$};
\node (gammanew) at (0.8,0,-0.8) {} ;
\draw [-latex, bend right] (inew) to (gammanew);
\node (Zeta) at (-0.2,2,0.2) {} ;
\node (zeta) at (1.4,0,1.2) {$\tred{\k_{j-1}}$} ;
\node (gammanewnew) at (1.6,0,-1.6) {} ;
\node (inewnew) at (0,0,2.4) {};
\node (varphi) at (2.2,0,4.2) {$\tred{\k_{j-1}-\frac{\p}{2}}$};
\node (gammanewnewnew) at (2.2,-1.2,-1.6) {$\tred{\frac{\pi}{2}}$} ;
\node (x2new) at (3.0,0,3.0) 
{};
\draw [-latex, bend right] (inewnew) to (x2new);
\draw [-latex, bend right] (x2new) to (gammanewnew);

 \end{tikzpicture}
 \caption{The references ${\rm F}_j$ and the $\cP$--coordinates $\k_{j-1}$, $j=2$, $\ldots$, $n$.}\label{P1}
 \begin{tikzpicture}
\draw [ultra thick, ->] (0,0,0) -- (4,0,0);
\node (j) at (4.7,0,0) {$\mathbf C_1\times {\bm\nu}_1$};
\draw [ultra thick,->] (0,0,0) -- (0,4,0);
\node (k) at (0,4.2,0) {$\mathbf C_1$};
\draw [ultra thick,->] (0,0,0) -- (0,0,5);
\node (i) at (0,0,5.5) {${\bm\nu}_1$};
\node (inew) at (0,0,1.8) {};
\draw [blue, ultra thick, ->] (0,0,0) -- (3.57,0,2.75) 
;
\node (x2) at (4.0,0,3.0) 
{$\tblue{\mathbf P_1}$};
\node (x2newnew) at (2.0,0,1.5) 
{};
\node (ell2) at (2.0,0,1.5) 
{};
\node (g1) at (2,0,2.7) {};
\node (ell3) at (2.5,0,2.7) {
};

\draw [-] (-3.0,0,3.0) -- (0,0,0);
\draw [blue, ultra thick,->] (0,0,0) -- (3.0,0,-3.0);
\node (gamma) at (1.0,0,-1.0) {};
\node (gammaNEW) at (3.2,0,-3.2) {$\tblue{{\mathbf n}_1}$};
\node (gammanew) at (0.8,0,-0.8) {} ;
\node (gammanewnew) at (1.6,0,-1.6) {} ;
\node (inewnew) at (0,0,2.4) {};
\node (varphi) at (2.2,0,4.2) {$\tred{\vartheta_1-\frac{\p}{2}}$};
\node (gammanewnewnew) at (2.5,0,-0.8) {$\tred{\frac{\pi}{2}}$} ;
\node (x2new) at (3.0,0,2.25) 
{};
\draw [-latex, bend right] (inewnew) to (x2new);
\draw [-latex, bend right] (x2new) to (gammanewnew);
\node (Zeta) at (-0.2,2,0.2) {$\tred{\Theta_1}$} ;
\node (zeta) at (1.2,0,1.2) {} ;
 \end{tikzpicture}
 \caption{The reference ${\rm F}_1$ and the $\cP$--coordinates $\Theta_1$, $\vartheta_1$.}\label{P2}
 \begin{tikzpicture}
\draw [ultra thick, ->] (0,0,0) -- (4,0,0);
\node (j) at (4.7,0,0) {$\mathbf P_j\times {\mathbf n}_j$};
\draw [ultra thick,->] (0,0,0) -- (0,4,0);
\node (k) at (0,4.2,0) {$\mathbf P_j$};
\draw [ultra thick,->] (0,0,0) -- (0,0,5);
\node (i) at (0,0,5.5) {${\mathbf n}_j$};
\node (inew) at (0,0,1.8) {};
\draw [blue, ultra thick, ->] (0,0,0) -- (3.5,3.5,3.5);
\node (C) at (3.8,3.8,3.8) {$\tblue{\mathbf S_{j-1}}$};
\node (Cnew) at (1.7,2.2,1.9) {$\tred {\chi_{j-2}}$};
\draw [dashed, -] (3.5,3.5,3.5) -- (3.5,0,3.5);
\draw [dashed, -] (0,0,0) -- (3.5,0,3.5);
\draw [dashed, -] (0,3.5,0) -- (3.5,3.5,3.5);
\draw [-] (-3.0,0,3.0) -- (0,0,0);
\draw [blue, ultra thick,->] (0,0,0) -- (3.0,0,-3.0);
\node (gamma) at (3.4,0,-3.2) {$\tblue{{\bm\nu}_{j-1}}$};
\node (gammanew) at (0.8,0,-0.8) {} ;
\draw [-latex, bend right] (inew) to (gammanew);
\node (Zeta) at (-0.4,2,0.2) {$\tred{\Theta_j}$} ;
\node (zeta) at (1.2,0,1.2) {$\tred{\vartheta_j}$} ;
%
%
\node (knew) at (0,2.2,0) {};
\node (xnew) at (2.2, 2.2, 2.2) {};
\node (incli) at (1.5, 3.2, 2.2) {$\tred{\iota_{j-1}}$};
\draw [-latex, bend right] (xnew) to (knew);
%
%
\node (gammanewnew) at (1.6,0,-1.6) {} ;
\node (inewnew) at (0,0,2.4) {};
\node (varphi) at (2.2,0,4.2) {$\tred{\vartheta_{j}-\frac{\p}{2}}$};
\node (gammanewnewnew) at (2.2,-1.2,-1.6) {$\tred{\frac{\pi}{2}}$} ;
\node (x2new) at (3.0,0,3.0) 
{};
\draw [-latex, bend right] (inewnew) to (x2new);
\draw [-latex, bend right] (x2new) to (gammanewnew);

 \end{tikzpicture}
 \caption{The references ${\rm G}_{j}$ and the $\cP$--coordinates $\Theta_j$, $\vartheta_j$, $\chi_{j-2}$, $j=1$, $\ldots$, $n$. When $j=2$, take $\chi_0:=\Theta_1$; when $j=1$, disregard ${\mathbf S}_0$, $\bm\nu_0$ and $\chi_{-1}$.}\label{P3}
 \end{figure}

The $\cP$--coordinates have been described in \cite{pinzari18}. Here, as in the case of $\cK$,  we change\footnote{The coordinates named in  \cite{pinzari18} $\Theta_0$, $(\Theta_j)_{1\le j\le n-1}$, $\chi_0$, $(\chi_j)_{1\le j\le n-2}$, $\chi_{n-1}$, $\Lambda_j$ here are denoted, respectively, as $\chi_n$, $(\Theta_{n-j+1})$, $\chi_{n-1}$, $(\chi_{n-j-1})$, $\Theta_1$, $\L_{n-j+1}$. An analogue change of notations holds for the conjugated coordinates.    }  notations a little bit and denote them as
\beqa{Peri}\cP=(\bm\Theta, \bm\chi, \bm\L, \bm\vartheta, \bm\k, \bm\ell)\in \real^n\times \real_+^n\times\real_+^n\times \torus^n\times \torus^n\times \torus^n\eeqa
where $\bm\Lambda$, $\bm\ell$ are as in \equ{Delaunay variables}, while
\beqano 
\begin{array}{lll}
\dst\bm\Theta=(\Theta_1,\ldots,\Theta_{n}),\quad &\bm\vartheta=(\vartheta_1,\ldots,\vartheta_{n})\\\\
\dst\bm\chi=(\chi_1,\ldots,\chi_{n}),&\bm\k=(\k_1,\ldots,\k_{n})
\end{array}
 \eeqano
are defined as follows. Consider a phase space where the Kepler Hamiltonians \equ{KeplerHam} take negative values. Let ${\mathbf S}_j$ be as in \equ{partial sums} and ${\mathbf P}_j$ the perihelia of the instantaneous ellipses generated by  \equ{KeplerHam}, assuming they are not circles. The coordinates $\bm\L$, $\bm\ell$ are the same as in Delaunay, while, roughly,  $(\bm\Theta, \bm\chi,  \bm\vartheta, \bm\k)$ in \equ{Peri} are defined as the $(\hat{\bm\Theta}, \hat{\bm\chi},  \hat{\bm\vartheta}, \hat{\bm\k})$ of $\cK$, ``replacing ${\mathbf x}_j$ with  ${\mathbf P}_j$'' (see Figures \ref{P1}, \ref{P2}, \ref{P3}). Exact definitions are below.

\noindent
Define the {\it$\cP$-nodes} \beq{good nodes} \widetilde{\bm\n}_j:=\left\{
    \begin{array}
        {llll} \dst {\mathbf k}\times {\mathbf C}\quad &j=n\\
        \\
        \dst {\mathbf P}_{j+1}\times {\mathbf S}_{j} &j=1,\ldots, n-1
    \end{array}
    \right.\qquad \widetilde{\mathbf n}_j:=\dst{\mathbf S}_{j}\times {\mathbf P}_{j} \qquad j=1,\ldots, n. \eeq    
    Then the $\cP$--coordinates are
    
    \beqa{belle*}
    \begin{array}
        {llllrrr} \dst \Theta_{j}:=\left\{
        \begin{array}
            {lrrr} \dst {\mathbf S}_{j}\cdot {\mathbf P}_{j} \\
            \\
            \dst |{\mathbf C}_{1}|
        \end{array}
        \right.& \vartheta_{j}:=\left\{
        \begin{array}
            {lrrr} \dst\a_{{\mathbf P}_{j}}(\widetilde{\mathbf n}_{j}, \widetilde{\bm\nu}_{j-1})\qquad& 2\le j\le n\\
            \\
            \dst \a_{{\mathbf C}_{1}}(\widetilde{\bm \nu}_{1}, \widetilde{\mathbf n}_1)&j=1
        \end{array}
        \right.\\
        \\
        \dst\chi_{j}:=\left\{
        \begin{array}
            {lrrr}Z:={\mathbf C}\cdot {\mathbf k}\\
            \\
           C:=|{\mathbf C}
            |\\
            \\
            |{\mathbf S}_{j+1}|
        \end{array}
        \right. & {\k}_{j}:=\left\{
        \begin{array}
            {lrrr}
            \zeta:=\a_{{\mathbf k}}({\mathbf i}, \widetilde{\bm\n}_n)
       \qquad& j=n\\
            \\
           \gamma:= \a_{{\mathbf S}_{n}}(\widetilde{\bm\n}_{n}, \widetilde{\mathbf n}_{n})&j=n-1\\
            \\
            \a_{{\mathbf S}_{j+1}}(\widetilde{\bm\n}_{j+1}, \widetilde{\mathbf n}_{j+1})&1\le j\le n-2\ (n\ge 3)

        \end{array}
        \right. 
    \end{array}
    \eeqa

\noindent
To prove that \equ{Peri} are canonical, we consider the map
$$\phi_{\cD_{e\ell, aa}}^{\cP}:\quad \cD_{e\ell, aa}=(\mathbf Z, \mathbf G, \bm\L, \bm\zeta, {\mathbf g}, \bm\ell)\to{\cal P}=({\bm\Theta}, {\bm\chi}, \bm\L,  {\bm\vartheta}, {\bm\k}, \bm\ell)$$
relating action--angle Delaunay \equ{Delaa} and $\cP$ and its projection
$$\hat\phi_{\cD_{e\ell, aa}}^{\cP}:\quad \cD_{e\ell, aa}=(\mathbf Z, \mathbf G,  \bm\zeta, {\mathbf g})\to{\cal P}=({\bm\Theta}, {\bm\chi},  {\bm\vartheta}, {\bm\k})$$
which is independent of $\bm\L$, $\bm\ell$ (even though this will not be used). \begin{lemma}\label{projection4}
$\hat\phi_{\cD_{e\ell, aa}}^{\cP}$ coincides with the map $\hat\phi_{\cD_{e\ell}}^\cK$ in \equ{projection1}.
\end{lemma}

\noindent
Combining Lemmas \ref{projection2} and \ref{projection4}, we have
\begin{lemma}
The map $$\phi^{\cP}_{\cD_{e\ell, aa}}:\quad  \cD_{e\ell, aa}=(\mathbf Z, \mathbf G, \bm\L, \bm\zeta, {\mathbf g}, \bm\ell)\to {\cal P}=({\bm\Theta}, {\bm\chi}, \bm\L,  {\bm\vartheta}, {\bm\k}, \bm\ell)$$
verifies
$$
{\bm\Theta}\cdot d{\bm\vartheta}+{\bm\chi}\cdot d{\bm\k}+{\bm\L}\cdot d{\bm\ell}={\mathbf Z}\cdot d\bm\zeta+{\mathbf G}\cdot d{\mathbf g}+{\bm\L}\cdot d{\bm\ell}\,.
$$
\end{lemma}
\paragraph{\bf Explicit expression of the $\cP$--map} We now provide the explicit formulae of the map which relates the coordinates \equ{belle*} to the coordinates $(\by_1, \ldots, \by_n, \bx_1, \ldots, \bx_n)$. We shall prove that such map has the expression
\beqa{Pmap}\left\{
\begin{array}{ll}\dst \bx_j=\bx^n_j:={\cal R}_j^n\tilde \bx_j\\\\
\dst \by_j=\by^n_j:={\cal R}_j^n\tilde \by_j
\end{array}
\right.\eeqa where \beqa{RrPmap}
\left\{\begin{array}{lll}
\dst{\cal R}^n_j:= {\cal T}_n {\cal S}_n\cdots  {\cal T}_{j+1} {\cal S}_{j+1} {\cal T}_{j} {\cal S}_{j}\\\\
 \dst\tilde \bx_j:=a_j\Big((\cos\xi_j-e_j) {\mathbf k}+\sqrt{1-e_j^2}\sin\xi_j \tilde{\mathbf Q}_j\Big)
 \\\\
\dst\tilde \by_j:=\frac{\mu_j n_j a_j}{1-e_j\sin\xi_j}\Big(-\sin\xi_j {\mathbf k}+\sqrt{1-e_j^2}\cos\xi_j \tilde{\mathbf Q}_j\Big)\end{array}
\right. \eeqa 
where
$ \cT_j$, $ \cS_j$ have  the expressions 
\beqa{CCPmap} && {\cal T}_{j}:=\left\{\begin{array}{ll}\cR_3(\zeta)\cR_1(\iota_n)\quad &j=n\\\\
\cR_3( \vartheta_{j+1})\cR_1(\iota_{j})&1\le j\le n-1
	\end{array}
\right.
\qquad {\cal S}_{j}:=\left\{
\begin{array}{lll}\dst\cR_3( \k_{j-1})\cR_1({\rm i}_j),\quad &2\le j\le n\\\\
\dst\cR_3( \vartheta_1)\cR_1(\frac{\p}{2}),\quad &j=1
\end{array}
\right. \eeqa
with
\beqa{good incli*Pmap}\dst&&\left\{\begin{array}{lll}\dst\cos\iota_{n}=\frac{Z}{ \chi_{n-1}}\quad &\\\\
\dst\cos\iota_{j}=\frac{ \Theta_{j+1}}{ \chi_{j-1}}& 2\le j\le n-1\ (n\ge 3)\\\\
\dst\cos\iota_{1}=\frac{ \Theta_{2}}{ \Theta_{1}}& 
\end{array}
\right.\quad \left\{\begin{array}{lll}\dst \cos{\rm i}_j:=\frac{ \Theta_{j}}{ \chi_{j-1}},\quad2\le j\le n\\\\
\dst {\rm i}_1=\frac{\p}{2}\end{array}
\right.\eeqa
and
$$\tilde{\mathbf Q}_j=\frac{\tilde{\mathbf C}_j}{C_j}\times\mathbf k$$
with
\beqano
&&C_j=|\mathbf C_j|=\left\{\begin{array}{ll}\dst
\sqrt{\chi_{j-1}^2+\chi_{j-2}^2-2\Theta_{j}^2+2\sqrt{\chi_{j-1}^2-\Theta_{j}^2}\sqrt{\chi_{j-2}^2-\Theta_{j}^2}}\cos\vartheta_j\quad &j=2\,,\ldots\,,n\\\\
\dst\Theta_1 &j=1
\end{array}
\right.\nonumber\\
&& \tilde {\mathbf C}_j:=\left\{\begin{array}{ll}\dst {\cal S}^{-1}_j \Big( \chi_{j-1}{\mathbf k}- \chi_{j-2} {\cal S}_{j} {\cal T}_{j-1}{\mathbf k}\Big)=\tilde \bx_j\times \tilde \by_j\quad &i=2, \ldots, n\\\\
\dst  \Theta_1 {\cal S}^{-1}_1{\mathbf k} &j=1
\end{array}
\right.\nonumber\\
&&e_j=\sqrt{1-\frac{C_j^2}{\L_j^2}}
\eeqano
$a_j$ as in \equ{Delaunay variables}, $n_j=\sqrt{\frac{M_j}{a_j^3}}$ the mean motion, and $\xi_j$ the eccentric anomaly, solving
$$\xi_j-e_j\sin \xi_j=\ell_j\,.$$
These formulae are easily obtained using the well--known relations
$$\bx_j=a_j\Big((\cos\xi_j-e_j) {\mathbf P}_j+\sqrt{1-e_j^2}\sin\xi_j {\mathbf Q}_j\Big)$$
$$\by_j:=\frac{\mu_j n_j a_j}{1-e_j\sin\xi_j}\Big(-\sin\xi_j {\mathbf P}_j+\sqrt{1-e_j^2}\cos\xi_j {\mathbf Q}_j\Big)$$
with $\mathbf P_j$ the $j^{\rm th}$ perihelion and ${\mathbf Q}_j=\frac{{\mathbf C}_j}{C_j}\times\mathbf P_j$, and  the relations  which relate ${\mathbf C}_j$, ${\mathbf P}_j$, ${\mathbf Q}_j$ to $\cP$, which, similarly to how done for  $\cK$, are:
$${\mathbf C}_j=\cR^n_j\tilde {\mathbf C}_j\,,\quad {\mathbf P}_j=\cR^n_j{\mathbf k}\,,\quad {\mathbf Q}_j=\cR^n_j\tilde {\mathbf Q}_j\,.$$

 \subsection[The  behavior of $\cK$ and $\cP$ under reflections]{The  behavior of $\cK$ and $\cP$ under reflections}\label{P-map vs rotations and reflections} 
 The maps $\cK$ and $\cP$ have a nice behavior under reflections, which turns to be useful if they are applied to Hamiltonians which are reflection--invariant.

\nl We denote as \beqa{x*}
\bx^*=(x_{1}, -x_{2}, x_{3})\eeqa
the 
vector obtained from $\bx=(x_{1}, x_{2}, x_{3})$ by reflecting its second coordinate, and as
$${\cal R}_2^-\Big((\by_1,\ldots, \by_n), (\bx_1,\ldots, \bx_n)\Big):=\Big((\by^*_1,\ldots, \by^*_n), (\bx^*_1,\ldots, \bx^*_n)\Big)$$
the simultaneous reflection of the second coordinate of all the $\by_j$ and all the $\bx_j$ in the system of Cartesian coordinates $(\by, \bx)=\Big((\by_1,\ldots, \by_n), (\bx_1,\ldots, \bx_n)\Big)$. We aim to show that 
\begin{lemma}\label{reflections} Using $\cK$,  the reflection ${\cal R}_2^-$ is obtained by changing
$$\Big((\hat\Theta_2\,,\ldots\hat\Theta_n\,,Z)\,,\ (\hat\vartheta_2\,,\ldots\,, \hat\vartheta_n\,,\zeta)\Big)\to \Big((-\hat\Theta_2\,,\ldots\, -\hat\Theta_n\,,-Z)\,,\ (-\hat\vartheta_2\,,\ldots\,,-\hat\vartheta_n\,,-\zeta)\Big)$$
Similarly, using $\cP$, it is obtained by changing
$$\Big((\Theta_2\,,\ldots\Theta_n\,,Z)\,,\ (\vartheta_2\,,\ldots\,, \vartheta_n\,,\zeta)\Big)\to \Big((-\Theta_2\,,\ldots\, -\Theta_n\,,-Z)\,,\ (-\vartheta_2\,,\ldots\,,-\vartheta_n\,,-\zeta)\Big)$$
\end{lemma}
\proof We prove for $\cK$. We write \equ{x*} as
$$\bx^*=\cI_2^-\bx\qquad \cI_2^-=\left(
\begin{array}{rrr}
1&0&0\\
0&-1&0\\
0&0&1
\end{array}
\right)$$
Now use the formulae in  \equ{Kmap}--\equ{good incli*} and that
$$\cI_2^-\cR_3(\alpha)=\cR_3(-\alpha)\cI_2^-\,,\qquad \cI_2^-\cR_1(\beta)=\cR_1(\pi-\beta)\cI_2^-$$
and finally that the change
$$(\hat\Theta_2\,,\ldots\,\hat\Theta_n\,,Z)\to (-\hat\Theta_2\,,\ldots\, -\hat\Theta_n\,,-Z)$$
acts on the functions in \equ{good incli*} as
 $$(\iota_1\,,\ldots\,, \iota_n\,, {{{\rm i}}}_2\,,\ldots\, {\rm i}_n)\to (\pi-\iota_1\,,\ldots\,, \pi-\iota_n\,, \pi-{{\rm i}}_2\,,\ldots\, \pi-{\rm i}_n)\,.$$
 The proof for $\cP$ is similar. $\quad \square$
 
 \vskip.1in
 \noindent
 Lemma \ref{reflections} reflects on the Hamiltonian \equ{Helio} as well as in all Hamiltonians which are $\cR_2^-$--invariant as follows.
 \begin{lemma}
 Let $\cH(\by, \bx)$ be $\cR_2^-$--invariant.  Using the coordinates $\cK$, the manifolds
 $$\hat\Theta_j=0\,,\quad \hat\vartheta_j\in \{0\,,\p\}\quad j=2\,,\ldots\,,n\quad Z=0\,,\quad \zeta\in \{0\,,\p\}$$  are equilibria.
  Similarly, using the coordinates $\cP$, the manifolds
 $$\Theta_j=0\,,\quad \vartheta_j\in \{0\,,\p\}\quad j=2\,,\ldots\,,n\quad Z=0\,,\quad \zeta\in \{0\,,\p\}$$
 are equilibria.
 \end{lemma}
\section{Applications}
\subsection{Arnold's Theorem}\label{Arnold's Theorem}
Here we retrace the main ideas of the proof of Theorem \ref{Arnold Theorem} given in \cite{chierchiaPi11b}. Such proof uses on the coordinates \equ{Deprit coordinates}.
The first step is to switch from the coordinates \equ{Deprit coordinates} to a new set of coordinates which are well fitted with the close--to--be--integrable form of the Hamiltonian \equ{HelioNEW}. Then we modify the coordinates \equ{Deprit coordinates} to the following form
\beqa{Depaa}{\cal D}_{ep, aa}=(\bm\L, \mathbf G,\bm\Psi, \bm\ell,\bm\g,\bm\psi)\eeqa
 which we call {\it action--angle Deprit coordinates}, where  $\bm\Psi=(\Psi_1, \ldots, \Psi_n)$, $\bm\psi=(\psi_1, \ldots, \psi_n)$ are left unvaried, while $\bm\L=(\L_1, \ldots, \L_n)$, $\mathbf G=(\G_1, \ldots, \G_n)$,, $\bm\ell=(\ell_1, \ldots, \ell_n)$, $\bm\gamma=(\gamma_1, \ldots, \gamma_n)$ are obtained replacing the quadruplets
$(R_i, G_i, r_i, \varphi_i)$ with the quadruplets $(\L_i, \G_i, \ell_i, \g_i)$ (with $G_i=\G_i$), through the symplectic maps (depending on $\mu_i$, $M_i$)
$$(R_i, G_i, r_i, \varphi_i)\to (\L_i, \G_i, \ell_i, \g_i)$$
which integrate Kepler Hamiltonian \equ{KeplerHam}. This step is necessary to carry the integrable part in \equ{HelioNEW} to the form
$$h_{\textrm{\scshape k}}(\bm\L)=\sum_{1\le i\le n}\left(-\frac{\mu_i^3M_i^2}{2\L_i^2}\right)\,.$$
Recall that the new angles $\gamma_i$ provide the direction of the perihelion of the instantaneous ellipse generated by \equ{KeplerHam}, however they
 have a different meaning compared to the analogous angles
$g_i$ appearing in the set of Delaunay coordinates \equ{Delaunay variables}, as, by construction, the $\g_i$'s are measured
 {\it relatively to the nodes $\bm\n_i$} in \equ{nodes} (because the $\varphi_i$ were), while the angles $g_i$  in the Delaunay set are measured relatively to $\bar{\bm n}_i$ in \equ{barni}.\\
The $3n-2$ degrees of freedom Hamiltonian which is obtained is still singular.  Singularities appear when the coordinates are not defined and in correspondence of collisions among the planets. The latter case will be later excluded through a careful choice of the reference frame. The singularities of the coordinates appear when the some of the convex angles ({\it Deprit inclinations})\\
\beqa{Dincli} {i^*_j}:=({\mathbf S}_{j}, {\mathbf S}_{j+1})\quad j=1\,,\ldots, n\,,\quad {\mathbf S}_{n+1}:={\mathbf k}\eeqa take the values $0$ or $\pi$, because in such situations the angle $\psi_j$  is not defined (see Figures \ref{Deprit1}, \ref{Deprit2}, \ref{Deprit3}) and when the instantaneous orbits of some of the Kepler Hamiltonians  \equ{KeplerHam} is a circle, because in that case, the corresponding $\gamma_i$ is not defined. Such singularities are important from the physical point of view, because the eccentricities and the inclinations of the planets of the solar system are very small, hence the system is in a configuration pretty close to the singularity. To deal with this situation, a regularization similar to the Poincar\'e regularization \equ{Poinc reg} of Delaunay coordinates has been introduced in \cite{chierchiaPi11b}. Note that, in principle, there are $2^n$ singular configurations (corresponding to any choice of $i^*_j\in\{0\,,\pi\}$, besides $e_j=0$ for some $j$). Here we discuss the case $i_j=0$ for some $j$. Another regularization will be discussed in Section \ref{Coexistence of stable and whiskered tori}.

\paragraph{{\rm\scshape rps} coordinates and Birkhoff normal form}
The {\scshape rps} variables are given by $(\bm\L,\bm\l,{\mathbf z}):=(\bm\L,\bm\l,\bm\eta,\bm\xi,{\mathbf p},{\mathbf q})$  with (again) the $\L$'s as in \equ {Delaunay variables} and
 \beqa{reg var} 
 \ \  
 \l_i=\ell_i+\g_i+\psi_{i-1}^n\  \ &&\arr{\eta_i=\sqrt{2(\L_i-\G_i)}\ \cos\big(\g_i+\psi_{i-1}^n\big)\\
\xi_i=-\sqrt{2(\L_i-\G_i)}\ \sin\big(\g_i+\psi_{i-1}^n\big)
}\nonumber\\
\\
&&\arr{\dst p_i=\sqrt{2(\G_{i+1}+\Psi_{i-1}-\Psi_i)}\ \cos\psi_i^n\\
\dst q_i=-\sqrt{2(\G_{i+1}+\Psi_{i-1}-\Psi_i)}\ \sin\psi_i^n
\nonumber}
\eeqa
where
\beq{conv}\Psi_0:=\G_1\ ,\quad \G_{n+1}:=0\ ,\quad \psi_0:=0\ ,\quad \psi^n_i:=\sum_{i\le j\le n} \psi_j\ .\eeq  
Let   $\phi_{\textrm{\scshape rps}}$  denote the map
\beq{P map}
\phi_{\cC}^{\textrm{\scshape rps}}:\quad (\by,\bx)\to (\bm\L,\bm\l,{\mathbf z})\ .
\eeq
{\begin{remark}\rm
The coordinates \equ{reg var}  have been constructed as follows. First of all, we look for a linear and canonical transformation which replaces $\Psi_i$, $\G_i$, $\L_i$ with 
$$I_i:=\L_i-\G_i\,,\ J_i:=\G_{i+1}+\Psi_{i-1}-\Psi_i\quad \L_i\,,\ i=1\,,\ldots\,, n\,.$$
with the conventions in \equ{conv}. To find the coordinates $\alpha_i$, $\beta_i$, $\l_i$ respectively conjugated to $I_i$, $J_i$, $\L_i$ we impose the conservation of the standard 1--form:
\beqano
\sum_{i=1}^n(I_i d\alpha_i+J_i d\beta_i+\L_i d\l_i)&=&\sum_{i=1}^n((\L_i-\G_i)d\alpha_i+(\G_{i+1}+\Psi_{i-1}-\Psi_i) d\beta_i+\L_i d\l_i)\nonumber\\
&=&\sum_{i=1}^n\L_i d(\alpha_i+\l_i)+\sum_{i=1}^n\G_i d(-\alpha_i+\beta_{i-1})
+\sum_{i=1}^n\Psi_i d(-\beta_i+\beta_{i+1})
\eeqano
with $\beta_0:=0$, $\beta_{n+1}:=0$. This provides the following relations
$$\left\{
\begin{array}{lll}
\dst \alpha_i+\l_i=\ell_i\\\\
\dst -\alpha_i+\beta_{i-1}=\g_i\\\\
\dst -\beta_i+\beta_{i+1}=\psi_i
\end{array}
\right.$$
These equations may be solved recursively, and give
\beqa{new angles}\left\{
\begin{array}{lll}
\dst \l_i=\ell_i+\g_i+\psi_{i-1}^n\\\\
\dst \alpha_i=-(\g_i+\psi_{i-1}^n)\\\\
\dst \beta_i=-\psi_{i}^n
\end{array}
\right.\eeqa
Note that $\l_i$, $\alpha_i$, $\beta_i$ are in fact angles as the linear combinations at right hand sides of \equ{new angles} have integer coefficients.
As a second step, one defines
\beqa{polar}\left\{
\begin{array}{lll}
\dst \eta_i=\sqrt{2I_i}\cos\alpha_i\\\\
\dst \xi_i=\sqrt{2I_i}\sin\alpha_i
\end{array}
\right.\qquad \left\{
\begin{array}{lll}
\dst p_i=\sqrt{2J_i}\cos\beta_i\\\\
\dst q_i=\sqrt{2J_i}\sin\beta_i
\end{array}
\right.\eeqa
and obtains \equ{reg var} . The transformations \equ{polar} are well known to be canonical.
\end{remark}}

\noindent
The main point is that
\begin{lemma}[\cite{chierchiaPi11b}]
The map $\phi_\cC^{\textrm{\scshape rps}}$ can be extended to  a symplectic diffeomorphism  on a set  $\cP_{ \textrm{\sc rps}}^{6n}$ where the eccentricities $e_j$ and  and the angles $i_j^*$ in \equ{Dincli} are allowed to be  zero. In particular, \begin{itemize}
\item[\tiny \textbullet] 
 $e_j=0$ corresponds  to  the {\scshape rps} coordinates $\eta_j=0=\xi_j$; 
 \item[\tiny \textbullet] $i_j^*=0$ corresponds  to the  the {\scshape rps} coordinates  $p_j=0=q_j$. 
 \end{itemize}
\end{lemma}

\noindent
From the  definitions \equ{reg var}--\equ{conv} it follows that the variables
\beq{pn qn def}\arr{p_n=\sqrt{2(\Psi_{n-1}-\Psi_n)}\cos{\psi_n}=\sqrt{2(C-Z)}\cos{\zeta}\\ \\
 q_n=-\sqrt{2(\Psi_{n-1}-\Psi_n)}\sin{\psi_n}=-\sqrt{2(C-Z)}\sin{\zeta}}
\eeq
are integrals (as they are defined only in terms of the integral  ${\mathbf C}$), hence, cyclic for the Hamiltonian \equ{HelioNEW}.
Therefore, if 
 $\cH_{\textrm{\scshape rps}}$ denotes the planetary Hamiltonian expressed in {\scshape rps} variables, we have that
 \beq{HRPS}\cH_{\textrm{\scshape rps}}(\bm\L,\bm\l,\bar{\mathbf z}):= \cH\circ \phi^{\textrm{\scshape rps}}_\cC=  h_{\textrm{\scshape k}}(\L)+\m f_{\textrm{\scshape rps}}(\bm\L,\bm\l,\bar{\mathbf z})\eeq
where $\cH$ is as in \equ{HelioNEW} and  $\phi_{\textrm{\scshape rps}}$ as in \equ{P map} has $3n-1$ degrees of freedom, as it depends on $\bm\L,\bm\l,\bar{\mathbf z}$, where
 $$\bar {\mathbf z}=(\bm\eta, \bar{\mathbf p}, \bm\xi, \bar{\mathbf q})\quad {\rm with}\quad \bar{\mathbf p}=(p_1, \ldots, p_{n-1})
 $$
We denote as $a_i=\frac{1}{M_i}\left(\frac{\L_i}{\m_i}\right)^2$ the semi--major axis associated to $\L_i$. The next result solves the problem of the construction of the Birkhoff normal form for the Hamiltonian \equ{HelioNEW}, mentioned in Section \ref{sec: AT intro}.

 \begin{theorem}[\cite{chierchiaPi11b, chierchiaPi11c}]\label{planetary normal form} For any $s\in \natural$ there exists an open set  $\cA\subset \{a_1<\cdots<a_n\}$, a set  $\cM^ {6n-2}_\varepsilon\subseteq\cA\times\torus^n\times{\real^{4n}}$ containing the strip $\cM^ {6n-2}_0=\cA\times\torus^n\times \{0\}_{\real^{4n}}$, a positive number $\varepsilon$ and a symplectic map (``Birkhoff transformation'')
 \beq{birkhoff transf}\Phi_{\textrm{\scshape b}}:\quad (\bm\L,{\mathbf l},\bar{\mathbf{w}})\in\cM^ {6n-2}_\varepsilon \to (\bm\L,\bm\l,\bar{\mathbf z})\in\Phi_{\textrm{\scshape b}}(\cM^ {6n-2}_\varepsilon )
 \eeq
 which carries the  Hamiltonian \equ{HRPS} into 
 \beq{birkhoff planetary}
\cH_{\textrm{\scshape b}}(\bm\L,{\mathbf l},\bar{\mathbf{w}}):=\tilde\cH_{\textrm{\scshape rps}}\circ\Phi_{\textrm{\scshape b}}=h_{\textrm{\scshape k}} (\bm\L)+\m f_{\textrm{\scshape b}}(\bm\L,\mathbf l, \bar{\mathbf w})\eeq
where the average $f_{\textrm{\scshape b}}^{\rm av}(\L,w):=\int_{\torus^n}f_{\textrm{\scshape b}}d l$ is in BNF of order $s$:
\beq{fb}
f_{\textrm{\scshape b}}^{\rm av}(\bm\L,\bar{\mathbf w})=C_0+\O\cdot \mathbf r+{\rm P}_s(\mathbf r)+{\rm O}(|\bar{\mathbf{w}}|^{2s+1})\quad \bar{\mathbf{w}}:=(\mathbf u,\mathbf v)\quad r_i:=\frac{u_i^2+v_i^2}{2}\ ,
\eeq
${\rm P}_s$ being homogeneous polynomial in $r$ of order $s$, parameterized by $\L$.
Furthermore, the normal form \equ{birkhoff planetary}--\equ{fb} is  non--degenerate, in the sense that, if $s\ge 4$,
the $(2n-1)\times (2n-1)$ matrix $\tau(\bm\L)$ of the coefficients of the monomial \beqa{torsion}\sum_{i, j=1}^{2n-1}\tau(\bm\L)_{ij}r_i r_j\eeqa
with degree 2 in ${\rm P}_s(\mathbf r)$
 is non singular, for all $\bm\L\in \cA$.
\end{theorem}

\vskip.1in
\noindent
Denote by $B_\e=B_\e^{2n_2}=\{y\in\real^{2n_2}: |y|<\e\}$ the $2n_2$--ball of radius $\e$ and let
\beq{defcP}
\cP_\e:=V\times \torus^{n_1}\times B_\e\,.\eeq
The second ingredient is a KAM theorem for properly--degenerate Hamiltonian systems. This has been stated and proved (with a proof of about 100 pages) by Arnold in \cite{arnold63}, who named it the {\it Fundamental Theorem}.
 Here we present a refined version appeared in  \cite{chierchiaPi10}. \begin{theorem}[Fundamental Theorem, V.I.Arnold, 1963]	\label{FT} Let \beq{pndham}
H({\mathbf I}, \bm\f, {\mathbf p}, {\mathbf q}):=H_0(\mathbf I)+\m P({\mathbf I}, \bm\f, {\mathbf p}, {\mathbf q})\ ,\eeq  be real--analytic on $\cP_\e$ and assume 
\begin{itemize}{\it 
\item[{\bf (A1)}] $\mathbf I\in V\to \partial_{\mathbf I} H_0$ is a diffeomorphism;

\item[{\bf (A2)}] $\dst P_{\!\rm av} (\mathbf p,\mathbf q;\mathbf I)=P_0(\mathbf I) +\sum_{i=1}^{n_2}\O_i(\mathbf I)r_i+\frac{1}{2}\sum_{i,j=1}^{n_2} \b_{ij}(\mathbf I)r_i r_j+o_4$ where $\dst r_i:=\frac{p_i^2+q_i^2}{2}$ and 
$o_4/|(\mathbf p,\mathbf q)|^4\to 0$ as $(\mathbf p,\mathbf q)\to 0$;

\item[{\bf (A3)}]   The matrix $\beta(\mathbf I)=(\beta_{ij}(\mathbf I))$ is non--singular for all $\mathbf I\in V$.

}\end{itemize}
 Then, there exist positive numbers $\e_*$, $\m_*$, $C_*$ and $b$ such that, for 
\beq{A5'}
0<\e<\e_*\ , \quad 0<\m<\m_*\ ,\quad
\m<\su{C_* (\log \e^{-1})^{2b}}\ ,
\eeq
one can find a 
set $\cT\subset  \cP$  formed by the union of $H$--invariant $(n_1+n_2)$--dimensional tori,  on which the $H$--motion is analytically conjugated to  linear Diophantine quasi--periodic motions. The  set $\cT$  is of positive Liouville--Lebesgue measure and satisfies
\beq{impmeasest}
\meas \cP_\e>\meas \cT> \Big(1- C_* \Big(\sqrt{\m}\ ( \log \epsilon^{-1})^b+ \sqrt{\epsilon}\Big) \Big) \meas \cP_\e\ .
\eeq
\end{theorem}

\noindent
An application of Theorem \ref{FT} with $n_0=n$, $n_1=2n-1$ to the system in \equ{birkhoff planetary} with $s=4$ now leads to the proof of Theorem \ref{Arnold Theorem}.

\vskip.1in
\noindent\subsection{Global Kolmogorov tori}\label{Global Kolmogorov tori}

The quasi--periodic motions of Theorem \ref{Arnold Theorem}  provide  almost circular and almost planar orbits.  This is because 
the normal form of Theorem \ref{planetary normal form}  is constructed around the strip $\cM^ {6n-2}_{0}$, and the origin corresponds to zero eccentricities and zero mutual inclinations.  The question whether similar motions may exist outside such regime is therefore natural and important from the physical point of view. To this end, one has to understand that the Birkhoff normal form (assumption {\bf (A2)} of Theorem \ref{FT}) is used in the proof only to construct a  {reasonable} integrable approximation for the whole Hamiltonian, in fact given by
$$H_{\rm int}(\mathbf I, \mathbf r)=H_0(\mathbf I)+\m\left(
P_0(\mathbf I) +\sum_{i=1}^{n_2}\O_i(\mathbf I)r_i+\frac{1}{2}\sum_{i,j=1}^{n_2} \b_{ij}(\mathbf I)r_i r_j
\right)$$
Therefore, a possible construction of full dimensional quasi--periodic motions outside the small eccentricities and small inclinations regime  should start from a different integrable approximation. In this  section we describe an approach in such direction, where we look 
 at the first terms of the series expansion of the $\bm\ell$--averaged $f$ with respect to a small parameter.
The small parameter will be taken to be the inverse distance between the planets (the idea goes back to S. Harrington \cite{harrington69}). In addition, the use of the coordinates $\cP$ will allow to construct $(3n-2)$--dimensional quasi--periodic motions without singularities when the inclinations become zero.  Recall that the tori of Theorem \ref{Arnold Theorem} may be reduced to $(3n-2)$ frequencies (as shown in \cite{chierchiaPi11b}), in a 
almost  co--planar, co--centric configuration, but away from it, due to singularities.

\noindent
Here we discuss the following result.
\begin{theorem}[Global Kolmogorov tori in the planetary problem, \cite{pinzari18}]\label{Global Kolmogorov tori in the planetary problem}
Fix numbers $0<\underline e_i<\ovl e_i<0.6627\ldots$, $i=1,\cdots,n$. There exists a number ${\rm N}$ depending only on $n$ and a number $\a_0$ depending on $\underline e_i$, $\ovl e_i${, and} $n$ such that, if $\a<\a_0$, $\m\le \a^{\rm N}$, in a domain of planetary motions where the semi-major axes $a_n<a_{n-1}<\cdots<a_1$ are spaced as follows 
\beqa{asymptotics}
a_i^-\le a_i\le a_i^+\qquad {\rm with}\qquad a_{i}^\pm:= \frac{a_n^\pm}{\a^{\frac{1}{3}(2^{n+1}-2^{i+1}+i-n)}}
\eeqa
there exists a positive measure set $\cK_{\m, \a}$, the density of which in phase space can be bounded below as $${\rm dens}(\cK_{\m, \a})\ge 1-(\log\a^{-1})^{\rm p}\sqrt\a,$$ consisting of quasi-periodic motions with $3n-2$ frequencies where the planets' eccentricities $e_i$ verify $$\underline e_i\le e_i\le \ovl e_i.$$
\end{theorem}

\noindent
Let us consider a general set of coordinates ${\cal C}=(\bm\L, \bm\ell, \mathbf u, \mathbf v)$
which puts the Kepler Hamiltonians \equ{KeplerHam} into integrated form and hence carries the Hamiltonian \equ{HelioNEW} to

 $$\cH_{\cal C}(\bm\L, \bm\ell, \mathbf u, \mathbf v):=\cH\circ{\cal C}=-\sum_{j=1}^n\frac{{\mu}^3_j{M}^2_j}{2 \L_i^2}+\m f_{{\cal C}}(\bm\L,  \mathbf u, \mathbf v),$$ where  $$f_{{\cal C}}(\bm\L, \bm\ell, \mathbf u, \mathbf v):=\sum_{1\le i<j\le n}\bigg(\frac{\by_i\cdot \by_j}{m_0}-\frac{m_im_j}{|\bx_i-\bx_j|}\bigg)\circ\cC\,.$$ 
 
\nl We denote \beq{average K}\ovl{ f_{\cal C}}(\bm\L,  \mathbf u, \mathbf v):=\frac{1}{(2\p)^n}\int_{\torus^n}f_{{\cal C}}(\bm\L, \bm\ell, \mathbf u, \mathbf v)d\bm\ell, \eeq so that \beqano
\begin{array}
    {lll} \dst f_{{\cal C}}=\sum_{1\le i<j\le n}f_{{\cal C}}^{ij},\qquad\qquad\qquad &\dst \ovl{ f_{\cal C}}=\sum_{1\le i<j\le n}\ovl {f_{{\cal C}}^{ij}}\\
    \\
    \dst f_{{\cal C}}^{ij}:=\left(\frac{\by_i\cdot \by_j}{m_0}-\frac{m_im_j}{|\bx_i-\bx_j|}\right)\circ\cC,&\dst \ovl {f_{{\cal C}}^{ij}}:=\frac{1}{(2\p)^n}\int_{\torus^n}{ f_{\cal C}^{ij}}d\ell_1\cdots d\ell_n.
\end{array}
\eeqano

 \noi For such any  $\cC$ one always has, as a consequence of the motion equations of \equ{KeplerHam}, the following identities
\beqa{yi} 
&&\frac{1}{2\p}\int_\torus \frac{1}{\bx_j} d\ell_j=\frac{1}{a_j}\nonumber\\
&&\frac{1}{2\p}\int_\torus \by_j d\ell_j=\frac{\mu_j}{2\p}\int_\torus \dot\bx_j d\ell_j
= 0
\nonumber\\
&&\frac{1}{2\p}\int_\torus \frac{\bx_j}{|\bx_j|^3}d\ell_j=\frac{1}{2\p\mu_jM_j}\int_\torus \dot\by_j d\ell_j
=0
\eeqa
with $a_j$ the semi--major axes.  Consider now the average $\ovl{ f_{\cal C}}(\bm\L,  \mathbf u, \mathbf v)$ in \equ{average K} with respect to $\bm \ell$. Due to the fact that $\by_j$ has zero-average, one has that only the Newtonian part contributes to $\ovl{ f_{\cal C}}(\bm\L,  \mathbf u, \mathbf v)$: $$\ovl {f_{\cal C}}=-\sum_{1\le i<j\le n}\frac{m_im_j}{(2\p)^2}\int_{\torus^2}\frac{d\ell_id\ell_j}{|\bx_i-\bx_j|}. $$ We now consider any of the contributions to this sum

\beqa{fijC}\ovl {f_{{\cal C}}^{ij}}=-\frac{m_im_j}{(2\p)^2}\int_{\torus^2}\frac{d\ell_id\ell_j}{|\bx_i-\bx_j|}\qquad 1\le i<j\le n \eeqa
and expand {any such} terms \[ \ovl {f_{{\cal C}}^{ij}}=\ovl {f_{{\cal C}}^{ij}}^\ppo+ \ovl {f_{{\cal C}}^{ij}}^\ppu+ \ovl {f_{{\cal C}}^{ij}}^\ppd+\cdots\] where $$ \ovl {f_{{\cal C}}^{ij}}^\pph:=-\frac{m_im_j}{(2\p)^2}\int_{\torus^2}\frac{1}{h!}\frac{d^h}{d\varepsilon^h}\frac{1}{| \bx_i-\varepsilon\bx_j|}\Big|_{\varepsilon=0}d\ell_id\ell_j$$ is proportional to {$\frac{1}{a_i}(\frac{a_j}{a_i})^h$}. Then the formulae in \equ{yi} imply that the two first terms of this expansion are given by $$\ovl {f_{{\cal C}}^{ij}}^\ppo= -\frac{m_im_j}{a_{i}},\qquad \ovl {f_{{\cal C}}^{ij}}^\ppu= 0.$$

\nl Namely, whatever is the map $\cC$ that is used, the first non--trivial term  is the double average of the second order term, which is given by \[ \ovl {f_{{\cal C}}^{ij}}^\ppd(\bm\L,  \mathbf u, \mathbf v)=-\frac{m_im_j}{(2\p)^2}\int_{\torus^2}\frac{3(\bx_i\cdot \bx_j)^2-|\bx_i|^2|\bx_j|^2}{|\bx_i|^5}d\ell_id\ell_j.\]

\nl
Using Jacobi coordinates, S. Harrington noticed that
\begin{lemma}[\cite{harrington69}]\label{HarringtonL}
If $n=2$, $ \ovl {f_{{\cal J}}^{12}}^\ppd$  depends on one only angle: the perihelion argument of the inner planet, hence is integrable.
\end{lemma}

\noindent 
When $n=2$, Lemma \ref{HarringtonL} provides   
an effective good starting point to construct quasi--periodic motions without the constraint of small eccentricities and inclinations, because in that case one can take, as initial approximation,
\beqa{HarringtonH}\cH_{Harr}=-\sum_{j=1}^2\frac{{\mu}^3_j{M}^2_j}{2 \L_i^2}+\m \left(-\frac{m_1m_2}{a_{2}} +\ovl {f_{\cal J}^{12}}^\ppd(\L_1, \L_2, \G_1, \G_2, \g_1)\right)\eeqa
The motions of $\cH_{Harr}$ have indeed widely studied in the literature, after \cite{harrington69}.
When $n>2$, the argument does not seem to have an immediate extension using Deprit coordinates (which, as said, are the natural extension of Jacobi reduction). 
The generalization of \equ{HarringtonH} for such a case is
\beqano\cH_{\cD_{ep, aa}}=-\sum_{j=1}^n\frac{{\mu}^3_j{M}^2_j}{2 \L_i^2}+\m\sum_{1\le i<j\le n} \left(-\frac{m_im_j}{a_{j}} +\ovl {f_{\cal J}^{ij}}^\ppd\right)\eeqano
It turns out that, even looking at the nearest neighbors interactions \beqa{HarringtonH1}\cH_{nn}=-\sum_{j=1}^n\frac{{\mu}^3_j{M}^2_j}{2 \L_i^2}+\mu\sum_{i=1}^{n-1}\left(-\frac{m_{i}m_{i+1}}{a_{i}}+\ovl {f_{{\cal D}_{ep}}^{i, i+1}}^\ppd\right)
\eeqa
the terms $\ovl {f_{{\cal D}_{ep}}^{i, i+1}}^\ppd$ with $1\le i\le n-2$ depend on two angles: $\gamma_i$ and $\psi_{i-1}$, so the effective study of the unperturbed motions of \equ{HarringtonH1} is involved. Using the $\cP$--coordinates
\beqa{HnnPeri} \cH_{nn}=-\sum_{j=1}^n\frac{{\mu}^3_j{M}^2_j}{2 \L_i^2}+\mu\sum_{i=1}^{n-1}\left(-\frac{m_{i}m_{i+1}}{a_{i}}+\ovl {f_{{\cal P}}^{i, i+1}}^\ppd\right)
\eeqa
 one has that the terms ${\ovl {f_{\cP}^{i, i+1}}}^{(2)}$ with $1\le i\le n-2$ depend on $3$ angles: $\k_{i-1}$, $\vartheta_{i}$ and $\vartheta_{i+1}$, but the dependence upon $\k_{i-1}$ and $\vartheta_{i}$ is
at a higher order term. This is shown by the following formula, discussed in \cite{pinzari18}:
\beqa{ovl f} {\ovl {f_{\cP}^{i, i+1}}}^{(2)}&=&m_{i} m_{i+1} \frac{a_{i+1}^2}{4a_{i}^3}\frac{\L_{i}^3}{\chi_{i-1}^2(\chi_{i-1}-\chi_{i-2})^3}\Big[ \frac{5}{2}(3\Theta_{i+1}^2-\chi_{i-1}^2)\nonumber\\
    &-&\frac{3}{2}\frac{4\Theta_{i+1}^2-\chi_{i-1}^2}{\L_{i+1}^2}\Big(\chi_{i}^2+\chi_{i-1}^2-2\Theta_{i+1}^2+2\sqrt{(\chi_{i}^2-\Theta_{i+1}^2)(\chi_{i-1}^2-\Theta_{i+1}^2)}\cos{\vartheta_{i+1}}\Big)\nonumber\\
    &+&\frac{3}{2}\frac{(\chi_{i-1}^2-\Theta_{i+1}^2)(\chi_{i}^2-\Theta_{i+1}^2)}{\L_{i+1}^2}\sin^2{\vartheta_{i+1}}\nonumber\\
    &+&{\rm O}(\Theta_{i}^2+(\vartheta_{i}-\vartheta^0_i)^2)\Big] \qquad i=1\,,\ldots\,,\ n-1\eeqa
    where $\chi_0:=\Theta_1$, $ \chi_{-1}:=0$, $\vartheta^0_i\in\{0\,,\p\}$ and the ${\rm O}(\Theta_{i}^2+(\vartheta_{i}-\vartheta^0_i)^2)$ term vanishes identically when $i=1$.

\vskip.1in
\noindent
{We denote as
\beq{HP}{\rm H}_\cP({\rm X}_\cP,\ell)={\rm h}_{\rm fast}^0(\L)+\m f_\cP({\rm X}_\cP, \ell)\qquad {\rm X}_\cP:=(\Theta,\chi, \L,\vartheta, \k)\eeq
where
\beq{hk0}{\rm h}_{\rm fast}^0(\L):=-\sum_{j=1}^n\frac{\mu_j M_j^2}{2\L_j^2}\ ,\eeq
the $(3n-2)$--dimensional Hamiltonian \equ{Helio} expressed in $\cP$--coordinates.
}
The proof of Theorem \ref{Global Kolmogorov tori in the planetary problem} is based on  {three} steps: {in step 0 we compute the holomorphy domain of $\HH_\cP$}; in the  step  1 the Hamiltonian is transformed to a similar one, but with a much smaller remainder. In  step  2, a well fitted KAM theory is applied. Note that, as the terms of the unperturbed part are smaller and smaller as and when the distance from the sun increases, such KAM theory will be required to take such different scales into account.

\paragraph{Step 0: Choice of the holomorphy domain}
{
A typical practice, in order to use perturbation theory techniques, is  to extend Hamiltonians governing dynamical systems to the complex field, and then to study their holomorphy properties.\vskip.1in
\nl
It can be proven that a  domain of  holomorphy for the perturbing function $f_\cP$ in \equ{HP}, regarded as a function of complex coordinates
can be chosen as
$${\mathbb  D}_{\cP}:={\cal T}_{\Theta^+,\vartheta^+}\times\big({\cal X}_\theta\times \ovl\torus^n_{{s}}\big)\times\big({\cal A}_\theta\times \ovl\torus^{n}_{{s}}\big)\ ,$$
where, for given positive numbers 
\beqano
\Theta_{j}^+\ ,\quad \vartheta_{j}^+\ ,\quad{\rm G}_i^\pm\ ,\quad \L_i^\pm\ ,\quad \theta_i\ ,\quad s
\eeqano
with  $i=1$, $\cdots$, $n$, $j=1$, $\cdots$, $n-1$,
\beqa{domain***}
{\cal T}_{\Theta^+,\vartheta^+}&:=&\Big\{(\ovl\Theta, \ovl\vartheta)=(\Theta_2,\cdots, \Theta_{n}, \vartheta_2,\cdots,\vartheta_{n})\in \complex^{n-1}\times \torus_\complex^{n-1}:\ \nonumber\\
&& |\vartheta_j-\p|\le {\vartheta^+_j}\ ,\quad |\Theta_{j}|\le {\Theta_j^+}\ ,\quad \forall\ j=2,\cdots, n\Big\}\nonumber\\\nonumber\\
{\cal X}_\theta&:=&\Big\{
(\Theta_1, \ovl\chi)=\Big(\Theta_1, (\chi_1,\cdots, \chi_{n-1})\Big)\in \complex^n:\ {\rm G}_{j}^-\le |\chi_{j-1}-\chi_{j-2}|\le{\rm G}_{j}^+\ ,\nonumber\\
&&|\Im (\chi_{j-1}-\chi_{j-2})|\le \theta_{j}\  \forall\ j=1,\cdots, n\Big\}\nonumber\\\nonumber\\
{\cal A}_\theta&:=&\Big\{\L=(\L_1,\cdots, \L_n)\in \complex^n:\quad \L_j^-\le | \L_j|\le \L_j^+\ ,\quad |\Im \L_j|\le \theta_j\nonumber\\
&& \forall\ j=1,\cdots, n\Big\}\nonumber\\
 \ovl\torus_{{s}}&:=&\torus+{\rm i}[-{s}, {s}]\eeqa
with $\chi_{-1}:=0$, $\chi_{0}:=\Theta_1$, and
\beqa{Choice of parameters}
 \L_i^\pm&:=&\m_i\sqrt{M_ia_i^\pm}\ ,\quad {\rm G}_i^+:=\ovl\cC_i^* \L_i^-\ ,\qquad {\rm G}_i^-:=\underline\cC_i^* \L_i^+\ ,\quad\Theta_j^+:=s{\rm G}_1^-\ ,\qquad 
\vartheta_j^+:=\cD_i\frac{{\L_i^-}}{{\rm G}_1^+}\nonumber\\ 
 \theta_i&:=&s\sqrt{\L_i^-}
\eeqa
with $s\in (0, 1)$ arbitrary, $\cD_i$, $\underline\cC_i^*$, $\ovl\cC_i^*$ depending only on $m_0$, $\ldots$, $m_n$,  $a_i^\pm$ as in \equ{asymptotics}.
}
\paragraph{Step 1: Normal Form Theory}

\begin{definition}\label{def: Diophantine numbers} \rm Given $m$, $\n_1$, $\cdots$, $\n_m\in \natural$, $\n:=\n_1+\cdots+\n_m$; $\g_1$, $\ldots$, $\g_m$, $\t\in \real_+$. We call {\it m--scale Diphantine set}, and denote it as ${\cal D}_{\g_1, \ldots, \g_m, \t}$, the set of $\omega=(\o_1,\cdots,\o_m)$, with $\omega_j\in {\mathbb R}^{\n_j}$ such that,    for any $k=(k_1,\cdots,k_m)\in \integer^\n\setminus\{0\}$, with $k_j\in \integer^{\n_j}$, the following inequalities hold:
    \beq{dioph2sc} |\o\cdot k|=\bigg|\sum_{j=1}^m\o_j\cdot k_j\bigg|\geq \left\{
    \begin{array}
        {l} \dst\frac{\g_1}{|\tk|^{\t}}\quad {\rm if}\quad k_1\neq0;\\
        \ \\
        \dst \frac{\g_2}{|k|^{\t}}\quad {\rm if}\quad k_1= 0,\quad k_2\neq 0;\\
        \ \\
        \;\;\;\vdots\\
        \ \\
        \dst \frac{\g_m}{|k_m|^{\t}}\quad {\rm if}\quad k_1=\cdots=k_{m-1}= 0,\ \cdots,\ k_{m}\neq 0.
    \end{array}
    \right. \eeq
\end{definition}
The set  ${\cal D}_{\g_1, \ldots, \g_m, \t}$ reduces to the usual diophantine set taking $\g_j=\gamma$ $\forall\ j$. The first multi--scale Diophantine set was proposed by Arnold in \cite{arnold63} with $m=2$.

\begin{proposition}\label{multi scale normal form}
    \label{exponential average} Let ${  \m}_j$, ${  M}_j$ be as in \equ{masses} and ${\rm m}_j:=\sum_{i=1}^{j-1} m_i$, with $j=2,\cdots, n$, ${\chi_0:=\Theta_1}$. There exists a number ${  c}$, depending only on $n$, $m_0,\cdots,m_n$, $a_1^\pm$, $\underline e_j$, $\ovl e_j$, and a number $0<\ovl{ {c}}<1$, depending only on $n$ such that, for any fixed positive numbers $\ovl\g<1<\bar K$, $\a>0$ verifying \beqa{bar K} &&\bar K\le \frac{{  c} }{\a^{3/2}}\eeqa and \beqa{small secular}\frac{1}{{  c}}\max\Big\{ \m(\frac{a_n^+}{a_1^-})^{5}\frac{\bar K^{2\bar\t+2}}{\bar\g^2},\ \frac{\bar K^{2(\bar\t+1)}\a}{\bar\g^2}\Big\}<1\eeqa there exist natural numbers $\n_1,\cdots,\n_{2n-1}$, with $\sum_j\n_j=3n-2$, open sets $B_j^*\subset B^{2}_{\varepsilon_j}, {\cal X}^*\subset {\cal X}$, positive real numbers \mbox{$\g_1> \cdots >\g_{2n-1} \varepsilon_1, \cdots, \varepsilon_{n-1}, \ovl r_1, \cdots, \ovl r_{n-1}, \widetilde r_1, \cdots, \widetilde r_{n}$}, a domain $${  D}_{\rm n}:=B_{\sqrt{2\ovl r}}\times {\cal X}_{\ovl r}\times\cA_{\widetilde r} \times\torus^n_{\ovl{ {c}}s}\times \torus^n_{\ovl{ {c}}s}$$ a sub-domain of the form $${  D}^*_{\rm n}:=B^*_{\sqrt{2\ovl r}}\times {\cal X}^*_{\ovl r}\times\cA_{\widetilde r} \times\torus^n_{\ovl{ {c}}s}\times \torus^n_{\ovl{ {c}}s}$$ verifying \beq{good set***}\meas{  D}^*_{\rm n}\ge\big(1-\frac{\bar\g}{\ovl{  c}}\big)\meas{  D}_{\rm n}\eeq a real-analytic transformation $$\phi_{\rm n}:\quad (p,q,\chi,\L,\k,\ell)\in {  D}^*_{\rm n}\to {  D}_\cP$$ which conjugates $\cH_\cP$ to $$\cH_{\rm n}(p,q, \chi,\L,\k,\ell) :=\cH_\cP\circ\phi_{\rm n}={\rm h}_{ {fast},  {sec}}(p,q,\chi,\L)+\m\,{f}_ {exp}(p,q, \chi,\L,\k,\ell) $$ where ${f}_ {exp}(p,q, \chi,\L,\k,\ell)$ is independent of $\k_{n-1}$, and the following holds.

\paragraph{1.} The function ${\rm h}_{ {fast},  {sec}}(p,q,\chi,\L)$ is a sum $${\rm h}_{ {fast},  {sec}}(p,q,\chi,\L)={\rm h}_{ {fast}}(\L)+\m\, {\rm h}_{ {sec}}(p,q,\chi,\L)$$ where, if \beqano \hat{\rm y}_i:=\bigg(\frac{p_2^2+q_2^2}{2},\ \cdots,\ \frac{p_{i+1}^2+q_{i+1}^2}{2},\ \chi_{0},\ \cdots,\ \chi_{i},\ \L_1,\ \cdots,\ \L_{i+1}\bigg)\qquad i=1\,,\ldots\,,n-1 \eeqano then ${\rm h}_{ {fast}}$ and ${\rm h}_{ {sec}}$ are given by $${\rm h}_{ {fast}}(\L)=-\sum_{j=1}^n\frac{{  m}_j^3{  M}_j^2}{2\L_j^2}-\m\sum_{j=1}^{n-1}\frac{{  M}_j{  m}_j^2m_j{\rm m}_j}{\L_j^2} ,\qquad {\rm h}_{ {sec}}(p,q,\chi,\L)=\sum_{i=1}^{n-1}{\rm h}_{ {sec}}^i(\hat{\rm y}_i)$$ where the functions ${\rm h}_{ {sec}}^i$ have an analytic extension on ${  D}_{\rm n}$ and verify \[{  c}\frac{(a_{j+1}^+)^2}{(a_{j}^-)^3}\le | {\rm h}_{ {sec}}^j(\hat{\rm y}_j) |\le \frac{1}{  c}\frac{(a_{j+1}^+)^2}{(a_{j}^-)^3}.\]

\paragraph{2.} The function ${f}_ {exp}$ satisfies $$|{f}_ {exp}|\le \frac{1}{  c}\frac{e^{-{  c}\bar K}}{a_{n}^-}.$$ \paragraph{3.} If $\zeta$ is $\hat{\rm y}_{{n-1}}$ deprived of  ${\chi_{n-1}=C}$, the frequency-map $$\zeta\to\o_{ {fast},  {sec}}(\zeta):= \partial_{\zeta}{\rm h}_{ {fast},  {sec}}(\zeta)$$ is a diffeomorphism of $\P_\zeta(B^*_{\sqrt{2\ovl r}}\times {\cal X}^*_{\ovl r}\times\cA^*_{\widetilde r})$ and, moreover, it satisfies \equ{dioph2sc}, with $m=2n-1$, $\t=\bar\t>2$, and
 \begin{eqnarray*}
 \n_j&:=&\left\{
    \begin{array}
        {llll} \dst\ 1& j=1,\cdots, n\\
        \\
        \dst \ 2\qquad &j=3,\ n=2\\
        \\
        \dst \ 3 & j=n+1,\ n\ge 3 \\
        \\
        \dst \ 2& n+2\le j\le 2n-2,\ n\ge 4 \\
        \\
        \dst\ 1& j=2n-1,\ n\ge 3
    \end{array}
    \right.
 \end{eqnarray*}
  \beqa{nugamma}  
    \omega_j&:=&\left\{
    \begin{array}
        {llll} \dst\partial_{\L_j}{\rm h}_{ {fast},  {sec}}& j=1,\cdots, n\\
        \\
        \dst \partial_{(\frac{p_{2}^2+q_{2}^2}{2},\chi_0)}\,{\rm h}_{ {fast},  {sec}}\qquad &j=3,\ n=2\\
        \\
        \dst \partial_{(\frac{p_{2}^2+q_{2}^2}{2}, \chi_{1},\chi_{0})}\,{\rm h}_{ {fast},  {sec}}& j=n+1,\ n\ge 3 \\
        \\
        \dst \partial_{(\frac{p_{j-n+1}^2+q_{j-n+1}^2}{2}, \chi_{j-n})}\,{\rm h}_{ {fast},  {sec}}& n+2\le j\le 2n-2,\ n\ge 4 \\
        \\
        \dst \partial_{\frac{p_{n}^2+q_{n}^2}{2}}\,{\rm h}_{ {fast},  {sec}}& j=2n-1,\ n\ge 3
    \end{array}
    \right.\nonumber\\
    \nonumber\\
    \g_j&:=& \left\{
    \begin{array}
        {llll} \dst\frac{1}{a_j^-}\frac{\ovl\g}{\theta_j} \qquad &1\le j\le n\\
        \\
        \dst\frac{\m(a_{2n-j+1}^+)^2}{(a_{2n-j}^-)^3}\frac{\ovl\g}{\theta_{j-n}} &n+1\le j\le 2n-1
    \end{array}
    \right.\eeqa

\paragraph{4.} The mentioned constants are \beqano \varepsilon_j:={  c}\,\sqrt{\theta_j},\quad \ovl r_j:=\frac{\theta_j\ovl\g}{\bar K^{\bar\t+1}} ,\quad \widetilde r_i:={  c}\,\theta_j \eeqano with $\bar\t>2$.
\end{proposition}

\paragraph{Step 2: KAM theory}
\begin{theorem}
    [Multi-scale KAM Theorem, \cite{pinzari18}]\label{two scales KAM} 
    Let $m,\ell,\n_1,\cdots,\n_m\in \natural$, $\n:=\n_1+\cdots+\n_m\ge \ell$, $\t_*>\n$, $\g_1\ge \cdots\ge\g_m>0$, $0<4s\leq \bar{s}<1$, $\r_1, \cdots, \r_\ell, r_1, \cdots, r_{\n-\ell}, \varepsilon_1, \cdots, \varepsilon_\ell>0$, $B_1, \cdots, B_\ell\subset \real^2$, $D_j:=\{\frac{x^2+y^2}{2}\in \real: (x,y)\in B_j\}\subset \real$, $B:=B_1\times\cdots\times B_\ell\subset \real^{2\ell}$, $D:=D_1\times\cdots\times D_\ell\subset \real^\ell$, $C\subset \real^{\n-\ell}$, $A:=D_\r\times C_r$.
     Let \beqano {\rm H} (\mathbf p,\mathbf q, \mathbf I,\bm\psi)={\rm h}(\mathbf p,\mathbf q, \mathbf I)+{f}(\mathbf p,\mathbf q, \mathbf I,\bm\psi) \eeqano be real-analytic on $B_{\sqrt{2\r}}\times C_r\times \torus_{\bar{ s}+s}^{\n-\ell}$, where ${\rm h}(\mathbf p,\mathbf q, \mathbf I)$ depends on $(\mathbf p,\mathbf q)$ only via \[ J(\mathbf p,\mathbf q):=\Big(\frac{p_1^2+q_1^2}{2},\ \cdots,\ \frac{p_\ell^2+q_\ell^2}{2}\Big).\] Assume that $\o_0:=\partial_{J(\mathbf p,\mathbf q, \mathbf I)} {\rm h}$ is a diffeomorphism of $A$ with non singular Hessian matrix $U_{ 1}:=\partial^2_{(J(\mathbf p,\mathbf q, \mathbf I)}{\rm h}$ and let $U_k$ denote the $ (\n_k+\cdots+\n_m)\times \n$ submatrix of $U$, \ie, the matrix with entries $(U_k)_{ij}=U_{ij}$, for $\n_{1}+\cdots+\n_{k-1}+1\leq i\leq \n$, $1\leq j\leq \n$, where $2\le k\le m$. Let \beqano &&  {\rm M}_k\geq\sup_{A}|U_k|,\quad \bar {\rm M} \geq\sup_{A}|U^{-1}|,\quad \pertnorm\geq|{f}|_{\r,\bar{ s}+s}\nonumber\\
    &&\bar {\rm M}_k\geq \sup_{A}|T_k|\quad {\rm if}\quad \dst U^{-1}=\left(
    \begin{array}
        {lrr} T_1\\
        \vdots\\
        T_m
    \end{array}
    \right)\qquad 1\le k\le m.\eeqano Define \beqano && \dst K:=\frac{6}{s}\ \log_+{\left(\frac{\pertnorm {\rm M}_1^2\,L}{\gamma_1^2}\right)^{-1}}\quad {\rm where}\quad \log_+ a :=\max\{1,\log{a}\}\\
    && \dst \hat\r_k:=\frac{\g_k}{3{\rm M}_kK^{\t_*+1}},\quad \hat\r:=\min\left\{\hat\r_1,\ \cdots,\ \hat\r_m,\ \r_1,\ \cdots,\ \r_\ell,\ r_1,\ \cdots ,\ r_{\n-\ell}\right\}\\
    \\
    && \dst L:=\max \Big\{\bar {\rm M} , \ {\rm M}_1^{-1},\ \cdots,\ {\rm M}_m^{-1}\Big\} \\
    && \hat E:=\frac{E L}{\hat\r^2}. \eeqano Then one can find two numbers $\hat c_\n>c_\n$ depending only on $\n$ such that, if the perturbation ${f}$ {is} so small that the following ``KAM condition'' holds \[ \hat c_\n\KAM<1, \] for any $\o\in\O_*:=\o_0({D})\cap\cD_{\g_1,\cdots,\g_m,\t_*}$, one can find a unique real-analytic embedding \beqano \phi_\o:\quad \vartheta=(\hat\vartheta,\bar\vartheta)\in\torus^{\n } &\to&(\hat v(\vartheta;\o),\hat\vartheta+\hat u(\vartheta;\o), \cR_{\bar\vartheta+\bar u(\vartheta;\o)}w_1,\ \cdots,\ \cR_{\bar\vartheta+\bar u(\vartheta;\o)}w_\ell)\nonumber\\
    &&\in \Re C_r\times \torus^{\n-\ell}\times \Re B^{2\ell}_{\sqrt{2r}} \eeqano where $r:= c_\n \KAM \hat\r$ such that ${\rm T}_\o:=\phi_{\o}(\torus^\n)$ is a real-analytic $\n$-dimensional ${\rm H}$-invariant torus, on which the ${\rm H}$-flow is analytically conjugated to $\vartheta\to \vartheta+\o\,t$. Furthermore, the map $(\vartheta;\o)\to\phi_\o(\vartheta)$ is Lipschitz and one-to-one and the invariant set $\dst {{\rm K}}:=\bigcup_{\o\in\O_*}{\rm T}_\o$ satisfies the following measure estimate \[
	\meas\Big(\!\Re({D}_r)\times\torus^\td\setminus{{\rm K}}\Big)
        \leq  c_\n\Big(\!\meas({D}\setminus{D}_{\g_1,\cdots, \g_m,\t_*}\times\torus^\td)+\meas(\Re({D}_r)\setminus{D})\times\torus^\td\Big),
    \] where ${D}_{\g_1,\cdots, \g_m,\t_*}$ denotes the $\o_0$-pre-image of $\cD_{\g_1,\cdots, \g_m,\t_*}$ in ${D}$. Finally, on $\torus^\n\times \O_*$, the following uniform estimates hold 
    \begin{align*}
    | v_k(\cdot;\o)-I_k^0(\o)|
    &\leq c_\n \Big(\frac{\bar {\rm M}_k}{\bar {\rm M}}+\frac{{\rm M}_k}{{\rm M}_1}\Big)\KAM\,\hat\r \\
    |u(\cdot;\o)| &\leq c_\n\KAM\,s
    \end{align*}
    where $v_k$ denotes the projection of $v=(\hat v, \bar v)\in \real^{\n_1}\times\cdots\times\real^{\n_m}$ over $\real^{\n_k}$, $\dst\bar v_k:=\frac{|w_k|^2}{2}$ and $I^0(\o) = (I^0_1(\o),\cdots, I^0_\n(\o)) \in D$ is the $\o_0$-pre-image of $\o\in\O_*$.
\end{theorem}

\noindent
Theorem \ref{two scales KAM}  generalizes Theorem 3 in \cite{chierchiaPi10} and  hence the Fundamental Theorem of \cite{arnold63}, to which Theorem 3 in \cite{chierchiaPi10}  is inspired.

\paragraph{Proof of Theorem \ref{Global Kolmogorov tori in the planetary problem}} Let $$ \bar\g:={\ovl{ c}}\sqrt\a(\log\a^{-1})^{\bar\t+1},\quad \bar K=\frac{1}{\widetilde{ c}}\log\frac{1}{\a}$$ where $\ovl{ c}$ is as in \equ{good set***} and $\widetilde{ c}$ will be fixed later. We {aim to} apply Theorem \ref{two scales KAM} to the Hamiltonian $\cH_{\rm n}$ of Proposition \ref{exponential average}, with these choices of $\bar\g$ and $\bar K$. To this end, we take \beqano &&{\rm M}_j=\left\{
\begin{array}
    {llll} \dst\frac{1}{{ c}_1a_j^-\theta_j^2}\qquad\ &1\le j\le n\\
    \\
    \dst\frac{\m(a_{2n-j+1}^+)^2}{{ c}_1(a_{2n-j}^-)^3\theta_j^2} & n+1\le j\le 2n-1
\end{array}
\right. \qquad L=\bar{\rm M}=\frac{1}{{ c}_2}\,\frac{\theta_1^2(a_{n}^+)^3}{\m(a_{n-1}^-)^2} \nonumber\\
&&E=\frac{1}{{ c}_3}\frac{\m}{a_n^-}e^{-{ c}\bar K} \qquad\qquad\qquad\qquad\qquad \qquad\qquad\quad K=\frac{1}{{ c}_4}\log_+\Big(\frac{1}{\ovl\g^2}\frac{(a_n)^3}{(a_{n-1}^-)^3}e^{-{ c}\bar K}\Big)^{-1}\nonumber\\
&&\hat \r_j=\arr{\dst{{ c}_5}\frac{\ovl\g\theta_j}{K^{\t_*+1}}\quad 1\le j\le n\\
\\
\dst{ c}_5\frac{\ovl\g\theta_{j-n}}{{}K^{\t_*+1}}\quad n+1\le j\le 2n-1 }\qquad \qquad\hat\r:=\frac{\theta_1\ovl\g}{\hat K^{\t_*+1} }\quad \t_*>3n-2\nonumber\\
&&\hat E=\frac{1}{{ c}_6}\frac{1}{\ovl\g^2}\frac{(a^{ +}_n)^3}{(a_{n-1}^-)^3}e^{-{ c}\bar K}\hat K^{2(\t_*+1)}\nonumber\\
\eeqano where $\hat K:=\max\{K,\bar K\}$. The number $\frac{1}{\ovl\g^2}\frac{(a_{n-1})^3}{(a_n^-)^3}$ can be bounded by $\frac{1}{\a^N}$ for a sufficiently large $N$ depending only on $n$. Hence, if $\widetilde{ c}<\frac{ c}{N}$ and $\a<{ c}_6$, we have $\hat E<1$ and the theorem is proved. $\quad\square$

\subsection{On the co--existence of stable and whiskered tori}\label{Coexistence of stable and whiskered tori}
In this section we discuss how the use two different sets of coordinates may lead to prove the co--existence of stable and unstable motions. Specifically, we deal with the following situation, which we shall refer to as {\it outer, retrograde configuration} ({\sc orc}):

\vskip.1in
\noindent
{\it Two planets describe  almost  co--planar orbits, revolving around their common sun,  in opposite sense. The outer planet has a lower angular momentum and retrograde motion, as seen from the total angular momentum of the system.  }
\vskip.1in
\noindent
  We aim to discuss the  following
  \begin{theorem}\label{thm: coexistence}
 \item[{\rm 1.}] {\it There exists a {8--dimensional} region $\cD_{\rm s}$ in the phase space almost completely filled with a positive measure set of five--dimensional 
{\sc kam} tori, in {\sc orc} configuration;
  
  \item[{\rm 2.}] There exists a {8--dimensional}  region $\cD_{\rm u}$ in the phase space including a {6--dimensional, hyperbolic} invariant region $\cD^0_{\rm u}$ consisting of co--planar, retrograde motions for the outer planet.

    \item[{\rm 3.}]     $\cD_{\rm s}$ and $\cD^0_{\rm u}$  have a non--empty intersection}.
    
      \end{theorem}

      \noindent
Theorem \ref{thm: coexistence} leads to the following conjecture, which is likely to be proved somewhere.
\begin{conjecture}
Full dimensional quasi--periodic motions and hyperbolic 3--dimensional tori co--exist in $\cD_{\rm s}$.
  \end{conjecture}
  
  \noindent
  The proof of statements 1. and 2. in Theorem \ref{thm: coexistence} relies on the use of two different sets of coordinates for the Hamiltonian \equ{HelioNEW} with $n=2$:
\beqa{3BP}{\cal H}_{3BP}&=\,&\frac{|\by_1|^2}{2\mu_1}-\frac{\mu_1 M_1}{|\bx_1|}+\frac{|\by_2|^2}{2\mu_2}-\frac{\mu_2 M_2}{|\bx_2|}+\mu\left(\frac{\by_1\cdot \by_2}{ m_0}-\frac{m_1 m_2}{|\bx_i-\bx_j|}\right)
\eeqa

\paragraph{Proof of 1.}

We consider the coordinates \equ{Depaa} with $n=2$. {It will turn to be useful to work with regularizing complex  coordinates, which we denote as 
\beqa{real and complex}
&& \textrm{\sc rps}_\p^{\complex}:=(\bm\L,\bm\l,   \mathbf t,  \mathbf t^*,  T,  T^*)=(\L_1,\L_2,\l_1,\l_2,   t_1,t_2,  t_3, t_1^*,  t_2^*, t_3^*, T,  T^*)\eeqa
}
and  define via the formulae
\beqa{PR}
\arr{
 \L_1=\L_1\\
 \L_2=\L_2\\
 t_1=-\ii\sqrt{\L_1-\G_1}\,e^{{\rm i}(-\g_1+\g+\zeta)}\\
t_2=\sqrt{\L_2-\G_2}\,e^{{\rm i}(\g_2+\g+\zeta)}\\
 t_3=-\ii\sqrt{C-\G_2+\G_1}\,e^{{\rm i}(\g+\zeta)}\\
T=\sqrt{{C} -\ZZ}\,e^{{\rm i}\zeta}
}\qquad 
\arr{
\l_2=\ell_2+\g_2+\g+\zeta\\
 \l_1=\ell_1+\g_1-\g-\zeta\\
  t_1^*=-\sqrt{\L_1-\G_1}\,e^{-{\rm i}(-\g_1+\g+\zeta)}\\
t_2^*=-{\rm i}\sqrt{\L_2-\G_2}\,e^{-{\rm i}(\g_2+\g+\zeta)}\\
t_3^*=-\sqrt{C-\G_2+\G_1}\,e^{-{\rm i}(\g+\zeta)}\\
T^*=-{\rm i}\sqrt{{C} -\ZZ}\,e^{-{\rm i}\zeta}
}
\eeqa
We also define, for later need, $\eta_1$, $\eta_2$, $p$, $\xi_1$, $\xi_2$, $q$    via\beqa{PR1}
t_2&:=&
  \frac{\eta_2-{\rm i}\xi_2}{\sqrt2}\qquad t_1:=
  \frac{\ii\eta_1-\xi_1}{\sqrt2}\qquad\ \ \ t_3:=
  \frac{\ii p-q}{\sqrt2}\qquad\ \ \ T:=
  \frac{P-{\rm i}Q}{\sqrt2}\nonumber\\\nonumber\\
  t^*_2&:=&
  \frac{\eta_2+{\rm i}\xi_2}{\sqrt2{\rm i}}\qquad  t^*_1:=
  \frac{\ii\eta_1+\xi_1}{\sqrt2{\rm i}}\qquad \
   t_3^*:=
\frac{\ii p+q}{\sqrt2{\rm i}}\qquad T^*:=
  \frac{P+{\rm i}Q}{\sqrt2{\rm i}}.  \eeqa

\noindent
  Observe that 
   \beq{singularities1}\cM_\p:=\big\{(\bm\L, \bm\l, \mathbf t, \mathbf t^*):\ (\mathbf t, \mathbf t^*)=(0,0)\big\}\eeq
corresponds to co--circular, co--planar orbits for the two planets, with the outer planet in retrograde motion.

   \noindent
   We denote as
   \beqa{3BPRPSpi}
  \cH_{\textrm{\sc rps}^\complex_\p}=-\frac{\mu_1^3 M_1^2}{2\L_1^2}-\frac{\mu_2^3 M_2^2}{2\L_2^2}+\mu f_{\textrm{\sc rps}^\complex_\p}(\bm\L,\bm\l,   \mathbf t,  \mathbf t^*)
   \eeqa
  the expression of the Hamiltonian  \equ{3BP} using the coordinates $\textrm{\sc rps}_\p^{\complex}$ in \equ{real and complex}, which, similarly to the  prograde case, $\cH_{\textrm{\sc rps}^\complex_\p}$ is independent of $(T, T^*)$. Abusively, we shall continue calling $\textrm{\sc rps}_\p^{\complex}$ the coordinates  \equ{real and complex} deprived of $(T, T^*)$.
  
  \noindent
     \noindent
We now define a domain where letting the  $\textrm{\sc rps}_\p^{\complex}$ coordinates vary.
First of all, we observe that {\sc orc} configuration can be realized only if the planetary masses are tuned with the semi--major axes. More precisely, that, if we denote as  ``2'' and ``1''  the inner\footnote{Compared to \cite{pinzari18a}, here ``2'' and ``1''  are exchanged, in order to keep uniform notations along the paper.}, outer planet; 
as $a_2$, $a_1$, the semi--major axes of their respective instantaneous orbits around the sun;  $\a_-$, $\a_+$, with $0<\a_-<\a_+<1$,  two numbers such that the semi--axes ratio $\a:=\frac{a_2}{a_1}$ verifies

 \beq{a}\a_-<\a<\a_+\ ,\eeq  
 then  the following inequality needs to be satisfied
 \beq{masses ratio}
 \frac{m_2}{m_1}\sqrt{\a_-}>1\ .\eeq

\nl
Indeed, since the motions are almost--circular, the lenghths of the angular momenta of the planets, $ {C}_1$, $ {C}_2$ are arbitrarily close to the action coordinates $\L_1$, $\L_2$  related to their semi--major axes, which in turn are related to the semi--axes and the mass ratio via
 \beqno
 1<\frac{C_2}{C_1}\sim \frac{\L_2}{\L_1}=\frac{{  \m}_2}{{  \m}_1}\sqrt{\frac{{  M}_2}{{  M}_1}}\sqrt{\a}\eeqno
 where ${  \m}_i$, ${  M}_i$ are as in \equ{massesNEW}. 
This inequality does not make conflict with \equ{a} if one assumes that 
\beqa{kpm}k_\pm:=\frac{{  \m}_2}{{  \m}_1}\sqrt{\frac{{  M}_2}{{  M}_1}\a_\pm}>1\ .\eeqa
whence the necessity of \equ{masses ratio}.

\noindent
We then fix the domain as follows. 
The coordinates $\L_1$, $\L_2$ will be taken to vary in the set
\beq{L0}\cL:=\Big\{\L=(\L_1,\L_2):\ \L_-\le \L_1\le \L_+\ ,\ k_-\L_1\le \L_2\le k_+\L_1\Big\}\eeq
with $k_\pm$ as in \equ{kpm}, and $0<\L_-<\L_+$ to be chosen later.

\nl
The coordinates $\l=(\l_1,\l_2)$ will be taken to run in the torus $\torus^2$.

\nl
As for the coordinates $(\mathbf t,\mathbf t^\star)$,
we take a domain of the form
$$\cU_{\rm s}:=\Big\{(\mathbf t, \mathbf t^\star)\in \complex^6:\quad |(\mathbf t,\mathbf t^\star)|\le \varepsilon\Big\}$$

\noindent
The domain for $\textrm{\sc rps}_\p^{\complex}$ will then be
\beqa{DS}\cD_{\rm s}=\cL\times {\mathbb T}^2\times \cU_{\rm s}\,.\eeqa

\noindent
The following statement is a more precise version of statement 1. in Theorem \ref{thm: coexistence}.
  \begin{theorem}[\cite{pinzari18a}]\label{stable tori}
There exist two numbers   $0<\varepsilon_+<\varepsilon_0$, $0<\a_+<1$, such that, for any $0<\varepsilon<\varepsilon_+$, $0<\a_-<\a_+$, $0<\L_-<\L_+$, one can find $\m_+(\varepsilon)>0$ such that, for any $0<\m<\m_+(\varepsilon)$, in the domain $\cD_{\rm s}$ there exists an invariant set ${\cal F}_{\varepsilon,\m}\subset \cD_{\rm s}$ with density going to $1$ as $\varepsilon\to 0$ which is foliated as
\beq{foliation}{\cal F}_{\varepsilon,\m}=\bigcup_{\omega}\cT_{\o,\varepsilon, \m}\eeq where $\cT_{\o,\varepsilon, \m}$ is  diffeomorphic to  $\torus^5$, where $\torus:=\real/(2\p\integer)$ is the standard, ``flat'' torus. Moreover, on  $\cT_{\o,\varepsilon, \m}$ the motions are quasi--periodic, in  {\sc orc} configuration, with suitable (``diophantine'') irrational frequencies.
\end{theorem}

  \noindent
Theorem \ref{stable tori} extends  Theorem \ref{Arnold Theorem} to {\sc orc} motions. As we briefly discuss below, even though the setting is similar, the extension is not completely trivial. Here we provide  a sketch of the proof.

  \vskip.1in
  	\noindent
   In \cite{pinzari18a} it is shown that $  \cH_{\textrm{\sc rps}^\complex_\p}$ is related to the
   Hamiltonian $  \cH_{\textrm{\sc rps}}$    in \equ{HRPS} with $n=2$ by a simple relation. 
   If, in order to avoid confusions, we equip with ``tildas'' the coordinates \equ{reg var} with $n=2$  
and denote as
\beqano
 \textrm{\sc rps}^{\complex}:=({\bm\L},\widetilde{\bm\l},   \widetilde{\mathbf t},  \widetilde{\mathbf t}^*,  \widetilde T,  \widetilde T^*)=(\L_1,\L_2,\widetilde \l_1,\widetilde \l_2,   \widetilde t_1,\widetilde t_2,  \widetilde t_3, \widetilde t_1^*,  \widetilde t_2^*, \widetilde t_3^*, \widetilde T,  \widetilde T^*)\eeqano
their complex version, defined via

   \beqa{PROld+}
\widetilde t_1&:=&
  \frac{\widetilde \eta_1-{\rm i}\widetilde \xi_1}{\sqrt2}\qquad \widetilde t_2:=
  \frac{\widetilde \eta_2-{\rm i}\widetilde \xi_2}{\sqrt2}\qquad\  \widetilde t_3:=
   \frac{\widetilde p-{\rm i}\widetilde q}{\sqrt2}\qquad\widetilde T:=
  \frac{\widetilde P-{\rm i}\widetilde Q}{\sqrt2}\nonumber\\\nonumber\\
  \widetilde t^*_1&:=&
  \frac{\widetilde \eta_1+{\rm i}\widetilde \xi_1}{\sqrt2{\rm i}}\qquad \widetilde  t^*_2:=
  \frac{\widetilde \eta_2+{\rm i}\widetilde \xi_2}{\sqrt2{\rm i}}\qquad \
   \widetilde t_3^*:=
\frac{\widetilde p+{\rm i}\widetilde q}{\sqrt2{\rm i}}\qquad\widetilde  T^*:=
  \frac{\widetilde P+{\rm i}\widetilde Q}{\sqrt2{\rm i}}  \eeqa
and, finally, introduce the   involution 
\beqa{involution}\phi_1^-\big(\L_1,\L_2,\l_1, \l_2,t,t^*, T, T^*\big):=\big(-\L_1,\L_2,-\l_1, \l_2,t,t^*, T, T^*\big)\,.\eeqa 
Then we have
\begin{proposition}[\cite{pinzari18a}]\label{signs1} $\cH_{{\rm rps}^\complex_\p}=\cH_{{\rm rps}^\complex}\circ\phi_1^-$. \end{proposition}
In particular, the coefficients of the expansion 
\beqa{quadratic retrograde expr}
f^{\rm av}_{\textrm{\sc rps}^\complex_\p}=C_0(\bm\L)+\ii \mathbf t_h\cdot \s(\bm\L) \mathbf t^*+\ii\varsigma(\bm\L)t_3 t_3^*+{\rm O}_4(\mathbf t,\mathbf t^*;\bm\L)\eeqa
of $f^{\rm av}_{\textrm{\sc rps}^\complex_\p}$ are obtained from the corresponding coefficients $\widetilde \s(\bm\L)$, $\widetilde \varsigma(\bm\L)$
computed in \cite{chierchiaPi11b} by applying the projection on $(\bm\L, \bm\l)$ of the  transformation in \equ{involution}. This immediately provides

\beqa{coefficients1}
\left\{
\begin{array}{lll}\dst\s(\L_1, \L_2)=\widetilde\s(-\L_1, \L_2)=\left(
  \begin{array}{ccc}
-\frac{\rm  s}{\L_1}&-{\rm i}\frac{\widetilde{\rm  s}}{\sqrt{\L_1\L_2}}\\
-{\rm i}\frac{\widetilde{\rm  s}}{\sqrt{\L_1\L_2}}&\frac{\rm  s}{\L_2}
  \end{array}
\right)\\\\
\dst
\varsigma(\L)=\widetilde\varsigma(-\L_1, \L_2)=
-\Big(\frac{1}{\L_2}-\frac{1}{\L_1}\Big){\rm  s} \end{array}
\right.\eeqa
with
 \beqa{coefficients2}
 {\rm  s}:=- m_1m_2\frac{\a}{2a_1}{b^{(1)}_{3/2}}(\a)\qquad \widetilde{\rm  s}:=m_1m_2\frac{\a}{2a_1}{b^{(2)}_{3/2}(\a)}\qquad \quad \alpha=\frac{a_2}{a_1}
\eeqa
where
$b^{(j)}_s(\a)$'s being the Laplace coefficients\footnote{The Laplace coefficients defined via the Fourier expansion $$\frac{1}{\big(1-2\a\cos\theta+\a^2\big)^s}=\sum_{k\in \integer}b^{(k)}_s(\a)e^{{\rm i} k\theta}\ \qquad {\rm i}:=\sqrt{(-1)}\ .$$}.
It is to be remarked, from the formulae in  \equ{coefficients1}--\equ{coefficients2} that the matrix $\sigma$ is symmetric but {\it not} real. This is a remarkable difference with the prograde case studied in \cite{fejoz04, chierchiaPi11b}, which, in particular, does not ensure ``a priori'' the reality of its eigenvalues.
However, the following turns true:
\begin{lemma}
The eigenvalues of the $(2\times 2)$ matrix $\sigma(\bm\L)$ in \equ{quadratic retrograde expr} are real. Hence, $(\mathbf t, \mathbf t^*)=(\mathbf 0, \mathbf 0)\in \real^3\times \real^3$ is an elliptic equilibrium point for $f^{\rm av}_{\textrm{\sc rps}^\complex_\p}$.
\end{lemma}

\proof  The  eigenvalues of $\s$ can be explicitly computed:
\beq{elleq1}\s_{1}, \s_2=\frac{\tr\s}{2}\pm\frac{1}{2}\sqrt{(\tr\s)^2-4\det\s}\ .\eeq
{Since $\tr\s=\Big(\frac{1}{\L_2}-\frac{1}{\L_1}\Big){\rm  s}$ is real, we} have  to check that the discriminant 
\beqano
\D:=(\tr\s)^2-4\det\s=(\frac{1}{\L_2}-\frac{1}{\L_1})^2{\rm s}^2+\frac{4}{\L_1\L_2}\big({\rm  s}^2-\widetilde{\rm  s}^2\big)
\eeqano
is positive.
Recalling that the Laplace coefficients  verify
$$b^{(j)}_s(\b)> b^{(j+1)}_s(\b)\quad \textrm{for all}\quad s>0,\quad j\in \integer,\quad 0<|\b|<1,$$
(see~Ref.\cite{fejoz04} for a proof), one has
\beq{elleq2}
{\rm  s}^2-\widetilde{\rm  s}^2=(m_1m_2\frac{\a}{a_1})^2\big((b^{(1)}_{3/2}(\a))^2-(b^{(2)}_{3/2}(\a))^2\big)> 0.\eeq
and we have the assertion. $\quad\square$

\vskip.1in
\noindent
The formulae in \equ{coefficients1}--\equ{coefficients2} show that, as in the prograde case, the eigenvalues of $\sigma(\bm\L)$ and the number $\varsigma(\bm\L)$ verify, identically
\beqa{HR}\s_1+\s_1+\varsigma\equiv0\eeqa
By analogy with the latter identity in \equ{Herman resonance}, we shall refer to \equ{HR} as {\it Herman resonance}.
The asymptotic values of  the eigenvalues $\sigma_1$, $\sigma_2$ and $\varsigma$ in the well--spaced regime \equ{L0} can be computed directly from \equ{elleq1}--\equ{elleq2}, or from the corresponding ones in \cite{fejoz04, chierchiaPi11b} applying the transformation \equ{involution}. In any case, the result is
$$\left\{
\begin{array}{lll}
\dst \s_1=+\frac{3}{4\L_1}\frac{a_2^2}{a_1^3}+{\rm O}(\frac{a_2^3}{a_1^4\L_1})\\\\
\dst\s_2=-\frac{3}{4\L_2}\frac{a^2_2}{a_1^3}+{\rm O}(\frac{a_2^3}{a_1^4\L_2})\\\\\dst\varsigma=\frac{3}{4}\frac{a^2_2}{a_1^2}\left(\frac{1}{\L_2}-\frac{1}{\L_1}\right)+{\rm O}(\frac{a_2^3}{a_1^4\L_2})
\end{array}
\right.
$$
It shows that there is no other resonance besides Herman resonance in \equ{HR},
 provided the semi--axes are well spaced. Recall the definition of $\cL$ in \equ{L0}.
\begin{lemma}\label{lem: HR}
For any $K>0$, there exist $\L_\pm$, $\alpha_\pm$ such that the triple $\Omega^\complex(\bm\L):=\big(\s_1(\bm\L), \s_2(\bm\L), \varsigma(\bm\L)\big)$ verifies
\beqa{nonresuptoHR}\Omega^\complex(\bm\L)\cdot k\ne 0\quad \forall k\in \integer^3\,,\ 0<|k|\le K\,,\ k\ne N(1,1,1)\quad \forall\ \bm\L\in \cL\eeqa
with some $N\in \integer$.
\end{lemma}

\noindent
At first sight, Lemma \ref{lem: HR} might seem an obstruction towards the construction of the Birkhoff normal form for the Hamiltonian \equ{3BPRPSpi}. However, as in the prograde case, the conservation of the angular momentum lenghth     \beqa{CRPS}C=\L_2-\L_1-{\rm i}\mathbf t\cdot\mathbf t^*\eeqa is of great help. Indeed, by the commutation of $ f_{\textrm{\sc rps}^\complex_\p}$ and $C$, it turns out that, in the Taylor expansion \equ{quadratic retrograde expr}, only monomials with literal part
${\mathbf t}^{\mathbf a}{\mathbf t^*}^{\mathbf a^*}$
verifying 
\beqa{Sum ai}\sum_{i} a_i=\sum_{i} a^*_i\eeqa
appear. In \cite{chierchiaPi11c} it is shown  that, because of \equ{Sum ai}, then \equ{nonresuptoHR} is  sufficient  for constructing a Birkhoff normal form  (i.e., Theorem \ref{planetary normal form} with $n=2$) for the  Hamiltonian \equ{3BPRPSpi}. Moreover, the torsion matrix (i.e.,  the matrix $\tau(\bm\L)$ defined via  \equ{torsion}) for this case can be computed from the analogue one from the prograde case again applying \equ{involution} to the torsion of the prograde problem.  The computation is omitted (see \cite{pinzari18a} for the details), apart for stating that it is non--singular. An application of Theorem \ref{FT} then leads to the proof of Theorem \ref{stable tori}.

\paragraph{Proof of 2.} 
As a second set of coordinates, we use the $\cP$--coordinates defined in Section \ref{The reduction of perihelia}.
In the case $n=2$, they reduce to
$$\cP=(Z, C, \bm\Theta, \bm\L, \zeta, \k_2 \bm\vartheta, \bm\ell)$$
with
$$\bm\L=(\L_1, \L_2)\,,\ \bm\Theta=(\Theta_1, \Theta_2)\,,\ \bm\ell=(\ell_1, \ell_2)\,,\ \bm\vartheta=(\vartheta_1, \vartheta_2)$$
We denote as
$$\cH_{\cP}=-\sum_{j=1}^2\frac{{\mu}^3_j{M}^2_j}{2 \L_i^2}+\mu f_{\cP} (\bm\L, \bm\Theta, \bm\ell, \bm\vartheta; C)$$
the four--degrees--of--freedom Hamiltonian \equ{3BP} written using $\cP$--coordinates, which is independent of $Z$, $\zeta$ and $\k_2$.
\\
The manifols
\beqa{equilibrium}\cD^0_{\rm u}:=\Big\{(\bm\L, \bm\Theta, \bm\ell, \bm\vartheta; C):\quad (\Theta_2,\vartheta_2)=(0,0)\Big\}\eeqa
corresponds to retrograde motions. It is invariant as $f_{\cP} $ has an equilibrium on it and includes, in particular, the manifold $\cM_\p$ in  \equ{singularities1}.

\noindent
We establish  a suitable domain (including $\cD^0_{\rm u}$) for the coordinates $\cP$   where $\cH_{\cP}$ is regular. We  check below that the following domain is suited to the scope: 
\beqa{DC}
\cD_{\cP}({C} )&:=&\Big\{(\bm\L,  \Theta_1)\in \cA({C} )\Big\}\times\Big\{ (\bm\ell, \vartheta_1)\in {\mathbb T}^3\Big\}\times\Big\{ (\Theta_2,\vartheta_2)\in \cB(\Theta_1,{C} )\Big\}
\eeqa
where
\beqa{B}
\cA({C} )&:=&\Big\{(\L_1, \L_2, \Theta_1): \ (\L_1,\L_2)\in \cL({C} ), \Theta_1\in \cG(\L_1,\L_2,{C} )\Big\}\nonumber\\
\cB(\Theta_1,{C} )&:=&\Big\{(\Theta_2,\vartheta_2): \ |\Theta_2|< \frac{1}{2}\min\{{C} ,\Theta_1\}, |\vartheta_2|< \frac{\p}{2}\Big\}\nonumber\\
\cL({C} )&:=&\Big\{\bm\L: \ \bm\L\in \cL,\quad \L_2>{C} +\frac{2}{ c}\sqrt{\a_+}\L_1\Big\}\nonumber\\
\cG(\L_1,\L_2,{C} )&:=&\Big({C} _-,{C} _+\Big),\qquad {C} _-:=\frac{2}{ c}\sqrt{\a_+} \L_1\qquad {C} _+:=\min\Big\{\L_2-{C} , \L_1\Big\}.
\eeqa
with $\cL$ is as in \equ{L0}, 
 while $c$ is an arbitrarily fixed number in $(0,1)$.
We need to establish two kinds of conditions.

\subparagraph{\it a) existence of the perihelia}

We need that  the planets' eccentricities $e_1$, $e_2$ stay strictly confined in $(0,1)$. Namely, that the following inequalities are satisfied: 
\beq{C1C2}0<\Theta_1<\L_1\,,\qquad 0<C_2<\L_2\eeq
with $C_2:=|\mathbf C_2|$, $\mathbf C_2$ as in~\equ{C}. The expression of $C_2$ using $\cP$ is
$$C_2=\sqrt{C^2+\Theta_1^2-2\Theta_{2}^2+2\sqrt{(C^2-\Theta_{2}^2)(\Theta_1^2-\Theta_{2}^2)}\cos{\vartheta_{2}}}$$
We observe that $C_2$ may vanish only for
$ (\Theta_2,\vartheta_2)= (0 ,\p)$. Since we deal with the equilibrium \equ{equilibrium}, the occurrence of this equality  is  automatically excluded, limiting the values of the coordinates $(\Theta_2,\vartheta_2)$ in the set $\cB$ in~\equ{B} since in this case
\beq{C1lowerbound}C_2^2\ge \frac{3}{4}{C} ^2.\eeq
\\
Moreover, the two right inequalities in~\equ{C1C2} are satisfied taking 
\beqno\Theta_1<\min\Big\{\L_2-{C} , \L_1\Big\}={C} _+\eeqno
where we have used the triangular inequality $C_2=|\mathbf C-\mathbf C_1|\le |\mathbf C|+|\mathbf C_1|={C} +\Theta_1$.
\subparagraph{\it b) non--collision conditions} 

We have to exclude possible encounters of the planets with the sun and  each other.
Collisions of the inner planet with the sun are excluded by~\equ{B}. Indeed, using~\equ{C1lowerbound},
$$1-e_2^2=\frac{C_2^2}{\L_2^2}\ge \frac{3}{4}\frac{{C} ^2}{\L_2^2}$$
whence the minimum distance of the inner planet with the sun $a_2(1-e_2)$ is positive.
In order to avoid planetary collisions, it is typical to ensure the following inequality: 
$$a_2(1+e_2)<{c^2}a_1(1-e_1)$$
with $0<c<1$.
A sufficient condition for it is
	\beqno\Theta_1\ge \frac{2}{{c}}\sqrt{\a_+} \L_1={C} _-.\eeqno
 Indeed, if this inequality is satisfied, one has
$$a_2(1+e_2)<2a_2<\frac{a_1}{2}\frac{\Theta_1^2{c^2}}{\L_1^2}=\frac{a_1}{2}(1-e_1^2){c^2}<a_1(1-e_1){c^2}.$$

\subparagraph{The hyperbolic equilibrium \cite{pinzari18a}} 
 By the formulae \equ{HnnPeri}--\equ{ovl f} with $n=2$, the $\bm\ell$ of
 $\cH_{\cP}$  is given by
\beqano \ovl\cH_{\cP}=-\sum_{j=1}^2\frac{{\mu}^3_j{M}^2_j}{2 \L_i^2}+\mu\left(-\frac{m_{1}m_{2}}{a_{1}}+\ovl {f_{{\cal P}}^{12}}^\ppd\right)+\frac{\mu}{a_1}{\rm O}\left(\frac{a_2^2}{a_1^2}\right)
\eeqano
with
\beqano {\ovl {f_{\cP}^{12}}}^{(2)}&=&m_1 m_2 \frac{a_2^2}{4a_1^3}\frac{\L_{1}^3}{\Theta_1^5}\Big[ \frac{5}{2}(3\Theta_{2}^2-\Theta_1^2)\nonumber\\
    &-&\frac{3}{2}\frac{4\Theta_{2}^2-\Theta_1^2}{\L_{2}^2}\Big(C^2+\Theta_1^2-2\Theta_{2}^2+2\sqrt{(C^2-\Theta_{2}^2)(\Theta_1^2-\Theta_{2}^2)}\cos{\vartheta_{2}}\Big)\nonumber\\
    &+&\frac{3}{2}\frac{(\Theta_1^2-\Theta_{2}^2)(C^2-\Theta_{2}^2)}{\L_2^2}\sin^2{\vartheta_{2}}\,.\eeqano
 We shall now prove that, restricting the domain \equ{DC} a little bit, so that the manifolds \equ{equilibrium} are hyperbolic for ${\ovl {f_{\cP}^{12}}}^{(2)}$. We fix the following domain 
 \beqa{DU}\cD_{\rm u}:=\cA_{\rm u}\times \cB_{\rm u}\times\torus^3\eeqa
 with
\beqa{assumptions0}
\cA_{\rm u}({C} )&:=&\Big\{(\L_1,\L_2)\in {\cal L}_{\rm u}({C} ),\quad \Theta_1\in{\cal G}_{\rm u}(\L_1,\L_2,{C} )\Big\}\nonumber\\
\cB_{\rm u}({C} )&:=&\Big\{(\Theta_2,\vartheta_2): \ |\Theta_2|< \frac{{C} }{2}, |\vartheta_2|< \frac{\p}{2}\Big\}
\eeqa 
where
\beqa{assumptions}\cL_{\rm u}({C} )&:=&\Big\{\L=(\L_1,\L_2)\in \cL: \ 5\L_2^2{C}  -({C} +\frac{2}{ c}\sqrt{\a_+}\L_2)^2 (4 {C} +\frac{2}{ c}\sqrt{\a_+}\L_2)>0,\nonumber\\
&& \hspace*{10em} \L_1>{C} , \L_2>\max\{{C} +\frac{2}{ c}\sqrt{\a_+}\L_1, 2{C} \}\Big\}\nonumber\\
\cG_{\rm u}(\L_1,\L_2,{C} )&:=&\Big(\ovl{C} _-, \ovl{C} _+\Big)
\eeqa
where $\cL$ is as in~\equ{L0} and, if ${C} ^{\star}(\L_2,{C} )$ is the unique positive root of the cubic polynomial ${C} _2\to 5\L_2^2{C}  -({C} +{C} _2)^2 (4 {C} +{C} _2)$, then

\beq{assumptions3}\ovl{C} _-:=\max\{\frac{2}{ c}\sqrt{\a_+}\L_1, {C} \}\qquad \ovl{C} _+:=\min\{ \L_1, {C} ^{\star}\}.\eeq
Implicitly, we shall prove that
\beq{assumptions7}\ovl C_-< \ovl C_+\ .\eeq
We check that the coefficients in front of $\Theta_2^2$, $\vartheta_2^2$ in the   Taylor expansion about  $(\Theta_2, \vartheta_2)=(0, 0)$ have opposite sign in the domain \equ{DU}, so that the equilibrium manifold \equ{equilibrium} is hyperbolic. Indeed, the part of degree 2 in such expansion is
 \beqano
 m_1m_2\frac{a_2^2}{a_1^3}\frac{1}{8}\frac{\L_{1}^3}{\L_2^2\Theta_1^5} \times \Big[\frac{3 a}{C}\Theta_2^2+3C\Theta_1^2{b}
\vartheta_2^2+{\rm O}(\Theta_2^4+\vartheta_2^4)\Big]
\eeqano
where

\beq{a b} a:=5\L_2^2C -(C+\Theta_1)^2 (4C+\Theta_1)\quad {\rm and}\quad  b:=C- \Theta_1.\eeq
Both $\Theta_1\to a(\L_1,\Theta_1;{C} )$ and $\Theta_1\to b(\Theta_1;{C} )$, as functions of $\Theta_1$ decrease monotonically from a positive value (respectively, ${C} (5\L_2^2-4{C} ^2)$ and ${C} $) to $-\infty$ as $\Theta_1$ increases from $\Theta_1=0$ to $\Theta_1=+\infty$. The function $a(\L_1,\Theta_1;{C} )$ changes its sign for $\Theta_1$ equal to a suitable unique positive value ${C} ^{\star}(\L_2,{C} )$, while  $b(\Theta_1;{C} )$ does it for $\Theta_1={C} $. We note that (i) inequality ${C} <\min\{{C} _+, {C} ^{\star}\}$ follows immediately from the assumptions~\equ{assumptions} (in particular, the two last ones) and (ii), more generally, that
 ${C} ^{\star}\leq{C} $ is equivalent to $\L_2\leq2{C} $. Since, for our purposes, we have to exclude ${C} ^{\star}={C} $ (otherwise, $a(\L_1,\Theta_1;{C} )$ and $b(\Theta_1;{C} )$ would be simultaneously positive and simultaneously negative, and no hyperbolicity would be possible), we distinguish two cases. 
 
 \begin{itemize}
 \item[(a)] ${C} >\frac{2}{ c}\sqrt{\a_+} \L_1$ and ${C} +\frac{2}{ c}\sqrt{\a_+} \L_1<\L_2<2{C} $. In this case ${C} ^\star<{C} $. We show that no such $\cG_{\rm u}$ can exist in this case. In fact, since ${C} ^{\star}<{C} $, in order that the interval $({C} ^{\star}, {C} )$ and the set  ${\cal G}$ have a non-empty intersection,  one should have, necessarily, ${C} _+=\sup{\cal G}>{C} ^{\star}$, hence, in particular,  $\L_2-{C} >{C} ^{\star}$. Using the definition of ${C} ^{\star}$, this would imply $\L_2>2{C} $, which is a contradiction.

 \item[(b)] $\L_2>\max\{2{C} , {C} +\frac{2}{ c}\sqrt{\a_+} \L_1\}$. In this case ${C} <{C} ^{\star}<\L_2-{C} $. In order that the interval $({C} , {C} ^{\star})$ and the set  ${\cal G}$ have a non-empty intersection, we need
\beq{conditions}{C} _-<{C} ^{\star}\qquad {\rm and}\qquad {C} _+>{C} \eeq and such intersection will be given by the interval $\cG_{\rm u}$ as in~\equ{assumptions}. Note that the definition of $\ovl{C} _+$ does not include $\L_2-{C} $ in the brackets because, as noted, ${C} ^{\star}<\L_2-{C} $. 
But~\equ{conditions} are equivalent to
\equ{assumptions}. 
\end{itemize}

\paragraph{Proof of 3.}
Here we prove that
\begin{theorem}\label{coexistence of tori}
Let $\a_+<\frac{1}{16}$.
There exist  universal numbers $1<\underline k<\ovl k$ such that, if
$$ \a_-<\frac{\underline k^2}{\ovl k^2}\a_+\ ,\quad\frac{ \ovl k}{\sqrt{\a_+}}<\frac{{  \m}_2}{{  \m}_1}\sqrt{\frac{{  M}_2}{{  M}_1}}<\frac{\underline k}{\sqrt{\a_-}}$$
then $\cD_{\rm s}\cap \cD_{\rm u}^0$ is non--empty. The following values work:
\beq{ass1} \underline k=\frac{1}{4}\sqrt{\frac{3}{10}(69+11\sqrt{33})}\sim 1.57\ ,\quad \ovl k=2\ .\eeq 
\end{theorem}
\proof The sets $\cD_{\rm s}$ in \equ{DS} and $\cD^0_{\rm u}$ in \equ{equilibrium} are expressed with different sets of coordinates. To prove that $\cD_{\rm s}$  and $\cD^0_{\rm u}$ have a non--empty intersection, we need to use the same set for both. We choose to use the coordinates $\cP$, so we rewrite $\cD_{\rm s}$ in terms of $\cP$.

\begin{figure}[htp]
\centering{
\includegraphics
{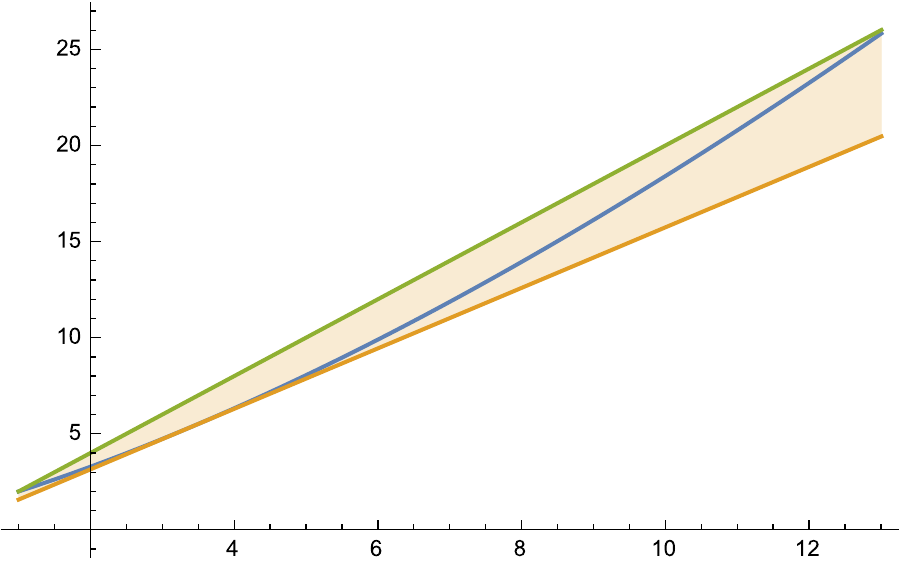}}
\caption{The blue curve is ${\cal C}$; the orange line has slope $\underline k$, the green one has slope $\ovl k$ (\textsc{Mathematica}).\label{coexistence1}
}
\end{figure}

\begin{figure}[htp]
\centering{
\includegraphics
{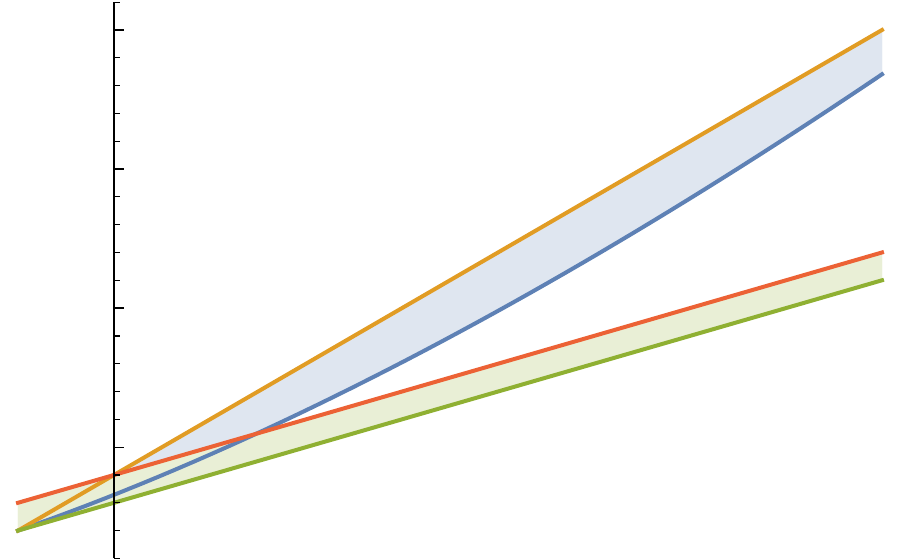}}
\caption{The blue strip corresponds to the set $\cL_1$, the green one to $\cL_2$ (\textsc{Mathematica}).  \label{coexistence}
}
\end{figure}

\begin{figure}[htp]
\centering{
\includegraphics
{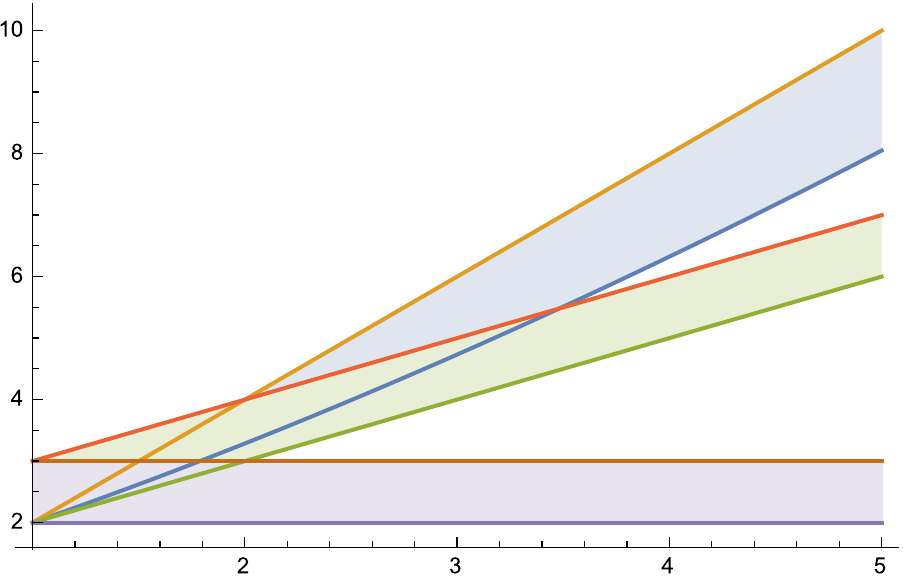}}
\caption{$\cL_1$: the blue region; $\cL_2$: the green region; $\cL_3$: the  violet region (\textsc{Mathematica}). \label{strip}
}
\end{figure}

\noindent
Using $\cP$, the set $\cD_{\rm s}$ becomes (at the expenses of diminishing $\varepsilon$, if necessary)
$$\cD_{\rm s}=\cA_{\rm s}\times \cB_{\rm s}\times\torus^3$$
where, if
\beq{assumptions1}\cL_{\rm s}({C} ):=\Big\{\L=(\L_1,\L_2)\in \cL_0:\ |\L_2-\L_1-{C} |<\varepsilon\Big\}\ ,\qquad {\cal G}_{\rm s}(\L_1):=\Big\{\Theta_1:\ 0<\L_1-\Theta_1<\varepsilon\Big\}\,,\ \eeq
then
\beq{assumptions22}\cA_{\rm s}:=\Big\{(\L_1,\L_2,\Theta_1):\ (\L_1,\L_2)\in \cL_{\rm s}\ ,\ \Theta_1\in \cG_{\rm s}(\L_1)\Big\}\,,\ \cB_{\rm s}:=\Big\{(\Theta_2,\vartheta_2):\ |(\Theta_2,\vartheta_2)|<\varepsilon\Big\}\ .\eeq
All we have to do is to check that the intersection
$\cA_{\rm s}\cap\cA_{\rm u}$
is non--empty. \\
 Recalling the definition of $\cA_{\rm u}$ in \equ{assumptions0}--\equ{assumptions}
and the definition of $\cA_{\rm s}$ in \equ{assumptions1}--\equ{assumptions22}, asserting that $\cA_{\rm s}\cap\cA_{\rm u}\ne \emptyset$ 
is equivalent to asserting that
$${\cal L}_{\rm s}( {C})\cap{\cal L}_{\rm u}( {C})\ne \emptyset$$
and
$${\cal G}_{\rm s}(\L_1)\cap{\cal G}_{\rm u}(\L_1,\L_2, {C})\ne \emptyset\quad \forall\ (\L_1,\L_2)\in {\cal L}_{\rm s}( {C})\cap{\cal L}_{\rm u}( {C})\ .$$
It will be enough to check that
\beq{assumptions5}{\cal L}_{\rm s}( {C})\cap{\cal L}_{\rm u}( {C})\cap{\cal L}_{{\rm su}}( {C})\ne \emptyset\eeq
and
\beq{assumptions4}{\cal G}_{\rm s}(\L_1)\cap{\cal G}_{\rm u}(\L_1,\L_2, {C})\ne \emptyset\quad \forall\ (\L_1,\L_2)\in {\cal L}_{\rm s}( {C})\cap{\cal L}_{\rm u}( {C})\cap{\cal L}_{{\rm su}}( {C})\ ,\eeq
where, if $\ovl {C}_\pm$ are as in \equ{assumptions3}, $\cL_{{\rm su}}$ is defined as
\beq{LSU}\cL_{{\rm su}}:=\Big\{(\L_1,\L_2):\ \ovl {C}_+=\L_1\Big\}\ .\eeq
Note that \equ{assumptions4} is certainly satisfied provided \equ{assumptions5} is, since in fact, for $(\L_1,\L_2)\in {\cal L}_{\rm s}( {C})\cap{\cal L}_{\rm u}( {C})\cap{\cal L}_{{\rm su}}( {C})$,
$${\cal G}_{\rm s}(\L_1)\cap{\cal G}_{\rm u}(\L_1,\L_2, {C})=\Big\{ \Theta_1:\ \max\{\ovl {C}_-,\L_1-\varepsilon\}< \Theta_1<\L_1\Big\}$$
which is well--defined by \equ{assumptions3}--\equ{assumptions7}.

\nl
On the other hand, in view of the definition of $\ovl {C}_+$ in \equ{assumptions3}, and of $ {C}^\star$ a few lines above, $\cL_{{\rm su}}$ in \equ{LSU} is equivalently defined as

\beq{L2}\cL_{{\rm su}}=\Big\{(\L_1,\L_2):\  5\L_2^2 {C} -( {C}+\L_1)^2 (4  {C}+\L_1)>0\Big\}\ .\eeq
Therefore, in view of this definition and the definitions of $\cL_{\rm s}$, $\cL_{\rm u}$ in \equ{assumptions} and \equ{assumptions1}, one sees that the set on the left hand side in \equ{assumptions5} is determined by inequalities

\beqa{allinequalities}
&&\L_-<\L_1<\L_+\nonumber\\
&&k_-\L_1\le  \L_2\le k_+\L_1\nonumber\\
&&5\L_2^2 {C} -( {C}+2\sqrt{\a_+}\L_2)^2 (4  {C}+2\sqrt{\a_+}\L_2)>0\nonumber\\
&& \L_1> {C}\nonumber\\
&&\L_2>\max\{ {C}+2\sqrt{\a_+}\L_1,\ 2 {C}\}\nonumber\\
&& |\L_2-\L_1- {C}|<\varepsilon\nonumber\\
&&5\L_2^2 {C} -( {C}+\L_1)^2 (4  {C}+\L_1)>0\eeqa

\nl
We  observe that no phase point\footnote{Inequalities $\L_1< {C}_\star$ (which is equivalent to \equ{L2}) and $ {C}_\star<\L_2- {C}$ (which is equivalent to $ \Theta_1> {C}$, in turn implied by the definition of $\cG_{{\rm su}}$ above) imply $\L_2-\L_1- {C}>0$.  
}  $(\L_1,\L_2)$ with $ \L_2-\L_1- {C}<0$ will ever satisfy \equ{allinequalities}, and that inequality $\L_2>2 {C}$ is implied by $\L_1> {C}$ and \equ{L2}. 
Then, we divide such inequalities in three groups, so as to rewrite the set \equ{assumptions5} as the intersection of the sets

\beqano
\widehat\cL_1&:=&\Big\{(\L_1,\L_2):\ \L_-<\L_1<\L_+,\ \L_1> {C} ,\ \L_2>2 {C}\ ,\nonumber\\
&&  \max\{k_-\L_1,\ ( {C}+\L_1)\sqrt{\frac{4  {C}+\L_1}{5 {C}} }\}<\L_2\le k_+\L_1\Big\}\nonumber\\
\widehat\cL_2&:=&\Big\{(\L_1,\L_2):\ 0<\L_2-\L_1- {C}<\varepsilon,\ \L_2> {C}+2\sqrt{\a_+}\L_1\ ,\ \L_1> {C}\Big\}\nonumber\\
\widehat\cL_3&:=&\Big\{(\L_1,\L_2):\ 5\L_2^2 {C} -( {C}+2\sqrt{\a_+}\L_2)^2 (4  {C}+2\sqrt{\a_+}\L_2)>0\ ,\ \L_2>2 {C}\Big\}
\eeqano

\nl
We now aim to choose the parameters $\L_\pm$, $k_\pm$ and $\a_+$ so as to find a non--empty intersection of the sets a above.

\nl
Let us denote as ${\cal C}$ the curve, in the $(\L_2,\L_1)$--plane, having equation
\beq{curve}{\cal C}:\qquad \L_2=( {C}+\L_1)\sqrt{\frac{4  {C}+\L_1}{5 {C}} }\eeq
 
\nl
Let
$$\L_1= k\L_2$$
be any straight line through the origin.
The  straight line
intersecting ${\cal C}$ into the point $(\underline{\L_1},\underline{\L_2})=( {C},2 {C})$ has  $\ovl k=2$, and intersects this curve, also in the higher point
$$(\ovl{\L_1},\ovl{\L_2})=(\frac{1}{2}(13+\sqrt{185}), (13+\sqrt{185})) {C}\ .$$
Any other line with $k>\ovl k$ has a lower intersection   $(\underline{\L_1}',\underline{\L_2}')$, with $\underline{\L_1}'< {C}$ and $\underline{\L_2}'<2 {C}$ and a higher intersection $(\ovl{\L_1}',\ovl{\L_2}')$   with $\ovl{\L_1}'>\ovl{\L_1}$ and $\ovl{\L_2}'>\ovl{\L_2}$.

\nl
The last straight line, in the plane $(\L_1,\L_2)$,   through the origin 
intersecting ${\cal C}$ is the tangent line, and  it is easy to compute (see below) that such a  tangent line has slope has slope $\underline k$ as in \equ{ass1} (Figure~\ref{coexistence1}).
We then conclude that, as soon as we choose $k_-<\underline k$, $k_+>\ovl k$, $\L_-<\underline{\L_1}$, $\L_+>\ovl{\L_1}$, we have   the inclusion
$$\widehat\cL_1\supset\cL_1:=\Big\{(\L_1,\L_2):\ ( {C}+\L_1)\sqrt{\frac{4  {C}+\L_1}{5 {C}} }\}<\L_2\le 2\L_1\Big\}\ .$$

\nl
Let us now turn to $\widehat\cL_2$. Since we are assuming $\a_+<\frac{1}{16}$, we conclude that the strip
$$\cL_2:=\Big\{(\L_1,\L_2):\ 0<\L_2-\L_1- {C}<\varepsilon\ ,\ \L_1> {C}\Big\}$$
is all included in the region
$$\widetilde\cL_2=\Big\{(\L_1,\L_2):\ \L_2> {C}+2\sqrt{\a_+}\L_1\ ,\ \L_1> {C}\Big\}$$
and this allows to conclude
$$\widehat\cL_2=\cL_2\cap\widetilde\cL_2=\cL_2\ .$$

\nl
Since the sets  $\cL_1$ and $\cL_2$ have a non--empty intersection, independently of $\a_+$  (see Figure~\ref{coexistence}), 
a fortiori, $\widehat\cL_1$ and $\widehat\cL_2$ have one:
$$\widehat\cL_1\cap\widehat\cL_2\supset\cL_1\cap\cL_2\ne \emptyset\ .$$
Observe, in particular, that $\cL_1\cap\cL_2$ (hence, $\widehat\cL_1\cap\widehat\cL_2$) has non--empty intersection with any strip
$\real\times \Big[2 {C},y\Big]$, with $y>2 {C}$ (see Figure~\ref{strip}).

\nl
On the other hand, it is immediate to check that $\widehat\cL_3$ includes the horizontal strip
$$\cL_3:=\Big\{(\L_1,\L_2):\ 2 {C}<\L_2<\frac{ {C}}{2\sqrt{\a_+}}\ ,\ \L_1\in \real\Big\}\qquad 0<\a_+<\frac{1}{16}$$
and so we conclude
$${\cal L}_{\rm s}( {C})\cap{\cal L}_{\rm u}( {C})\cap{\cal L}_{{\rm su}}( {C})=\widehat\cL_1\cap\widehat\cL_2\cap\widehat\cL_3\supset \cL_1\cap\cL_2\cap\cL_3\ne \emptyset$$
In order to complete the proof, it remains to prove that the tangent straight line to $\cC$ through the origin has slope has slope $\underline k$ as in \equ{ass1}.\\
 We switch to the homogenized variables
 $$x:=\frac{\L_1}{{C}}\qquad y=\frac{\L_2}{{C}}$$
 so that the curve $\cC$ in \equ{curve} becomes
 $$\widehat\cC:\qquad y=(1+x)\sqrt{\frac{4+x}{5}}\ .$$
 We look for a straight line through the origin $y=\underline k x$ with $\underline k>0$ which is tangent to $\widehat\cC$ at some point $(a,b)$, with $a>0$. 
 
 \nl
 The intersections between $\widehat\cC$ and any straight line through the origin $y= k x$ are ruled by a complete cubic equation, given by
\beq{1st}x^3+(6-5 k^2)x^2+9x+4=0\ .\eeq
In order that such an equation has a double solution $x=a$ for $k=\underline k$, one needs that, when $k=\underline k$, it can  factorized as
\beq{2nd}(x-a)^2(x-c)=0
\eeq
 Therefore, equating the respective coefficients of \equ{1st} and \equ{2nd} one finds the equations 

 $$\arr{\dst -(c+2a)=6-5\underline k^2\\\\
 \dst 2ac+ a^2=9\\\\
 \dst -a^2c=4
 }$$
Two last equations, allow to eliminate $b$ so as to obtain the equation for $a$
$$a^3-9 a-8=0$$
which has the following three roots:
$$a_0=-1\ ,\qquad a_{\pm}=\frac{1\pm\sqrt{33}}{2}\ .$$
The only admissible (positive) value is then
$$a=a_+=\frac{1+\sqrt{33}}{2}$$
and it provides the values
$$c=\frac{-17+\sqrt{33}}{32}\ ,\qquad \underline k=\frac{1}{4}\sqrt{\frac{3}{10}(69+11\sqrt{33})}\ .\quad \square$$

\appendix

\section{The $m_0$--centric reduction}\label{appendix}
The Hamiltonian of $1+n$ masses $m_0$, $\ldots$, $m_n$ interacting through gravity is
\beqa{manybody}\cH=\sum_{i=0}^{n}\frac{|\mathbf u_i|^2}{2m_i}-\sum_{0\le i<j\le n}\frac{m_i m_j}{|\mathbf v_i-\mathbf v_j|}\,.\eeqa
We switch from the position coordinates $\mathbf v_i$ to new new ones, denoted $\mathbf x_i$, where $\bx_0$ is the coordinate of $m_0$, while $\bx_i$ is the coordinate of $m_i$ relatively to $m_0$. The change is
\beqa{xi}
\mathbf v_i=\left\{\begin{array}{lll}
\mathbf x_0\quad &i=0\\
\mathbf x_i+\mathbf x_0\quad &i=1\,,\ldots\,,n\,.
\end{array}
\right.
\eeqa
As the change does not involve the $\mathbf u_i$'s, the  coordinates $\by_i$ conjugated to $\bm x_i$ may be computed imposing the conservation of the standard 1--form
$$\L=\sum_{i=0}^n\bm y_i\cdot d\bm x_i=\sum_{i=0}^n\mathbf u_i\cdot d\mathbf v_i\,.$$
We find
\beqano
\sum_{i=0}^n\mathbf u_i\cdot\mathbf v_i&=&\mathbf u_0\cdot\mathbf x_0+\sum_{i=1}^n\mathbf u_i\cdot(\mathbf x_i+\bx_0)\nonumber\\
&=&\left(\sum_{i=0}^n\mathbf u_i\right)\cdot \bx_0+\sum_{i=1}^n\mathbf u_i\cdot\mathbf x_i\,.
\eeqano
So we identify
$$\mathbf y_i=\left\{\begin{array}{lll}
\dst\sum_{i=0}^n\mathbf u_i\quad &i=0\\\\
\dst\mathbf u_i\quad &i=1\,,\ldots\,,n\,.
\end{array}
\right.$$
We recognize that $\by_0$ is the total linear momentum, which keeps constant along the motions of $\cH$, as $\cH$ is translation--invariant. Fixing a reference frame moving with the  the centre of mass of $m_0$, $\ldots$, $m_n$, we have $\by_0=0$ and hence
\beqa{yiNEW}\mathbf u_i=\left\{\begin{array}{lll}
\dst-\sum_{i=1}^n\mathbf y_i\quad &i=0\\\\
\dst\mathbf y_i\quad &i=1\,,\ldots\,,n\,.
\end{array}
\right.\eeqa
Replacing \equ{yiNEW} and \equ{xi} into \equ{manybody} we arrive at \equ{Helio}, with $\mu_i$, $M_i$ as in \equ{masses}.

\newpage
\addcontentsline{toc}{section}{References}
\def\cprime{$'$} \def\cprime{$'$}

\end{document}